\documentclass{amsart}
\usepackage{graphicx}
\usepackage{amssymb}
\usepackage[all]{xy}
\usepackage{kmath}
\usepackage[T1]{fontenc}
\usepackage{hyperref}
\usepackage[parfill]{parskip}
\usepackage{float} 

\usepackage{tikz}

\usepackage{tkz-euclide}
\usetikzlibrary{decorations.markings,arrows}
\tikzset{
    arrowMe/.style={
        postaction=decorate,
        decoration={
            markings,
            mark=at position .5 with {\arrow[thick]{#1}}
        }
    }
}

\usetikzlibrary{backgrounds,calc}
\tikzset{point/.style = {fill=black,circle,inner sep=0.7pt}}
\pgfarrowsdeclaredouble{<<s}{>>s}{stealth}{stealth}%   double stealth
\pgfarrowsdeclaretriple{<<<s}{>>>s}{stealth}{stealth}% triple stealth
\tikzstyle{vertex}=[circle,fill=black!25,minimum size=20pt,inner sep=0pt]
\tikzstyle{edge} = [draw]
\tikzstyle{weight} = [font=\tiny]
\newtheorem{thm}{Theorem}[section]
\newtheorem{cor}[thm]{Corollary}
\newtheorem{lem}[thm]{Lemma}
\newtheorem{prop}[thm]{Proposition}
\theoremstyle{definition}
\newtheorem{defn}[thm]{Definition}
\theoremstyle{remark}

\numberwithin{equation}{section}
\makeatletter
\def\subsection{\@startsection{subsection}{3}%
  \z@{.5\linespacing\@plus.7\linespacing}{.001\linespacing}%
  {\normalfont\itshape}}
\makeatother
\makeatletter
\def\subsubsection{\@startsection{subsubsection}{3}%
  \z@{.5\linespacing\@plus.7\linespacing}{.001\linespacing}%
  {\normalfont\itshape}}
\makeatother
\begin{document}
\title{Numerical properties of Koszul connections}
\author{Michel NGUIFFO BOYOM}%
\address{IMAG Alexander Grothendieck Research Institute UMR CNRS 5149 University of Montpellier }%
\email{boyom@math.univ-montp2.fr}%
%\textbf{Michel Nguiffo Boyom} \\
%\textbf{IMAG, Alexander Grotehdieck Institute. Montpellier France.}
%\thanks{.TBW}
%\commby{}%
%----------------------------------------------------------------

\maketitle
%----------------------------------------------------------------
%\section{Some numerical properties of Koszul connections}
%\abstract 
%\textbf{R\'esum\'e}
\newpage
\tableofcontents
\newpage

%\begin{abstract} 
%{\centering\Large
\section*{\bfseries Abstract}
We use the notation EX(S>M), EXF(S>M) and DL(S>M), where M is a smooth manifold and S is a geometric structure. EX(S>M) is the question whether S exists in M. EXF(S>M) is the question whether M admits S-foliations. DL(S>M) is the search of an invariant measuring how M is far from admitting S. For many major geometric structures, those questions are widly open. In this paper, we address EX(S>M), EXF(S>M) and DL(S>M) for affine structure and symplectic structure, left invariant affine structure, left invariant symplectic structure and bi-invariant riemanniann structure in Lie groups
%\end{abstract} 

\section{\bfseries Prologue}
%\begin{prologue} 

\begin{enumerate}
\item 
\end{enumerate}
A gauge structure is a pair $(M,\nabla)$ where $\nabla$ is a Koszul connection in $M$. The Lie algebra of infinitesimal gauge transformations of $M$ is the real vector space of sections of the vector bundle $Hom(TM,TM)$. it is denoted by $g(M)$. The gauge group $\mathit{G}(M)$ is the open subset of inversible element s of $g(M)$. In this work we use a pair of gauge structures $\left\{(M,\nabla),(M,\nabla^*)\right\}$ for defining three fundamental differential equations $FE(\nabla\nabla*)$, $FE^*(\nabla)$, $FE^{**}(\nabla).$ This work is mainly devoted to the global analysis of the pair $\left\{FE(\nabla\nabla^*), FE^*(\nabla)\right\}$ and their impacts on notable research topics: (a) The locally flat geometry; (b) the Symplectic geometry in manifolds and in Lie groups; (c) the bi-invariant Riemannian geometry in Lie groups; (d) the information geometry. We have been motivated by two difficult fundamental open problems in the global differential geometry, in the global differential topology and in their applications in the information geometry. Consider a finite dimensional manifold $M$ and a geometric structure $\mathcal{S}$. \textbf{$EX(\mathcal{S})$}: Arises the question whether $\mathcal{S}$ exists in $M$. $EXF(\mathcal{S})$: Arises the question whether $M$ admits non trivial $\mathcal{S}$-foliations, viz foliations whose leaves carry the structure $\mathcal{S}$. If $\mathcal{S}$ is the Riemannian structure then  solutions to $EX(\mathcal{S})$ exists in every finite dimensional manifold $M$. In the subcategory of positive definite Riemannian geometry the problem $EXF(\mathcal{S})$ is but a problem in the differential topology, viz the existence of non discrete foliations in $M$. However complications arise in the case of positive signatures $(p,q)$. The case of Riemannian geometry is rather rare. Up to now there does not exist any criteron for deciding whether a finite dimensional manifold admits symplectic structures. Thirty years ago Alexander K. Guts raised the question whether a Riemannian manifold $(M,g)$ admits Hessian structures. Independently, motivated by the complexity of statistical models S-I Amari raised the same problem. The main purposes of this paper are relatioships between the global analysis of the pair $\left\{(FE(\nabla\nabla^*), FE^*(\nabla)\right\}$ and solutions to the pair $\left\{ EX(\mathcal{S}),EXF(\mathcal{S})\right\}$. Other motivations of this research paper came from the Needs of Mathematic Structures in Complex Systems (CS) whose statistical models are domains of the information geometry, See the EE document \textbf{ReportonMathematicsfor DigitalScience: Opportunities and Challenges at the Interface of Big Data, High-Performance Computing and Mathematics}. We focus on geometric structures which significantly impact some outstanding topics in the fundamental mathematic and their applications. Instances are (1) the geometry of Koszul and its links with the geometric theory of Heat. (2) The Hessian geometry and its links with the information geometry. (3) The symplectic geometry and its links with the thermodynamics, integrable systems and the information geometry. (4) The bi-invariant Riemannian geometry in Lie groups and its links with the theory of probabilities in finite dimensional Lie groups. (5) The left invariant symplectic geometry in finite dimensional Lie groups and its links with locally flat geometry in finite dimensional Lie groups. To emphasize those choice we go to recall some still open problems $EX(\mathcal{S})$ and $EXF(\mathcal{S})$. \textbf{$EX(\mathcal{S}:LF(M))$ is the existence of locally flat structures in a differentiable manifold $M$}: J-L. Koszul, E.B. Vinberg, Y. Matsushima, J. Milnor, J. Smilie and others. \textbf{$EX(\mathcal{S}:Symp(M))$ is the existence of symplectic structures in a differentiable manifold $M$}: J-M. Souriau, B. Kostant, S. Stenberg, V. Guillemin, A. Lichnerowicz, and others. \textbf{$EX(\mathcal{S}:Hes(M))$ is the existence of Hessian structures in a differentiable manifold $M$}: J-L. Koszul, E.B. Vinberg, H. Shima, S-I. Amari, J. Armstrog-Amari and others. \textbf{$EX(\mathcal{S}: Hes(M,g))$ is the existence of Hessian structures in a Riemannian manifold $(M,g)$}: S-I. Amari, A.K. Guts, J.Armstrong-Amari and others. \textbf{$EX(\mathcal{S}:Hes(M,\nabla))$ is the existence of Hessian structures in a locally flat manifold $(M,\nabla)$}: J-L. Koszul, J. Vey and others. \textbf{$EX(\mathcal{S}:lf(G))$ is the existence of left invariant locally flat structures in a finite dimensional Lie group $G$}: L. Aulander, J. Milnor, J-L; Koszul, A. Nijenhuis and others.\\ \textbf{$EX(\mathcal{S}:symp(G))$ is the existence of left invariant symplectic structures in a finite dimensional Lie group $G$}: S. Stenberg, J-M. Souriau, J-L. Koszul and others.\\ \textbf{$EX(\mathcal{S}:br(G))$ is the existence of bi-invariant Riemannian metrics in a finite dimensional Lie group $G$}: E. Cartan and many others. The open problems we just listed have their avatars in the differential topology, namely the problems $EXF(\mathcal{S}: LF(M))$, $EXF(\mathcal{S}:Symp(M))$, $EXF(\mathcal{S}:Hes(M))$, $EXF(\mathcal{S}:Hes(M,g))$, $EXF(\mathcal{S}:lf(G)$. Another challenge to face is the moduli space of finite dimensional locally flat manifolds. We focus on the subcategory of geodesically complete locally flat manifolds. For those purpose we introduce the Canonical linear representation of the fundamental group of a gauge manfold. Permanently we are motivated to search for characteristic obstructions. Another approach to both $EX(\mathcal{S})$ and $EXF(\mathcal{S})$ is the combinatorial analysis. The aim is to answer the distancelike question \textbf{$DL(\mathcal{S}:M)$: How far from admitting a structure $\mathcal{S}$ is a manifold $M$?}. For many geometric structures we bring complete solutions to the problems $EX(\mathcal{S})$, $EXF(\mathcal{S})$ and $DL(\mathcal{S})$. We also implement the fundamental equation $FE^*(\nabla)$ in revisiting the fundamental conjecture of Gindikin-Piatecci-Shapiro-Vinberg. In view of the diversity of geometric structures we are concerned with there was no hint of thinking that all of the problems $EX(\mathcal{S})$, $EXF(\mathcal{S})$ and $DL(\mathcal{S})$ might be linked with only two differential equations $\clubsuit$ \\ 
\textbf{R\'esum\'e}. Une structure de jauge est une paire $(M,\nabla)$ o\`u $\nabla$ est une connexion de Koszul dans une vari\'et\'e $M$. L'alg\`ebre de Lie des transformations de jauge infinit\'esimales est l'espace vectoriel des sections du fibr\'e vctoriel $Hom_\mathbb{R}(TM,TM)$. Le groupe de jauge $\mathit{G}(M)$ est l'ouvert des sections inversibles de $Hom(TM,TM)$. Nous utilisons une paire de structures de jauge $\left\{((M,\nabla),(M,\nabla^*)\right\}$ pour introduire deux op\'erateurs differentiels fondamentales de second ordre. Nous en d\'eduisons deux syst\`emes d'\'equations aux d\'eriv\'ees partielles not\'ees $FE^*(\nabla\nabla^*)$ et $FE^*(\nabla)$. Ce travail est consacr\'e aux impacts de l'analyse globale de ces deux \'equations sur deux problems majeurs en analyse globale sur les vari\'et\'es diff\'erentiables. Les voici. On fixe une structure g\'eom\'etrique $\mathcal{S}$ et une vari\'et\'e diff\'erentiable $M$. Le probl\`eme $EX(\mathcal{S})$ est de savoir si la structure $\mathcal{S}$ exite dans $M$. Le probl\`eme $EXF(\mathcal{S})$ est de savoir si $M$ porte des feuilletages non discret dont les feuilles sont (diff'erentiablement) des $\mathcal{S}$-vari\'et\'es. Un cas exemplaire est celui de structure de vari\'et\'e Riemannienne. Le  probl\`eme $EX(\mathcal{S})$ est r\'esolu, la structure Riemannienne existe dans toute variété diff\'erentiable . Dans la cat\'egorie des vari\'et\'es Riemanniennes positives $EXF(\mathcal{S})$ est un probl\`eme de topologie diff\'erentielle: l'existence de feuilletage non discret dans $M$. Le cas de signature positive est moins docile. Eu \'egard au probl\`eme $EX(\mathcal{S})$ le cas de la structure Riemannienne est singulier. Pour des nombreuses structures g\'eom\'etriques importantes les probl\`emes $\mathcal{S}$ et $EXF(\mathcal{S})$ sont ouverts. A titre d'illustration citons (a) la structure symplectique dans une vari\'et\'e $M$, (b) la structure symplectique invariante \`a gauche dans un groupe de Lie $G$, (c) la structure Riemannienne bi-invariante dans un groupe de Lie $G$. A ce jour aucune obstruction caract\'eristique n'en est connue. Ce travail est consacr\'e aux impacts de l'ananlyse globale des \'equations $FE^*(\nabla\nabla^*)$ et $FE^*(\nabla)$ sur ces probl\`emees $EX(\mathcal{S})$ et $EXF(\mathcal{S})$. Outre ces d\'efi purement intelectuels que posent ces deux probl\`emes il y a des motivations pratiques n\'ees de besoin de Structures G\'eom\'etriques dans les Syst\`ems Complexes (,Big Data) dont l'\'etude des  mod\`eles statistiques est un des objets de la G\'eom\'etrie et de la Topologie de l'Information, cf document EE \textbf{RepotonMathematicsfor DigitalScience: Opportunities and Challenges at the Interfaces of Big Data, High Performence Computing and Mathematics}. Nous avons limit\'e l'ambition \`a une liste non exhaustive des structures dont l'importance r\'epose sur deux qualit\'es. (1) La premi\`ere qualit\'e est la richesse en probl\`emes g\'eom\'etriques internes non r\'esolus. (2) La seconde qualit\'e est la f\'econdit\'e de leurs impacts sur d'autres domaines de recherche fondamentale et de recherche appliqu\`ee. Pour certaines structures $\mathcal{S}$ nous donnons des solutions compl\`etes des probl\`emes $EX(\mathcal{S})$ et $EXF(\mathcal{S})$. En voici une liste non exhaustive. \\
$EX(\mathcal{S}:(LF(M))$: Existence de structures localement plates (ou affinement plates) dans $M$.\\
$EX(\mathcal{S}:(SY(M))$: Existence de structures symplectiques dans $M$.\\
$EX(\mathcal{S}:(He(M))$: Existence de structures Hessiennes dans $M$.
Des d\'efis subsidiaires.\\
$EX(\mathcal{S}:(He(M,g))$: Existence de structures Hessiennes dans une vari\'et\'e Riemannienne $(M,g)$. Ce probl\`eme connu de S-I Amari a \'et\'e r\'ecemment \'etudi\'e par S-I Amari et J. Armstrong. Ind\'ependamment Alexander K. Guts posa ce probl\`eme \`a l'auteur il y a une trentaine d'ann\'ees.\\
$EX(\mathcal{S}:(He(M,\nabla))$: Existence de structures Hessiennes dans une vari\'et\'e localement
plate $(M,\nabla)$. C'est un enjeu de la G\'eom\'etrie de localement plate hyperbolique. D'apr\`es Koszul une vari\'et\'e localement plate compacte est hyperbolique si et seulement si elle porte une m\'etrique Hessienne positive dont la classe de KV cohomologie est nulle. Smilie \textbf{L'analyse globale des \'equations $FE(\nabla\nabla^*)$ et $FE^*(\nabla)$ a des relais dans la g\'eom\'etrie des groupes de Lie}. Dans la cat\' egorie des groupes de Lie le probl\`eme $EX(\mathcal{S})$ concerne les structures g\'eom\'etriques qui sont soit invariantes \`a gauche soit bi-invariantes.\\
$EX(\mathcal{S}:lf(G))$: Existence de structure localement plate invariante \`a gauche dans $G$.\\
$EX(\mathcal{S}:sym(G))$: Existence des structures symplectiques invariantes \`agauche dans $G$.\\
$EX(\mathcal{S}: br(G))$: Existence de structure Riemanniennes bi-invariantes dans $G$.\\
L'analyse globale des \'equations fondamentales a conduit \`a l'introduction de deux functions $r^b$ et $s^b$. Ces fonctions sont d\'efinies $r^b$ dans l'espace des modules de structures de gauge, $s^b$ dans l'espace des transformations de gauge infinit\'esimales. Ces fonctions sont des invariant g\'eom\'etriques. Elles sont la sources des obstructions caract\'eristiques \`a la plupart des probl\`emes que nous avons `'etudi\'es. \textbf{De prime abord rien ne laisse penser que les probl\`emes ouverts mmentionn\'es ci-dessus puissent tous \^etre trait\'es par un m\^eme Invariant G\'eom\'etrique, Topologique ou Fonctionel}. Sous une autre perspective nous avons r\'eformul\'e le probl\`eme $EX(\mathcal{S})$ en termes de quasi-distance combinatoire. \textbf{$DL(\mathcal{S}(M)$: = trouver des Invariants qui EVALUENT l'\'eloignement de $M$ de la Possession de la structure $\mathcal{S}$}. Les outils pour affronter ce d\'efi sont venus de deux \'equations fondamentales. L'\'equation diff\'erentielle $FE^*(\nabla)$ est un avatar de la conjecture de Muray Gerstenhaber sur la th\'eorie de d\'eformation des vari\'et\'es localement plates. 
A l'attention des lecteurs peu habitués aux méthodes de l'analyse globale (Pseudogroupes de Lie, Equations de Lie, G\'eom\'etre de Sternberg, Kuranishi-Formalisme de Spencer), nous soulignons la disparit\'e des probl\`emes g\'eom\'etriques $EX(\mathcal{S})$ qui sont \'etudi\'es dans ce travail. $EX(\mathcal{S}: Hes(M,\nabla))$ est un SEDP homog\`ene de second ordre; $EX(\mathcal{S}: Symp(M))$ est du domaine de calcul ext\'erieur; $EX(\mathcal{S}:LF(M))$ est du domaine du calcu tensoriel. Eu \'egard \`a cette disparit\'e rien ne sugg\`ere que leurs solutions puissent 
d\'ependre des m\^emes invariants de la g\'eom\'etrie de gauge. C'est la performance des \'equations diff\'erentielles fondamentales $FE(\nabla\nabla^*)$ et $FE^*(\nabla)$. Nous avons d\'emontr\'e une analogue de la conjecture fondamentale de Gindikin-Piatecci-Shapiro-Vinberg sur la fibration de vari\'et\'es Kaehleriennes homog\`enes. Notre fibration est topologiquement plus fine et notre d\'emonstration plus courte que celle de Dorfmeister-Nakajima 
%\end{prologue}

\section{INTRODUCTION}

\subsection{The general concerns}
In the category of differentiable manifolds the question whether a given smooth manifold $M$ admits a given geometric structure $\mathcal{S}$ is generally a difficult problem. This problem is named $EX(\mathcal{S})$. A rare geometric structure which exists in every finite dimensional differentiable manifold is the Riemannian structure. For many geometric structures there exit known obstrtuctions to the existence, however there does not exit any characteristic obstruction. The question whether $M$ admits foliations whose leaves (smoothly) support $\mathcal{S}$ is another hard open problem. It is named $EXF(\mathcal{S})$.\\
The problems $EX(\mathcal{S})$ and $EXF(\mathcal{S})$ are the main concerns of this research paper. We involve notions and methods of several mathematical topics in studying those problems. To a pair of gauge structures gauge $\left\{(M,\nabla),(M,\nabla^*) \right\} $ we go to assign three differential operators $D^{\nabla\nabla^*}$, $D^\nabla$ and $D_\nabla$. We involve the global analysis of those differential operators in introducting remarkable Algebras Sheaves. Formally $EXF(\mathcal{S})$ is a problem in the differential topology. The existence of Riemannian foliations and symplectic foliations have been studied in 
\cite{Nguiffo Boyom(6)}. There the concerns were the transverse geometry of foliations. The concern of $EXF(\mathcal{S})$ is the Intrinsic Geometry of Foliations. \textbf{We recall the problems we go to focus on.}\\
\textbf{$EX(\mathcal{S}(M))$: The goal is to address the existence of $\mathcal{S}$ in $M$}.\\
\textbf{$EXF(\mathcal{S}(M))$: The goal is to address the existence of (differentiable) $\mathcal{S}$- foliations in $M$.}\\
We also reformulate $EX(\mathcal{S}(M))$ is terms of a distance-like challenge. This reformulation is denoted by $DL(EX(\mathcal{S}))$. The challenge $DL(\mathcal{S}(M))$ is the search of characteristic invariants which measure how far from addmitting $\mathcal{S}$. We go to give full attention to structures which have two outstanding properties.\\
(1) They are rich in competitive relevant open problems.\\
(2) They are rich in links with other domains of the fundamental research and with other domains the apllied research.  \\
We recall a non exhaustive list of significant open problems.\\
\textbf{$EX(\mathcal{S}: LF(M))$: The existence of locally flat structures in $M$}.\\
\textbf{$EX(\mathcal{S}: Symp(M))$: The existence of symplectic structures in $M$}.\\
\textbf{$EX(\mathcal{S}: Hes(M))$: The existence of Hessian structures in $M$}.\\
Those problems lead to many subsequent open problems the solutionas to which are linked with the applied research. Here are some important instances.\\
\textbf{$EX(\mathcal{S}: Hes(M,g))$: The existence of Hessian structures in a Riemannian manifold (M,g)}. The problem $EX(\mathcal{S}: Hes(M,g))$ is linked with many significant research topics. We go to recall some instnces which objects of current research activities.\\
$EX(\mathcal{S}: Hes(M,g))$ is linked with the information geometry, \cite{Amari}, \cite{Armstrong-Amari}, \cite{Ay-Tuschmann}, \cite{Nguiffo Boyom(6)}, \cite{Murray-Rice}.\\
$EX(\mathcal{S}: Hes(M,g))$ is linked with the thermodynamics \cite{Guts}.\\
$EX(\mathcal{S}: Hes(M,g))$ is linked with the geometric theory of Heat of Souriau \cite{Barbaresco}. Below is another subsequenr open problem.\\
\textbf{$EX(\mathcal{S}: Hes(M,\nabla))$: The existence of Hessian structures in a locally flat manifold $(M,\nabla)$}. This problem $EX(\mathcal{S}: Hes(M,\nabla))$ is linked with the hyperbolicity in locally flat geometry \cite{Kaup}, \cite{Koszul(1)}, \cite{Vinberg(2)}, \cite{Vey(2)} and others.\\
\subsection{The geometry of finite dimensional Cartan-Lie groups and abstract Lie groups}
The geometry of Cartan-Lie groups is a part of the global analysis on manifolds. A subsequent problem in this geometry is the existence of bi-invariant locally flat structures in Cartan-Lie groups and in abstract Lie groups. At one side every finite dimensional Lie group admits left invariant Riemannian metrics. This fact is far from being true for left invariant symplectic forms. The same claim is far from being true for bi-invariant Riemannian metrics. In this paper we address those problems. Before proceding we go to recall a few challenges in the category of finite dimensional Lie groups.\\
\textbf{$EX(\mathcal{S}: lf(G))$: The existence of left invariant locally flat structures in $G$}\\
\textbf{$EX(\mathcal{S}: br(G))$: The existence of bi-invariant Riemannian metrics in a Lie group $G$}.\\
\textbf{$EX(\mathcal{S}: symp(G))$: The existence of left invariant symplectic structures in a Lie group $G$}.\\
\subsection{The information geometry of complex systems}
In the context we go to deal with Big Data are complex systems. Every complex system $\Xi$ may admit many structures of measurable set. Every structure of measurable set $(\Xi,\Omega)$ may admit many statistical models $[\mathcal{E},\pi,M,D,p]$, \cite{Nguiffo Boyom(6)}. In section 24 we raise some ideas about the links of the problems $EX(\mathcal{S})$, $EXF(\mathcal{S})$ and $DL(\mathcal{S})$ with the topological-geometric statistical invariants of Big Data.

\subsection{The gauge group}  
Let $\mathit{G}(M)$ be the gauge group of the vector bundle $TM$ \cite{Petrie-Handal}. Let $\mathcal{MG}(M)$ be the category of gauge structures $(M,\nabla)$, $\nabla$ is a Koszul connection in $M$.\\
\textbf{Reminder: an element of $\mathit{G}(M)$ is an inversible section of the vector bundle $Hom(TM,TM)$, viz $\Phi(x) \in GL(T_xM) \forall x \in M$}. There is a natural action 
$$ \mathit{G}(M)\times \mathcal{MG}(M)\ni (\Phi, \nabla)\rightarrow  \Phi\circ\nabla\circ\Phi^{-1} \in \mathcal{MG}(M)$$
The Koszul connection $\Phi\circ\nabla\circ\Phi^{-1}$ is defined by
$$[\Phi\circ\nabla\circ\Phi^{-1}]_XY = \Phi(\nabla_X\Phi^{-1}(Y).$$
The moduli space of this action is denoted by 
$$\mathbb{MG}(M) = \frac{\mathcal{MG}(M)}{\mathit{G}(M)}.$$
The class of $\nabla$ is denoted by
$$[\nabla]\in \mathbb{MG}(M).$$
We endow the ring $\mathbb{Z}$ with the structure of trivial module of the group $\mathit{G}(M)$. We go to introduce an $\mathit{G}(M)$-equivariant $\mathbb{Z}$-valued function $r^b$ which is defined in $\mathcal{MG}(M)$. So $r^b$ becomes a function defined in the moduli space 
$$ \mathbb{MG}(M)\ni [\nabla]\rightarrow r^b([\nabla])\in \mathbb{Z}.$$
We plan to emphasize the impacts of the function $r^b$ on both $EX(\mathcal{S})$ and $EXF(\mathcal{S})$. We also emphasize the impact of $r^b$ on the combinatorial probem $DL(\mathcal{S})$.\\
\subsection{The overview  of the main results: Solutions to some problems EX(S)}
We give our full attention to the open problems $EX(\mathcal{S}: \Theta))$ that we have mentioned. Then the function $r^b$ enriches every challenge $DL(\mathcal{S}: \Theta)$ with a numerical invariant $r^b(\Theta)$ which has many magnificent statuses. In this overview we go to emphasize the status as characteristic obstruction to $EX(\mathcal{S}: \Theta)$.\\
Consider the problem \textbf{$EX(\mathcal{S}: LF(M))$}, then the function $r^b$ yields a numerical invariant $r^b(M)$ which has the following property
\begin{thm}
In a finite dimensional differentiable manifold $M$ the following assertions are equivalent\\
(1): $r^b(M) = O,$\\
(2): The manifold $M$ admits locally flat structures $\clubsuit$
\end{thm}
Up to nowadays it was hopless to search for such a characteristic obstruction. \cite{Medina-Saldarriaga-Giraldo}\\
The Hessian geometry is rich in relevant links with other research domains. About links with the information geometry the readers are referred to \cite{Barndorff-Nielsen}, \cite{Nguiffo Boyom(6)} and references therein. It is of great interest to know whether a given statistical model is isomorphic to an exponential model, \cite{Amari-Nagaoka}, \cite{Murray-Rice}. This still open problem is called the complexity of a statistical model. By \cite{Nguiffo Boyom(6)} the complexity problem is linked with the Hessian geometry. The subsequent problem is to measure how far from being an exponetial model is a given statistical model. Therefore our approach is the search of distance-like obstructions. We use the function $r^b$ is Riemannian structure or symplectic structure.\\
(1): The relative Riemannian Hessian defect of a Riemannian manifold $r^b(M,g)$,\\
(2): The relative affine Hessian defect of a locally flat manifold $r^b(M,\nabla)$,\\
(3): The absolute Hessian defect of a differentiable manifold $r^B(M)$.\\
Those invariants are the characteristic obstructions in the following meaning.
\begin{thm} (Answer an old question of Alexander K. Guts). In a finite dimensional Riemannian manifold $(M,g)$ the following statements are equivalent\\
$(1.1): r^b(M,g) = 0$,\\
$(1.2):$ the Riemannian manifold $(M,g)$ admits Hessian Riemannian structures $\clubsuit$
\end{thm}
As mentioned this theorem answers a question raised by Alexander K. Guts thirsty years ago \cite{Guts} and independently by S-I Amari \cite{Armstrong-Amari} and references therein.\\
Another source of $EX(\mathcal{S}: He(M,\nabla))$ is the geometry of hyperbolic locally flat manifolds \cite{Kaup}, \cite{Vey(1)}. In this context we have obtained the following result.
\begin{thm} In a locally flat manifold $(M,\nabla)$ the following assertions are equivalent\\
$(2.1): r^b(M,\nabla) = 0$,\\
$(2.2):$ the locally flat manifold $(M,\nabla)$ admits Hessian Riemannian structures $\clubsuit$
\end{thm}
By \cite{Nguiffo Boyom(6)} a locally flat manifold $(M,\nabla)$ admits Hessian structures if and only its 2nd KV cohomology space $H^2_{KV}(\nabla,C^\infty(M))$ contains inversibles symmetric class $[g]$. By \cite{Koszul(1)} for $(M,\nabla)$ being hyperbolic it is necessary that
$B^2_{KV}(\nabla,C^\infty(M)))$ contains an positive symmetric coboundary. This condition is sufficient if $M$ is compact. Those observations are useful for seeing that the geometry of Koszul is a vanishing theorem in the theory of KV cohomology, \cite{Nguiffo Boyom(3)}. In \cite{Armstrong-Amari} J. Armstrong and S-I Amari study the existence of Hessian structures in low dimensional manifolds. Their approach is based on the global analysis. They show that in 2-dimensional analytic manifolds the only formal obstruction is the Pontryagin form. However this approach only yields local solutions. Below we have mentioned the numerical invariant $r^B(M)$. The invariant $r^B(M)$ yields global solutions to $EX(\mathcal{S}: Hes(M))$.
\begin{thm} In a finite dimensional manifold $M$ the following assertions are equivalent\\
$(3.1): r^B(M) = 0$,\\
$(3.2):$ the manifold $M$ admits a Hessian structures $\clubsuit$
\end{thm}
 Those theorems have significant impacts on the information geometry. Currently the information geometry is a domain of international intense research. The exponential statistical models and their generalization are widely studied because of their optimal properties. Some years ago P. McCullagh raised the question as WHAT IS A STATITICAL MODEL? \cite{McCullagh}. The re-establishment of the theory of statistical models is the purpose of \cite{Nguiffo Boyom(6)}. 
Loosely speaking the complexity of a statistical model is the question whether a family of probability densities is an exponential family. We use the function $r^b$ for getting an invariant which measures how far from being an exponential family is a statistical model. About the natation used in the theorem below see the appendix to this paper.
\begin{thm} Consider a regular statistical model $\mathbb{M}$ whose Fisher information is denoted by $g$, namely
$$\mathbb{M} = [\mathcal{E},\pi,M,D,p].$$
The following assertions are equivalent\\
$(4.1): r^b(M,g) = 0,$\\
$(4.2):$ The model $\mathbb{M}$ is an exponential family $\clubsuit$
\end{thm}
In this paper we are concerned with some major problems in the locally flat geometry. Among those problems is the problem of the moduli space of isomorphic class of complete locally flat manifolds. For this purpose we introduce the notion of right spectrum of associative algebras. This notion highlights the links of the function $r^b$ with $EXF(\mathcal{S})$. Those links are useful for emphasizing the differential topological natute of the function $r^b$.\\
An affine Lie group is a pair $(G,\nabla)$ where $\nabla$ is a left invariant locally flat Koszul connection in the Lie group $G$. An abstract Lie group $G$ admits a structure of affine Lie group if and only if its Lie algebra is the commutator Lie algebra of a structure of Koszul-Vinberg algebra. A bi-invariant affine Lie group is a pair $(G,\nabla)$ where $\nabla$ is a bi-invariant locally flat Koszul connection in the Lie group $G$. Mutatis mutandis we get similar definition of affine Cartan-Lie group and bi-invariant affine Cartan-Lie group. In this paper we show that the Locally Flat Geometry is a Byproduct of the bi-invariant Locally Flat Geometry in the category of finite dimensional abstract Lie groups. We use the notion of right spectrum for reducing the affinely flat geometry to the geometry of the category of finite dimensional associative algebras. Those approachs are close to some pioneering ideas of E.B. Vinberg \cite{Vinberg(1)}, \cite{Vinberg(2)}. See also the work of Yozo Matsushima \cite{Matsushima}
\begin{thm} \textbf{The structural theorem}. In the category of gauge structures in a finite dimensional Lie group $G$\\
$(1):$ Every affine Lie group $(G,\nabla)$ is
the INITIAL object of a unique optimal semi-inductive system of finite dimensional bi-invariant affine Cartan-Lie groups
$$\left\{(\Gamma_q,\tilde{\nabla}^q)|\quad \tilde{h}_{pq}:(\Gamma_p,\tilde{\nabla}^p)\rightarrow(\Gamma_q,\tilde{\nabla}^q),\quad p\leq q \right\}.$$
Further that semi-inductive system of bi-invariant affine Cartan-Lie groups is the localization of an semi-inductive system of finite dimensional simply connected bi-invariant affine Lie groups
$$\left\{(G_q,\nabla^q)|\quad h_{pq}:(G_p,\nabla^p)\rightarrow (G_q,\nabla^q)\right\}.$$
$(2):$ In the category of gauge structures in a finite dimensional manifold $M$ every locally flat structure $(M,\nabla)$ is the FINAL object of a unique optimal semi-projective system of bi-invariant affine Cartan-Lie groups
$$ \left\{(\Gamma_q,\tilde{\nabla}^q)|\quad \pi^q_p:(\Gamma_q,\tilde{\nabla}^q)\rightarrow(\Gamma_p,\tilde{\nabla}^p,\quad p\leq q \right\}.$$
Further that semi-projective system of bi-invariant Cartan-Lie groups is the localization of a semi-projective system of finite dimensional simply connected bi-invariant affine Lie groups
$$\left\{(G_q,\nabla^q)\rightarrow (G_p,\nabla^p)\right\}$$
Up to affine isomorphisms a locally flat manifold $(M,\nabla)$ is well defined by its optimal semi-projective system $\clubsuit$
\end{thm}
The theorem just stated may be regarded as another nature of the locally flat geometry. The global geometry of locally flat manifolds is well understood. 
\begin{thm} The Cartan-Lie group of transformations of every locally flat manifold $(M,\nabla)$ is a bi-invariant affine Cartan-Lie group $(\Gamma,\tilde{\nabla})$. Further the connection $\nabla$ is induced by affine dynamic
$$(\Gamma,\tilde{\nabla})\times(M,\nabla)\rightarrow (M,\nabla).$$
Furthermore $(\Gamma,\tilde{\nabla})$ is the localization of a unique simply connected bi-invariant affine Lie group $(G_0,\nabla^0)$ whose infintesimal action in $M$ induces the Koszul connection $\nabla \clubsuit$
\end{thm}
The structural theorem has a relevant corollary that we go to state.
\begin{thm} The group of automorphisms of a (geodesically) complete locally flat manifold $(M,\nabla)$ is a bi-invariant affine Lie group $(G_0,\nabla^0)$. Furtermore the connection $\nabla$ derives from the transitive action
$$(G_0,\nabla^0)\times(M,\nabla)\rightarrow (M,\nabla)\clubsuit$$
\end{thm}
Without the statement of the contrary we go to deal with the category of finite dimensional geometric structures. The symbol $A\leftrightarrow B$ means that the category $A$ is equivalent to the category $B$. Thereby we have the following equivalence of categories\\
E.1: Cartan-Lie Groups (CLG) $\leftrightarrow$ Lie Algebras Sheaves (LAS),\\
E.2: Affinely flat Cartan-Lie Groups (ACLG) $\leftrightarrow$ KV Algebras Sheave (KVAS),\\
E.3: Bi-invariant Affine Cartan-Lie Groups (BCLG) $\leftrightarrow$ Associative Alegebras Sheaves (ASAS).\\
The equivalences we just listed are the localization of the following equivalence of categories.\\
E1.1: Simply connected Lie Groups (LG) $\leftrightarrow$ Lie Algebras (LA),\\
E2.2: Simply connected Affine Lie Groups(ALG) $\leftrightarrow$ Koszul-Vinberg Algebras (KVA),\\
E3.3: Simply connected bi-invariant affine Lie Groups (BALG) $\leftrightarrow$ Associative Algbras (ASA).\\
We also use the function $r^b$ for facing the challenge $DL(EX(\mathcal{S}))$ in the category of Lie groups and the left invariant geometry. Thus we obtain similar numerical constants enjoying expected properties in the category of left invariant gauge structures in Lie groups. \\
\begin{thm} In a finite dimensional abstract Lie group $G$ the following assertions are equivalent\\
$(1): r^b(G) = 0,$\\
$(2): G$ admits left invariant affinely flat structures $\clubsuit$
\end{thm}
We consider the left invariant Hessian geometry in finite dimensional Lie groups. Then the relative Hessian defects $r^b(G,g)$, $r^b(G,\nabla)$ and $r^B(G)$ have the same obstruction nature.\\
 The function $r^b$ is linked with the fundamental equation $FE^*(\nabla)$. We use the other  fundamental equation $FE(\nabla\nabla^*)$ for introducing another numerical functions denoted by $s^b$. The domain of $s^b$ is included in the Lie algebra of infinitesimal gauge transformations. The function $s^b$ imapacts both the symplectic geometry in of finite dimensional manifolds and the left invariant symplectic geometry in finite dimensional Lie groups. It also impacts the bi-invariant Riemannian geometry in the category of finite dimensional Lie. For instance\\
(a) the function $s^b$ is used for studying the problem $EX(\mathcal{S}: br(G)),$\\
(b) the function $s^b$ is used for studying both $EX(\mathcal{S}: Symp(M)$ and $EX(\mathcal{S}: symp(G))$.\\
(i)Classical examples of finite dimensional Lie groups carrying bi-invariant Riemannian metrics are semi simple Lie groups. Their Killing forms are non-degenrate. (ii) Classical examples of finite dimensional Lie groups carrying bi-invariant positive Riemannian metrics are semi simple compact Lie groups. Their Killing forms are negative definite. (iii) Classical examples of symplectic manifods are cotangent bundles endowed with the Liouville form. (iv) No semi simple Lie group admits left invariant symplectic structure. (v) No semi simple Lie group admits left invariant affinely flat structures
\subsection{Geometry of Lie groups}
We go to use the function $s^b$ for investigating the impacts of the differential equation $FE(\nabla\nabla^*)$ on the geometry of finite diemsional Lie groups.\\
\textbf{The bi-invariant Riemannian Geometry in Lie groups}
\begin{thm} In a finite dimensional Lie group $G$ the following properties are equivalent\\
$(1.1): s^b(G) = 0,$\\
$(2.1):$ The Lie group $G$ admits bi-invariant Riemannian metrics $\clubsuit$
\end{thm}
We recall that a Riemannian metric is a nondegenerate symmetric bilinear form $g$. It is called positive if $g$ is positive definite.\\
\begin{thm} In a finite dimensional Lie group $G$ the following properties are equivalent\\
$(1.2): s^b_+(G) = 0.$\\
$(2.2):$ The Lie group $G$ admits bi-invariant positive Riemannian metrics $\clubsuit$
\end{thm}
\textbf{The symplectic Geometry}\\
The function $s^b$ generates the geometric numerical invariant $s^{*b}(M)$ which may be regarded as a characteristic symplectic obstruction.\\
\begin{thm} In a finite dimensional manifold $M$ the following assertions are equilivalent\\
$(1.1):  s^{*b}(M) = 0,$\\
$(1.2):  M$ admits symplectic structures $\clubsuit$
\end{thm}
We also obtain a similar symplectic gap versus Lie groups.
\begin{thm} In a finite dimensional Lie group $G$ the following statements are equivalent\\
$(2.1): s^{* b}(G) = 0,$\\
$(2.2): G$ admits left invariant symplectic structures $\clubsuit$
\end{thm}
\textbf{Important comments}.\\
(Comment1): We take into account Lemma 4 as in \cite{Nguiffo Boyom(3)}, it has settled the algebraic topology of the locally flat geometry, viz the restricted theory of homology generated by the theory of deformation \cite{Gerstenhaber}, \cite{Koszul(1)}  .\\
(Comment2): The introduction of the functions $r^b$ is useful for completely solving
some outstanding geometric problems $EX(\mathcal{S}$ and $DL(\mathcal{S})$.\\
The intrinsinc nature of $r^b$ is also linked with the solution to $EXF(\mathcal{S})$. This is a topological impact of the fundamental equation $FE^*(\nabla)$.\\
(Comment3): By the structural theorem, the finite dimensional Affinely Flat Geometry is the Initial-Final Object of semi Inductive-Projective Systems of bi-invariant Affine geometry in the category of Cartan-Lie groups.\\
(Comment4): $DL(\mathcal{S})$ is completely solved in the category of affinely flat manifolds and Lie groups and in the category of symplectic manifolds and symplectic Lie groups.\\
\textbf{A partial conclusion}.\\
By $Comment1\cap Comment2\cap Comment3\cap Comment4$ we conclude that the problem $DL(\mathcal{S})$ is brought in completion in  both the Locally Flat Geometry and the Symplectic Geometry.\\
We rmind that our concerns include the moduli space of geodesically complete locally flat manifolds and the geometric completeness of geodesically complete locally flat manifolds. For those purpose we have introduced cononical linear representaions of fundamental groups of finite dimensional gauge manifolds. The moduli space of canonical linear representations of fundamental groups is linked with the moduli space of gauge structures. We use those representations for studying the moduli space of geodesically complete locally flat manifolds.\\
\begin{thm} There one to one correspondence between the moduli space (of isomorphism class) of the geodesically complete locally flat manifolds and the moduli space (of conjugacy class) of the canonical representations of their fundamental groups $\clubsuit$
\end{thm}
The canonical linear representations of fundamental groups provide another insight of the topology of finite dimensional locally flat manifolds.
\begin{thm} Generically, the fundamental group of a finite dimensional differentiable manifold has a canonical representation in the group of automorphisms of an object of the following categories.\\
Category.1.1: The category of semi-inductive systems of finite dimensional bi-invariant affine Cartan-Lie groups. \\
category.1.2: The category of semi-projective systems of finite dimensional bi-ivariant affine Cartan-Lie groups.\\
Category.2.1: The category of semi-inductive systems of finite dimensional simply connected bi-invariant affine Lie groups.\\
category.2.2: The category of semi-projective systems of finite dimensional simply connected bi-invariant affine Lie groups.\\
category.3.1: The category of semi-inductive systems of finite dimensional bi-invariant affine Lie groupoids  (as in \cite{Kumpera-Spencer}.\\
category.3.2: The category of semi-projective systems of finite dimensional bi-invariant affine Lie groupoids  (as in \cite{Kumpera-Spencer} $\clubsuit$
\end{thm}

\begin{figure}[H]
\centering
\includegraphics[scale=0.45]{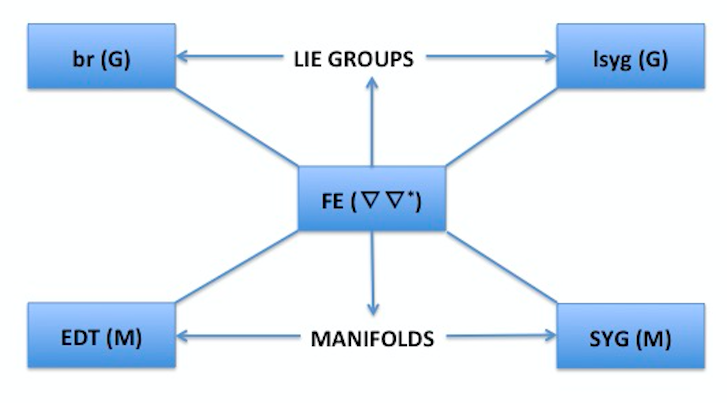}
\caption{Symplectic geometry (where EDT is Extrinsic geometry of foliations)}
\label{fig:symplectic}
\end{figure}

\begin{figure}[H]
\centering
\includegraphics[scale=0.45]{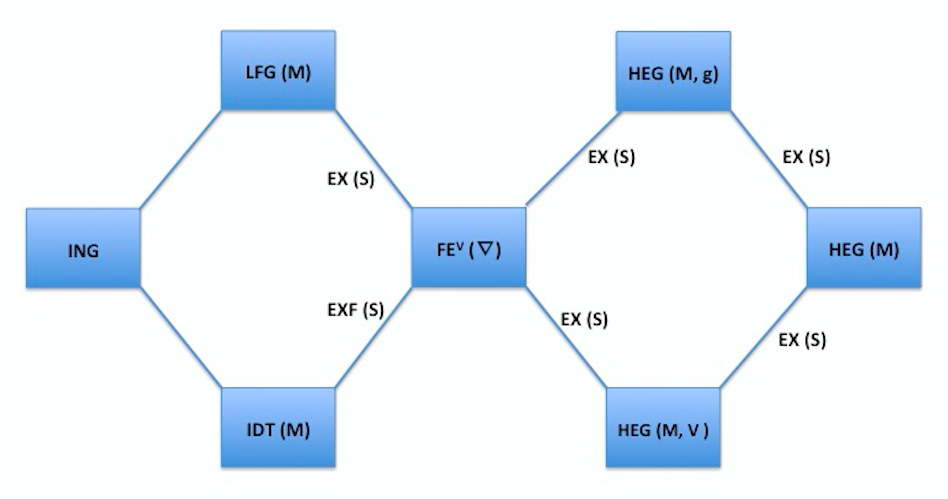}
\caption{Hessian geometry and statistical geometry (where IDT is Intrinsic geometry of foliations)}
\label{fig:hessian}
\end{figure}

\begin{figure}[H]
\centering
\includegraphics[scale=0.45]{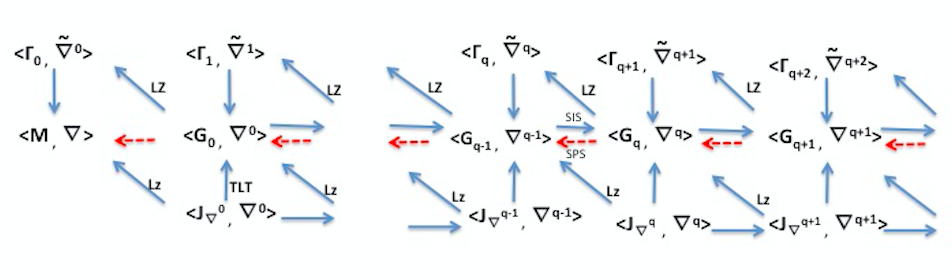}
\caption{Structural theorem in Local Flat Geometry  (where SIS is Semi Inductive System, SPS is Semi Projective System, LZ is Localization and Lz is Lie algebra action)}
\label{fig:edt}
\end{figure}

\begin{figure}[H]
\centering
\includegraphics[scale=0.45]{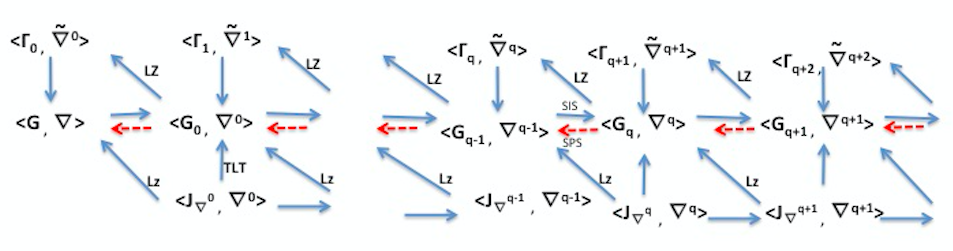}
\caption{Structural theorem in Affine Geometry in Lie Groups (where SIS is Semi Inductive System, SPS is Semi Projective System, LZ is Localization and Lz is Lie algebra action)}
\label{fig:lfg}
\end{figure}

\subsection{A comment on Figues} 
In FIGURE 1\\
\textbf{br(G)} stands for bi-invariant Riemannian geometry in a Lie group(G).\\
\textbf{lsyg(G)} stands for left invariant symplectic geometry in a Lie group G.\\
\textbf{EDT} stands for Extrinsic Differential Topology in $M$. This means the quantitative study of foliations and webs in a manifold $M$ (e.g. the existence of Riemannian foliations, the existence of symplectci foliations.\\
\textbf{SYG(M)} stands for quantitative study of symplectic structures in $M$, the existence problem.\\
\textbf{FIGURE 1 emphasizes some impacts of the differential equation $FE(\nabla\nabla)$}\\
In FIGURE 2\\
\textbf{LFG(M)} stands for Locally Flat Geometry in $M$.\\
\textbf{HEG(M)} stands for HEssian Geometry in a Riemannian manifold $(M,g)$.\\
\textbf{HEG(M)} stands for HEssian Geometry in a manifold $M$.\\
\textbf{$HEG(M,\nabla)$} stands for HEssian geometry in a locally flat manifold $(M,\nabla)$.\\
\textbf{IDT(M)} stands for Intrisinc Differential Topology in a manifolds $M$. This means foliations with prescribed geometric structures in leaves.\\
\textbf{ING} stands for INformation Geometry.\\
\textbf{TLT} stands for Third Lie Theorem.\\
\textbf{FIGURE 2 emphasizes some impacts of the differential equation $FE^*(\nabla)$}\\
In FIGURE 3 and FIGURE 4.\\
\textbf{$<G_q,\nabla^q>$} stands for the simply connected bi-invariant affine Lie group whose associatiove algebra is $ J_\nabla,\nabla>$.\\ 
\textbf{$(\Gamma_q,\tilde{\nabla}^q)$} is the bi-invariant affine Cartan-Lie group whose Associative Algebras Sheaf is $(\mathcal{J}_\nabla,\nabla)$\\ 
\textbf{LZ} stands for LocaliZation of Lie groups. \\
\textbf{FIGURE 3 is (homological algebra) structural theorem of the locally flat geometry in finite dimensional manifolds}.\\
\textbf{FIGURE 4 is the (homological algebra) structural theorem of the left invariant locally flat geometry in Lie finite dimensional Lie groups}
\textbf{FIGURE 1, FIGURE 3 and FIGURE 4 emphasize some impacts of the differential equation $FE^*(\nabla)$}

\newpage
\section{BASIC NOTIONS.}
Let M be an $m$-dimensional differentiable manifold and let $k$ be a non negative integer. The vector bundle of k-jets of sections of the tangent bundle $TM$ is denoted by $J^k(TM)$. At a point $x\in M$ the fiber $J^k_x(TM)$ is the quotient vector space
$$ J^k_x(TM) = \frac{\Gamma(TM)}{I^{k+1}_x(M)}.$$
Here $I_x(M) \subset C^\infty(M)$ is the ideal of differentiable functions which vanish at $x \in M$. Let $Diff_{loc}(M)$ be the Lie pseudogroup of local diffeomorphisms of $M$.\\
\begin{defn} Two local diffeomorphisms $\psi$ and $\phi$ are $k$-equivalent at $x$ is the following requirements are satisfied\\
(1) $\phi(x) = \psi(x)$,\\
(2) $[f\circ\phi^{-1}\circ\psi -f] \in I^{k+1}_x(M)\quad \forall f\in C^\infty(M)$ $\clubsuit$
\end{defn}
To be $k$-equivalent at $x \in M$ is an equivalence relation. The set of equivalence class at $x$ of local diffeomorphisms of $M$ is denoted by $J^k_x(M)$. \\
We put
$$J^k(M) = \cup_{x\in M}J^k_x(M).$$
Given a pair $(x,x*)\in M\times M$ we define $J^k_{xx*}$ by \\
$$  J^k_{xx*}(M) = \left\{[\phi]\in J^k(M) \phi(x) = x*.\right\}$$
By this definition of $J^k_{xx*}(M)$ we get
$$J^k(M) = \cup_{[(x,x*)\in M\times M]} J^k_{xx*}.$$
This presentation is useful for defining the map $(\alpha,\beta)$ of $J^k(M)$ in $M\times M$ by setting
$$[\alpha(\theta),\beta(\theta)] = (x,x*)\quad \forall \theta \in J^k_{xx*}.$$
The triple $[J^k(M),\alpha,\beta]$ is a transitive Lie groupoid the units elements of which are the points of $M$ \cite{Kumpera-Spencer}. It is called the $k^{th}$ Lie groupoid of $M$. The transitivity means that for every pair $(y,y*)\in M\times M$ there exists a class $[phi] \in J^k(M)$ subject to the requirements
$$\alpha([\phi]) = y,$$
$$\beta([\phi]) = y*.$$ The dimension of $M$ is denoted by $m$, we consider
$$ J^k_x(M) = \cap_{\left\{x^*\in M\right\}\left\{J^k_{xx^*}\right\}},$$
The group $J^k_{xx}$ is isomorphic to the general linear group of order $k$, namely 
$GL^k(m,\mathbb{R})$. The projection $\beta$ yields the $J^k_{xx}(M)$-principal bundle
$$  J^k_{xx}(M)\rightarrow J^k_x(M)\rightarrow M $$
which is but the bundle of linear frames of order $k$ of the manifold $M$.\\
For more details on the theory of Lie equations, Cartan-Lie groups and Lie groupoids the readers are referred to Victor Guillemin \cite{Guillemin}, A. Kumpera and D. Spencer \cite{Kumpera-Spencer}, V. Guillemin and S. Sternberg\cite{Guillemin-Stenberg}, Huranish and Rodrigues \cite{KU-RO} I.M. Singer and S. Stenberg \cite{Singer-Sternberg}.\\
In this paper the term [Cartan-Lie Group] means [Group of Lie and Cartan] as in the paper of I.M. Singer and S. Stenberg \cite{Singer-Sterberg}. By the formalism of Spencer the theory of transitive Cartan-Lie groups is linked with the theory of transitive Lie groupoids. This link is based on the theory of Lie equations. Among notable references are Kumpera-Spencer \cite{Kumpera-Spencer}, B. Malgrange \cite{Malgrange}, H. Goldschmidt\cite{Goldschmidt}.\\
Loosely speaking a transitive analytic Cartan-Lie group of order $k$ $\Gamma$ is the analytic solution to a transitive Lie equations $\mathcal{E}$. The later is a system of partial derivative equations of order $k$. The corresponding Lie groupoid is
$$J^k(\Gamma) = \cup_{[x\in M, \phi\in \Gamma]}j^k_x\phi.$$
 Thus $\Gamma$ is the zero of a differential operator $\mathbb{D}$ the symbol of which is $\mathbb{E}$, viz
$$\mathbb{D}(\Gamma) = 0 \leftrightarrow \mathcal{E}(J^K(\Gamma)) = 0. $$
According to our Basic Spelling  Words \\
$BSW.1: \Gamma$ is a Cartan-Lie group of $M$,\\
$BSW.2: \mathcal{E}$ is the Lie equation of $\Gamma$,\\
$BSW.3: J^k(\Gamma)$ is the Lie groupoid generated by $\Gamma,$\\
$BSW.4:$ Lie groups in the usual sense are (often) called abstract Lie groups.\\
\textbf{Reminder}\\
(1) For convenience we go to recall the definition of Lie equation.\\
Let $\mathcal{E}$ be a differential equation of order $k$ the local solutions to which form a subpseudogroup $\Gamma \subset Diff_{loc}(M).$ We have the inclusion
$$J^\infty(\Gamma)\subset \Pi(M).$$
(2) Then $\mathcal{E}$ is a Lie equation if and only if any local diffeomorphism $\phi$ satistfying
$$ j^\infty\phi \in J^\infty\Gamma$$
is an element of $\Gamma \clubsuit$
\subsection{Cartan-Lie groups and abstract Lie groups}
Let $J^k(M)$ be the Lie groupoid of k-jets of local diffeomorphisms of a smooth manifold $M$. Let $(k,k^*)$ be a pair of non negative integers such that
$$k < k^*.$$
We consider the canonical projection
$$J^{k^*}_k:\quad J^{k^*}(M)\rightarrow J^k(M).$$
Every canonical projection is a morphism of Lie groupoid. Further we have
 $$\pi^k_j\circ\pi^j_i = \pi^k_i, \quad\forall i < j < k.$$
So we have the projective (or inverse) limit
$$J^\infty(M) = lim_{k\rightarrow\infty} J^k(M).$$
\textbf{Reminder}\\
To make short we go to put\\
$$\mathcal{P}(M) = Diff_{loc}(M),$$
$$\Pi(M) = J^\infty(M),$$
$$\Pi^k(M) = J^k(M).$$
We have the projective system of Lie groupoids
$$\left\{\Pi^k(M), \pi^j_i | i \leq j \right\}.$$
\textbf{Remarks}\\
(RM.1): The Lie groupoid $\Pi(M)$ is the formal counterpart of the Cartan-Lie Group $\mathcal{P}(M)$.\\
(RM.2): A Cartan-Lie group of order $k$ of a manifold $M$ is a sub-pseudogroup $\Gamma \subset \mathcal{P}(M)$ elements of which form the complete solution to a Lie equation of order $k$.\\
(RM.3): The algebra of a Cartan-Lie group $\Gamma$ is the Lie Algebras Sheaf $\mathcal{L}(\Gamma)$ sections of which are infinitesimal generators of $\Gamma$, viz their local flows belong to $\Gamma$ $\clubsuit$
\begin{defn} Let $G$ be a finite dimensional (abstract) Lie group the Lie algebra of which is denoted by $\bar{G}$.\\
(1) An infinitesimal action of $G$ in a manifold $M$ is a Lie algebra homomorphism $L$ of $\bar{G}$ in the Lie algebra of vector fields in $M$. An infinitesimal action is called effective if the Lie algebra homomorphism is injective.\\
(2) The pseudogroup generated by the local flows of $L(\bar{G})$ is called the localization of the Lie $G$.\\ 
(3) A Cartan-Lie group $\Gamma$ is called finite dimensional if its Lie Algebras Sheaf is generated by an infintesimal action of a finite dimensional abstract Lie group $G$ $\clubsuit$
\end{defn}
\section{THE DIFFERENTIAL EQUATIONS}
\textbf{Fix a geometric structure $\mathcal{S}$. For convenience we recall that this reserach paper is devoted to three problems and to their impacts.}\\
(1) The existence problem $EX(\mathcal{S}$ whose framework is the global analysis on manifolds.\\
(2) The existence of restricted foliations $EXF(\mathcal{S})$ whose framework is the differential topology.\\
(3) The distancelike problem $DL(\mathcal{S})$ whose framework is the combinatorial analysis, viz the study of functions defined in the vertices of graphs.\\
We go to give some instances. Fix a positive Riemannian manifold $(M,g)$.\\
\textbf{$EX(\mathcal{S}: Kaeh(M,g))$}: arises the question whether $M$ admits an integrable almost complex tensor $J$ such that $(M,g,J)$ is a Kaehlerian manifold.\\
This problem impacts the global geometry and the topology of $M$. For instance (i) assume that $M$ is a compact solvmanifold non diffeomorphic to the flat torus, then $M$ does not admit any Kaehlerian structure, \cite{Benson-Gordon}, \cite{McDuff}, \cite{Nguiffo Boyom(7)} and references ibidem. (ii) Assume that $M$ is compact and its fundamental group $\pi_1(M)$ is a free group with two generators, then there is no solution to $EX(\mathcal{S}: Kaeh(M,g))$. 
\textbf{$EX(\mathcal{S}: Kaeh(M,\omega))$}: The problem is whether $M$ admits Kaehlerian structures with a prescribed symplectic structure $(M,\omega)$.\\
\textbf{$EXF(\mathcal{S}: Kaeh(M))$}: The question wheter a smooth manifold $M$ carries a Kaehlerian foliation. An instance is $\mathbb{R}^m\times \mathbb{C}^n$.\\
This section is devoted to settle the background materials which are used for studying the problems  $EX(\mathcal{S})$, $EXF(\mathcal{S})$ and $DL(\mathcal{S})$. \\
Three methods play significant roles. \\
$M.1$: The first method is the homological method. The ingredient is the theory of cohomology of Koszul-Vinberg algebroids \cite{Nguiffo Boyom(3)} \cite{Nguiffo Boyom(4)}. This theory is overviewed in Section 13.\\
$M.2$: The second method is the method of the information geometry \cite{Amari}, \cite{Amari-Nagaoka}. In many parts we involve the functor of Amari.\\
$M.3$: The third method has its source in the global analysis of two differential equations which are encoded by gauge structures, namely $FE(\nabla\nabla^*)$ and $FE^*(\nabla)$. Those equations are useful for introducing two functions $r^b$  and $s^b$. The unction $r^b$ is defined in the vertices of the graph whose vertices are gauge structures in a manifold. This function is invariant under the action of the gauge group $\mathit{G}(M)$. Thereby it goes down in the  moduli space of isomorphism class of gauge structures. This functon is introduced in Section 7. The global analysis of the fundamental equation $FE(\nabla\nabla^*)$ is source of the function $s^b$. The function $s^b$ is defined in the Lie algebra of infinitesimal gauge transformations of $TM$.
\subsection{Notation and definitions}
We go to overview some basic notions.\\
The ring of integers is denoted by $\mathbb{Z}$. The field of real numbers is denoted by $\mathbb{R}$. The field of complex number is denoted by $\mathbb{C}$. Differentiable manifolds are connected and para-compact. Without the statement of the contrary all geometric objects we are conerned with are differentiable. The class of differentiability is $C^\infty$. The associative commutative algebra of real valued differentiable functions in a differentiable manifold $M$ is denoted by $C^\infty(M)$.\\
Let $V$ be a vector bundle over an $m$-dimensional manifold $M$. The vector space of sections of $V$ is denoted by $\mathcal{V}$. It is a left module of the associative commutative algebra $C^\infty(M)$.\\
Let $\mathcal{I}_x$ be the ideal of $C^\infty(M)$ formed of functions which vanish at $x$. Given a non negative integer $k$ the space of $k-jet$ at $x \in M$ of $\mathcal{V}$ is the quotient vector space
$$J^k_x\mathcal{V} = \frac{\tilde{\mathcal{V}}}{\mathcal{I}^{k+1}_x\tilde{\mathcal{V}}}.$$
The vector space $J^k_x\mathcal{V}$ is canonically isomorphic to the vector space
$$ \sum^k_0 [S^j(T^\star_x M)]\otimes \mathcal{V}_x.$$
We denoted by $j^k$ the canonical projection
$$\tilde{\mathcal{V}} \ni s\rightarrow [s]\in J^k_x\mathcal{V}$$
That is to say 
$$j^k(s) = [s] \in J^k(\mathcal{V}).$$
In terms of local coordinate functions $j^k_x(s)$ is the Taylor expansion at $x$ of the section $s$.\\
We put
$$ J^k\mathcal{V} = \cup_{[x\in M]} J^k_x\mathcal{V}.$$
Then the family $\left\{j^k_x, x \in M\right\}$ yields the canonical differential operator
$$\mathcal{V} \ni s\rightarrow j^k(s)\in J^k\mathcal{V}.$$
We consider another vector bundle $W$ over the manifold $M$.
\begin{defn}
A $k^{th}$ order differential operator of $V$ in $W$ is a map $D$ of $\mathcal{V}$ in 
$\mathcal{W}$ admitting a factorization through $j^k$ $\clubsuit$
\end{defn}
\textbf{A comment.}\\
The factorization through $j^k$ means that there exists a mapping
$\tilde{D}$ of $J^k\mathcal{V}$ in $\mathcal{W}$ satisfying the identity
$$\tilde{D}\circ j^k = D.$$
Without the express statement of the contrary we shall be dealing with linear differential operator, viz $D$ is $\mathbf{R}$-linear.\\
\begin{defn}
The mapping $\tilde{D}$ is called the symbol of $D$ $\clubsuit$
\end{defn}
We consider the canonical linear projection $\pi^{k}_{k-1}$ of $J^{k}\mathcal{V}$ onto $J^{k-1}\mathcal{V}$, viz
$$\pi^{k}_{k-1}(j^{k+1}s) = j^{k-1}s.$$
The kernel of $\pi^k_{k-1}$ is denoted by $\tilde{J}^{k}\mathcal{V}$. It is isomorphic to
$$S^k (T^\star M)\otimes\mathcal{V}$$
The restriction to $Ker(\pi^{k}_{k-1})$ of $\tilde{D}$ is denoted by $\sigma(D)$. It is called the principal symbol of $D$.\\
The kernel of the principal symbol $\sigma(D)$ is a subset of
$$ S^k (T^\star M)\otimes [\mathcal{V}^\star\otimes\mathcal{W}],$$
it is called the geometric symbol of $D$.\\
Thus at $x\in M$ the kernel of the principal symbol is denoted by $a_x$. It is a vector subspace of $$Hom(S^k(T_xM),Hom(\mathcal{V}_x,\mathcal{W}_x))$$
\subsection{The Kuranishi-Spencer formalism}
We keep the notation used above. For $s \in \mathcal{V}$ and $x\in M$ one has
$$ D(s)(x) = \tilde{D}(j^k_xs).$$
Now  $\frac{\partial}{\partial x}$ stands for all of the partial derivatives with respect to a local coordinates of the point $x$.\\
(1) $J^k\mathcal{V}$ is a vector bundle over $M$.\\
(2) $\forall s \in \mathcal{V}$ we have $$j^1(j^ks) = j^{k+1}s.$$
Now we set
$$P^1(D)(s)(x) = \frac{\partial D(s)(x)}{\partial x}$$
\begin{defn}
$P^1(D)$ is called the first (Kuranishi) prolongation of the operator $D$.
\end{defn}
Inductively, for $q > 0$ one defines the $q^{th}$ prolongation of $D$ by setting
$$ P^{q}(D) = P^1(P^{q-1} (D))$$
The geometrical symbol of $P^q(D)$ is denoted $a^q$. It is easily seen that
$$a^{q+1} = Hom(S^{q+1}TM,Hom(\mathcal{V},\mathcal{W}))\cap Hom(TM,a^q).$$
\subsubsection{The Koszul-Spencer complex of differential operators}
A major problem in the global analysis is the question whether a formally integrable differential equation is analytically integrable. This question has had a long story. The answer to this question is known as the LEWY PROBLEM \cite{Nirenberg}. The Lewy problem deals with the integrability of the following System of Partial Differential Equations in $\mathbb{R}^3$
$$ \frac{\partial u(z,y)}{\partial\bar{z}} + \frac{\partial u(z,y)}{\partial y} = \frac{d\psi(y)}{dy}.$$
Here $\psi(y)$ is a complex valued differentiable function defined in $\left\{O\right\}\times \mathbb{R}$, $u(z,y)$ is a unknown complex valued differentiable function,
$u(z,y)$ is defined is $\mathbb{C}\times \mathbb{R}$. \\
This SPDE is formally integrable, however if $\psi(y)$ is not analytic there is no differentiable solution to this equation.\\
The formal integrability is controled by the cohomology of the Koszul-Spencer complex of the geometrical symbol. We go to outline the theory of Koszul-Spencer homology and to highlight its impact on the problem of formal integrability of SPDE.   \\
We consider the bi-graded vector space $$ C = \oplus_{p,q}C^{p,q}. $$
The $(p,q)$-homogeneous subspace is defined by
$$ C^{p,q} = \Omega^p(M)\otimes a^q.$$
Here $\Omega^p $ is the vector space of differential q-forms in $M$. We consider a monomial $\xi \in C^{p,q+1}$
$$\xi = \omega_1\wedge..\wedge \omega_p\otimes \omega^\star_1.\omega^\star_2...\omega^\star_{q+1}\otimes L.$$
The Koszul-Spencer coboundary operator $d$ maps $C^{p,q+1}$ in $C^{p+1,q}$. It is defined by
$$d\xi  = \sum^{q+1}_1 (-1)^j \omega_1\wedge\omega_2\wedge..\wedge\omega_p\wedge\omega^\star_j\otimes \omega^\star_1..\hat{\omega}^\star_j..\omega^\star_{q+1}\otimes L.$$
By direct calculations one gets\\
$$ d^2\xi = 0. $$ 
Thus we get the Koszul-Spencer homology complex $(C,d)$
$$ \rightarrow C^{p-1,q+1}\rightarrow C^{p,q}\rightarrow C^{p+1,q}\rightarrow $$
Its homology space at the level $C^{p,q}$ is denoted by $H^{p,q}(a)$.\\
\subsubsection{The involutivity after E. Cartan}
Let $V$ and $W$ be finite dimensional vector spaces and let $a$ be a vector subspace of the vector space $Hom(V,W)$.\\
\begin{defn} The first prolongation of $a$ is defined by
$$a^{(1)} = Hom(V,a)\cap Hom(\mathcal{S}^2(V),W) \clubsuit$$
\end{defn}
Now we consider a basis of $V$
$$\left\{v_1,v_2,..,v_m \right\}.$$
We consider a non negative integer $j \leq m$.\\
\begin{defn} We defined the subspace $\mathcal{a}_j$ by
$$a_j = \left\{A\in a | A(v_i) = 0,\quad 1\leq i\leq j \right\}\clubsuit$$
\end{defn}
The following inequality is due to Elie Cartan,\cite{Singer-Stenberg}, \cite{Guillemin-Stenberg}
$$ dim(a^{(1)}) \leq \sum^m_0 dim(a_j).$$
\begin{defn} With the notation just used
(1) a basis $[v_1,v_2,..,v_m]$ is called quasi regular for $a$ if
$$ dim(a^{(1)}) = \sum^m_0 dim(a_j) $$
(2) $a \subset Hom(V,W)$ is called involutive if there exists a quasi regular basis for it $\clubsuit$
\end{defn}
\subsubsection{The involutivity after J.P. Serre}
We consider the Koszul-Spencer complex of $a$
$$\rightarrow \wedge^p V\otimes a^{q +1}\rightarrow \wedge^{p+1}V\otimes a^q \rightarrow $$
\begin{defn} $a$ is called involutive if
$$ H^{p,0}(a) = 0\quad \forall p.$$
\end{defn}
The following theorem is due to Jean-Pierre Serre, see the appendix to \cite{Guillemin-Stenberg}
\begin{thm} The following assertions are equivalent\\
(1) $H^{p,q}(a) = 0 \forall p > 0,$\\
(2) there exists a quasi regular basis for $a$ $\clubsuit$
\end{thm}
\begin{defn}
A differential operator whose geometric symbol is involutive is called involutive $\clubsuit$
\end{defn}
\subsection{Some categories of geometric structures}
A Koszul connection is a first order differential operator $\nabla$ of $TM^{\otimes 2}M$ in $TM$. It is usually denoted by
$$ (X,Y)\rightarrow \nabla_XY.$$ 
It has the following properties
$$(1):\quad \nabla_{fX}Y) = f\nabla_XY),$$
$$(2): \quad \nabla_XfY) = df(X)Y + f\nabla_XY\quad \forall f \in C^\infty(M), \forall X, Y \in \mathcal{X}(M).$$
The torsion tensor $T^\nabla$ and the curvature tensor $R^\nabla$ are defined by
$$T^\nabla(X,Y) = \nabla_XY - \nabla_YX - [X,Y],$$
$$R^\nabla(X,Y,Z) = \nabla_X\nabla_YZ - \nabla_Y\nabla_XZ -
\nabla_{[X,Y]}Z.$$
\subsubsection{The category $\mathcal{LF}(M)$}
In the introduction we have raised a few major open problems in the pure differential geometry and in the applied differential geometry.\\
Let $\mathcal{LF}(M)$ be the category whose objects are gauge structures $(M,\nabla)$ having the following properties 
$(1):\quad T^\nabla (X,Y) = 0,$\\
$(2):\quad R^\nabla (X,Y) = 0 \forall X,Y.$
An object of $\mathcal{LF}(M)$ is called a locally flat structure in $M$.\\
We recall the open problems mentioned in the introduction.\\
\textbf{$EX(\mathcal{S}: LF(M))$: The existence of locally structures in $M$. This is a very old open problem. See \cite{Smilie} and references therein $\clubsuit$}
\subsubsection{The category $\mathcal{HGE}(M,g)$.}
\textbf{$EX(\mathcal{S}: Hes(M,g))$: The existence of Hessian structures in a Riemannian manifold $(M,g)$}\\
The problem $EX(\mathcal{S}: Hes(M,g))$ has a long story. Thirty years ago, in a private communication to the author Alexander K. Guts raised this problem. In the recent correspondence between Alexander K. Guts and the author the main concern was this old problem $EX(\mathcal{S}: Hes(M,g))$. Independently S-I Amari raised the same question many years ago. In the information geometry this question is linked with the comlexity of statistical models 
\cite{Nguiffo Boyom(6)}. Recently J. Armstrong and S-I Amari have studied this problem in \cite{Armstrong-Amari}. See also S-I Amari, J. Armstrong and A.K. Guts are interested in positive Riemannian metrics. In this paper we deal with the general case , viz including the pseudo-Riemannian geometry. We recall that the case of positive Riemannian metrics is also linked with the Geometry of Koszul. Indeed in the category of compact positive Riemannian manifolds the question $EX(\mathcal{S}: Hes(M,g))$ is equivalent to the question whether $\tilde{M}$ is diffeomorphic to a convex domain not containing any straight line in $\mathbb{R}^m$ \cite{Koszul(1)}, \cite{Koszul(4)}. For convenience we recall that a Hessian structure in Riemannian manifold $(M,g)$ is a triple $(M,g,\nabla)$. Here $\nabla$ is a Koszul connection which has the following properties\\
$(1):\quad T^\nabla = 0,$\\
$(2):\quad R^\nabla = 0,$\\
$(3):\quad (\nabla_Xg)(Y,Z) - (\nabla_Yg)(X,Z) = 0.$\\
The requirement (3) has a homological nature. It means that the metric tensor $g$ is a 2-cocycle of the KV complex of the locally flat manifold $(M,\nabla)$
\cite{Shima(2)} \cite{Nguiffo Boyom(6)}.\\
The category whose objects are Hessian structures in a given Riemannian manifold $(M,g)$ is denoted by $\mathcal{HGE}(M,g)$.
\subsubsection{The category $\mathcal{HGE}(M,\nabla)$.}
Let $(M,\nabla)$ be an object of the category $\mathcal{LF}(M)$.\\
$EX(\mathcal{S}: Hes(M,\nabla))$: \textbf{The existence of Hessian structures in a locally flat manifold $(M,\nabla)$.}\\
The problem is the search of a Riemannian metric tensor $g$ such that the triple $(M,\nabla,g)$ is a Hessian manifold. That is the central problem in the geometry of Koszul and its applications, \cite{Koszul(1)}, \cite{Kaup}, \cite{Vey(1)}, \cite{Barbaresco}, \cite{Shima(2)}.\\
Let $H^2_{KV}(M,\mathbb{R})$ be the $2nd$ scalar KV chomology of a locally flat manifold $(M,\nabla)$. Every solution to $EX(\mathcal{S}: Hes(M,\nabla))$, namely $(M,g)$ yields the cohomology class $[g] \in H^2_{KV}(M,\mathbb{R}$ which "measures" how far from being hyperbolic is $(M,\nabla)$. \cite{Kaup}, \cite{Koszul(1)}. See the recent attempt of J. Armonstrong and S-I Amari in \cite{Armstrong-Amari}. The category whose objects are Hessian structures in a locally flat manifold $(M,\nabla)$ is denoted by $\mathcal{HGE}(M,\nabla)$.
\subsubsection{The category $\mathcal{HGE}(M)$}
$EX(\mathcal{S}: Hes(M))$: \textbf{The existence of Hessian structures in a manifold $M$.}\\
The Hessian geometry has major impacts on many topics. Among those topics are the information geometry, The topology and the geometry of bouded domains. Sadely $EX(\mathcal{S}: Hes(M))$ is a difficult problem. The author does not know any reference of succeesful attempts. This general problem is studied in this paper. The category whose objects are Hessian structures in $M$ is denoted by $\mathcal{HGE}(M)$. We have already highlighted the role played by the Hessian geometry in the theory of statistical models. For more details the readers are referred to \cite{Amari}, \cite{Barndorff-Nielsen}, \cite{Murray-Rice}, \cite{Nguiffo Boyom(6)}$\clubsuit$
\subsubsection{The category $\mathcal{EXF(S)}$}
We are concerned with a few problems in the differential topology. We go to deal with restricted  foliations, viz foliations with restricted intrinsic geometric structures of leaves. Before pursuing we recall another formulation of the starting fundamental problems,\\
$EXF(\mathcal{S})(M)$: \textbf{The existence of $\mathcal{S}$-foliations in $M$.} In general $EXF(\mathcal{S})(M)$ is also a hard problem. If $\mathcal{S}$ stands for positive Riemannian structure then every foliation is a $\mathcal{S}$-foliation. In the category of positive Riemannian manifolds $EXF(\mathcal{S})$ is reduced to the existence of foliations. In applied differential geometry many interesting involutive distributions are singular distributions. Sadely the classical theorem of Frobenius does not work for all involutive singular distributions. Important instances of involutive singular distribution are the kernels of Fisher information of singular statistical models \cite{Nguiffo Boyom(6)}. The category of finite dimensional manifolds admitting $\mathcal{S}$-foliations is denoted by $\mathcal{EXF(S)}$.
The global analysis of the fundamental equation $FE^*(\nabla)$ will be used for exploring this category $\mathcal{EXF(S)}$. This highlight the range of impacts of that differential equation
\subsubsection{The differential operators $\left\{D^\nabla, D_\nabla, D^{\nabla\nabla^\star}\right\}$}
This subsubsection is devoted to the central source of the arsenal we go to use for addressing the central problem $EX(\mathcal{S})$. We go to use a pair of gauge structures for introducing three differential operators. Those operators yield fundamental differential equations whose global analysis helps to introduce the functionnal geometric invariants $r^b$ and $s^b$. We consider a pair of gauge structures $\left\{(M,\nabla),(M,\nabla^*)\right\}$. We go to introduce three differential operators which are denoted by $D_\nabla$, $D^\nabla$, and $D^{\nabla\nabla^*}$. Let $R^\nabla$ be the curvature tensor of $\nabla$. The Lie derivative of $\nabla$ in the direction $X$ is denoted by $L_X\nabla$. The inner product of $R^\nabla$ by $X$ is denoted by $\iota_XR^\nabla$. The differential operators $D^\nabla$ and $D_\nabla$ are defined in the tangent vector bundle $TM$. Their values belong to the vector bundle $hom(\otimes^2 TM, TM)$. Thus for a vector field $X$  $D^\nabla(X)$ and $D_\nabla(X)$ are sections of the vector bundle $TM\otimes T^{*\otimes 2}M$. Those sections are defined by
$$(1):\quad D^\nabla(X) = L_X\nabla - \iota_XR^\nabla,$$
$$(2):\quad D_\nabla(X) = \nabla^2X.$$
Thereby the complete expressions are
$$D^\nabla(X).(Y,Z) = (L_X\nabla)(Y,Z) - R^\nabla(X,Y).Z,$$
$$D_\nabla(X).(Y,Z) = \nabla_Y\nabla_ZX - \nabla_{\nabla_YZ}X\quad \forall Y, Z \in \mathcal{X}(M).$$
The differential operator $D^{\nabla\nabla^*}$ is defined in the Lie algebra of infintesimal gauge transformations of $M$. Its values also belong to the vector bundle $TM\otimes T^{*\otimes 2}M$. The differential operator $D^{\nabla\nabla^*}$ is defined by
$$(3):\quad D^{\nabla\nabla^*} (\psi) = \nabla^* \circ \psi - \psi \circ \nabla \quad \forall \psi,$$
\textbf{Reminder: Here $\psi$ is a section of the vector bundle $Hom(TM,TM)$}. So we have
$$ [D^{\nabla\nabla^*}(\psi)].(X,Y) = \nabla^*_X\psi(Y) - \psi(\nabla_XY)\quad \forall X, Y.$$
MIKE figure
TM - Vide- TM   
$D^\nabla - D_\nabla$
$vide-T^*M^{\otimes2} \otimes TM$
$D^{\nabla\nabla^*}$
$vide-T^*M\otimes TM-vide $
\section{THE FUNDAMENTAL EQUATIONS}
We go to introduce three differential equations whose global analysis impacts the problems $EX(\mathcal{S})$, $EXF(\mathcal{S})$ and $DL(\mathcal{S})$.  
\subsection{Notation}
We choose a system of local coordinate functions
$$\left\{x_1,..,x_m\right\},$$
we set
$$ \partial_i = \frac{\partial}{\partial x_i}$$
$$\nabla_{\partial_i}{\partial_j} = \sum_k\Gamma_{ij:k}{\partial_k}$$
$$\nabla^*_{\partial_i}{\partial_j} = \sum_k\Gamma^*_{ij:k}{\partial_k}$$
$$ X = \sum_jX^j\partial_j.$$
The operators $D^\nabla$, $D_\nabla$ and $D^\nabla\nabla^*$ are $\mathbb{R}$-linear. Further their values are sections of the same vector bundle $Hom(TM^{\otimes 2},TM)$. Our concerns are the kernels of those differential operators. We go to introduce the differential equations whose global analysis has impacts on the problems $EX(\mathcal{S})$, $EXF(\mathcal{S})$ and $ DL(\mathcal{S})$
\subsection{The fundamental gauge equation and its link with foliations}
\textbf{Reminder} The Lie algebra of infintesimal gauge transformations is denoted by $g(M)$. 
\begin{defn} The fundamental gauge equation is the first order differential equation $FE(\nabla,\nabla^*)$ which is defined  in $g(M)$ by
$$FE(\nabla,\nabla^*):\quad D^{\nabla\nabla^*} (\psi) = O \clubsuit$$
\end{defn}
In terms of local coordinate functions $\psi \in Hom(TM^{\otimes 2},TM)$ is the matrix 
$$ [\psi(x)] = [\psi_{kj}(x_1,..,x_m)].$$ 
Therefore $FE(\nabla,\nabla^*)$ is the following SPDE
$$ S_{ij:k}:\quad \frac{\partial\psi_{kj}}{\partial_i} - \Sigma_{1\leq t\leq m}[\Gamma_{ij:t}\psi_{kt} - \Gamma^\star_{it:k}\psi_{tj}] = 0. $$
 The sheaf of germs of local solutions to $FE(\nabla,\nabla^\star)$ is denoted by
 $\mathcal{J}_{\nabla\nabla^*}$. The vector space of sections of $\mathcal{J}_{\nabla\nabla^\star}$ is denoted by $J_{\nabla\nabla^*}$.\\
 \subsubsection{Riemannian foliations, symplectic foliations and the fundamental equation $FE(\nabla\nabla^*)$}
We go to relate the solutions to fundamental equation $FE(\nabla\nabla^*)$ and the differential topology. A solution to the first fundamental equation $\psi$ satisfies the relation
$$FE(\nabla\nabla^*):\quad \nabla^*\circ \psi - \psi\circ \nabla = 0.$$
\text{Reminder} In \cite{Nguiffo Boyom(6)} we have involved the sheaf $\mathcal{J}_{\nabla\nabla^*}$ in studying Riemannian foliations and symplectic foliations. There the study was quantitative. For convenience we go to recall the useful notation.\\
(i) The sheaf of germs of differential 2-forms in $M$ is denoted by $\Lambda_2(M)$. The vector space of sections of $\Lambda_2(M)$ is denoted by $\Omega^2(M)$. (ii) The sheaf of germs of symmetric bi-linear forms in $M$ is denoted by $S_2(M)$. The vector space of sections of $S_2(M)$ is denoted by $\mathcal{S}_2(M)$. (iii) The Lie algebra of infinifinitesimal gauge transformations of $M$ is denoted by $g(M)$. The bracket of two infinitesimal gauge transformations is defined by
 $$[\psi,\phi](x) = \psi(x)\circ\phi(x) - \phi(x)\circ\psi(x)\quad \forall x \in M.$$
\begin{defn} The map
 $$\mathcal{R}ie(M)\times g(M) \ni (g,\psi)\rightarrow (q(g,\psi),\omega(g,\psi) \in S_2(M) \times \Omega^2(M)$$
 is defined by 
 $$ [q(g,\psi)](X,Y) = \frac{1}{2}[g(\psi(X),Y) + g(X,\psi(Y))],$$
 $$ [\omega(g,\psi)](X,Y) = \frac{1}{2}[g(\psi(X),Y) - g(X,\psi(Y))] \clubsuit$$
 \end{defn}
\begin{defn} To every pair $(g,\psi) \in \mathcal{R}ie\times \mathit{g}(M)$ we assign the unique pair $(\Psi,\Psi^*) \subset \mathit{g}(M)$ which is defined by\\
 $$  g(\Psi(X),Y) = [q(g,\psi)](X,Y),$$
 $$  g(\Psi^*(X),Y) = [\omega(g,\psi)](X,Y).$$
 \end{defn}
 \begin{prop}
 For every triple $(g,\nabla,\psi)\in, \mathcal{R}ie(M)\times \mathcal{LC}(M)\times g(M)$ the following assertions are equivalent\\
 (1) $(q(g,\psi),\omega(g,\psi))\in S^\nabla_2(M)\times \Omega^\nabla_2(M),$\\
 (2) $\Psi$ and $\Psi^*$ are solutions to $FE(\nabla\nabla^g)$ $\clubsuit$
 \end{prop}
 \textbf{Hint}: Perform the formulas 
 $$[\nabla_X q(g,\psi)](Y,Z) = [q(g,\psi)](D^{\nabla\nabla^g}\Psi(X,Y),Z),$$
 $$[\nabla_X\omega(g,\psi)](Y,Z) = [\omega(g,\psi)](D^{\nabla\nabla^g}\Psi^*(X,Y),Z).$$
 \textbf{Reminder}: Assume that $\nabla \in \mathcal{SLC}(M)$. (a) Global sections of $\mathcal{S}^\nabla_2(M)$ are Riemannian foliations; (b) Global sections of $\Lambda^\nabla_2(M)$ are symplectic foliations \cite{Nguiffo Boyom(6)}. By the functor of Amari every pair $(g,\nabla)\in \mathcal{R}ie(M)\times \mathcal{SLC}(M)$ defines the fundamental equation $FE(\nabla\nabla^g)$. The sheaf $\mathcal{J}_{\nabla\nabla^g}$ is linked with the pair ($\Lambda^\nabla_2 (M), \mathcal{S}^\nabla_2(M))$ according to the splitting the exact sequence
$$O\rightarrow \Lambda^\nabla_2(M)\rightarrow \mathcal{J}_{\nabla\nabla^*}\rightarrow S^\nabla_2(M)\rightarrow 0.$$
Therefore the vector space of solutions to the fundamental equation $FE(\nabla\nabla^g)$ is linked with the pair ([Riemannian foliations], [symplectic foliations]) 
according to the splitting exact sequence
$$ 0 \rightarrow\Omega^\nabla_2(M)\rightarrow J_{\nabla\nabla^g}\rightarrow S^\nabla_2(M)\rightarrow 0.$$
Both $q((g,\psi)$ and $\omega(g,\psi$  are solutions to the problem $EXF(\mathcal{S}$, depending on $\mathcal{S}$ is the Riemannian structure or the symplectic structure. 
What we just discussed are impacts of methods of the Information Geometry and of the global analysis  on the differential topology.
\subsubsection{The fundamental equation $FE(\nabla\nabla^*)$ and the distancelike problem $DL(\mathcal{S})$}
We take into account the results in the last subsubsection, in particular the relationship between the fundamental equation $FE(\nabla,\nabla^g)$ and the intrinsic geometry of regualr foliations. The concern is to use that relationship for exploring the symplectic geometry in finite dimensional manifolds and the bi-invariant Riemannian geometry in finite dimensional Lie groups. We are mainly concerned by two requests,\\
 \textbf{$DL(\mathcal{S}: Sym(M))$: How far from admitting symplectic structures is a finite dimensional manifold $M$?}\\
\textbf{$DL(\mathcal{S}:br(G))$: How far from admitting bi-invariant Riemannian metrics is a finite dimensional Lie group $G$?}\\
In the litterature such as \cite{Pennec}, \cite{Medina} and \cite{Medina-Revoy} the authors are concerned with the problem $EX(\mathcal{S}:br(G))$. The case of bi-invariant positive Riemannian metrics is easy. Lie groups carrying positive Riemannian metrics are direct product of abelian Lie groups and semi-simple compact Lie groups. In this paper we deal with the general case, viz including the pseudo-Riemannian geometry.\\
Now we go to use the global ananlysis of the fundamental differential equation $FE(\nabla\nabla^g)$ for introducing a new $\mathbb{Z}$-valued function 
$$\mathcal{R}ie(M)\times\mathcal{SLC}(M) \ni (g,\nabla)\rightarrow s^b(g\nabla)\in \mathbb{Z}.$$ 
The reader will easily see that those function $s^b$ is defined in $\mathcal{R}ie(M)\times\mathbb{MG}(M)$. Here $\mathbb{MG}(M)$ is the moduli space of gauge structures in $M$.
\begin{defn} The function $s^b(g,\nabla)$ is defined by
$$s^b(g,\nabla) = \min_{[\psi\in J_{\nabla\nabla^g}]} \left\{dim(M) - rk(\Psi)\right\};$$
the function $s^{*b}(g,\nabla)$ is defined by
$$s^*b(g,\nabla) = \min_{[\psi \in J_{\nabla\nabla^g}}\left\{dim(M) - rk(\Psi^*)\right\}; $$
the numerical invariant $s^b(M)$ is defined by
$$s^b(M) = \min_{[(g,\nabla) \in \mathcal{R}ie(M)\times\mathcal{SLC}(M)]}\left\{s^b(\nabla,\nabla^g)\right\} \clubsuit$$
\end{defn}
\textbf{warning}: here $rk(\Psi(x))$ stands for the rank of the linear map
$$\Psi(x) : T_xM\rightarrow T_xM$$
The impacts of functions just defined will be addressed in the section devoted to both the symplectic geometry in finite dimensional manifolds and the bi-invariant Riemannian geometry in finite dimensional Lie groups. The next subsection is devoted to the affine geometry and its  appications.
\subsection{Two numerical fundamental equations}
\begin{defn} The first (numerical) fundamental equation of a gauge structure $(M,\nabla)$ is the 2nd order differential equation
$$FE^*(\nabla):\quad  D_\nabla(X) = 0\quad\clubsuit $$
\end{defn}
The local versus of the first fundamental equation $FE^*(\nabla)$ is the following system of partial derive equations
$$ \left\{\Theta^k_{ij}(X) = 0,\quad i,j,k \in \left\{1,..,m\right\}\right\}$$
\begin{defn} The second (numerical) fundamental equation of a Koszul connection $\nabla$ is the 2nd order differential equation
$$FE^{**}(\nabla):\quad D^\nabla(X) = 0\quad\clubsuit$$
\end{defn}
Both $FE*(\nabla)$ and $FE**(\nabla)$ have the same principal symbol which is denoted by $\sigma(\nabla).$
\subsubsection{A sketch of the global analysis of $FE^*(\nabla)$: the analytic integrability}
The geometric symbol of $FE^*(\nabla)$, namely $a$ is defined by the SPDE
$$\frac{\partial^2X^k}{\partial_i\partial_j} = 0.$$
Therefore both the equation $FE^*(\nabla)$ and the equation $FE^{**}(\nabla)$ are involutive and elliptic. Thereby we can implement the Kuranish-Spencer formalism, \cite{Spencer}, \cite{Kuranishi}, \cite{Goldschmidt}, \cite{Guillemin}, \cite{Malgrange}. Thus we can state the following proposition
\begin{prop} The PDE $FE^*(\nabla)$ and $FE^{**}(\nabla)$ are pointwise formally integrable $\clubsuit$
\end{prop}
The sheaf of germs of solutions to the first fundamental equation $FE^*(\nabla)$ is denoted by 
$\mathcal{J}_\nabla$; $J_\nabla$ is the vector space of sections of $ \mathcal{J}_\nabla$. The sheaf of germs of solutions to the second fundamenatal equation $FE^{**}(\nabla)$ is denoted by $\mathcal{J}(\nabla)$; $J(\nabla)$ is the space of sections of $\mathcal{J}_\nabla$. 
\begin{lem} The pair $\left\{\mathcal{J}_\nabla,\nabla\right\}$ is an Associative Algebras Sheaf (AAS) $\clubsuit$
\end{lem}
\textbf{Proof}\\
Let $\xi$ and $\xi^*$ be sections of $\mathcal{J}_\nabla$. Let $X$ and $Y$ be two vector fields in $M$. To make simple we put
$$X.Y = \nabla_XY.$$
Therefore we have
$$\nabla^2\xi = 0,$$
$$\nabla^2\xi^* = 0 .$$
Thus for all vectors $X, Y \mathcal{X}(M)$ we have
$$ X.(Y.\xi^\star) - (X.Y).\xi^\star = 0.$$
We calculate
$$ X.(Y.(\xi.\xi^\star)) - (X.Y).(\xi.\xi^\star) = X.((Y.\xi).\xi^\star) -
((X.Y).\xi).\xi^\star $$
$$ = [(X.(Y.\xi)) - (X.Y).\xi].\xi^\star $$
Since both $\xi$ and $\xi^\star$ are sections of the sheaf $\mathcal{J}_\nabla$ we have
$$\nabla^2\xi = 0,$$
$$\nabla^2\xi^\star = 0.$$
In final we have
$$\nabla^2(\xi.\xi^\star) = 0.$$
The lemma is proved $\clubsuit $
At a point $x \in M$ the vector subspace of $T_xM$ that is spanned by locall solutions to $FE*(\nabla$ is denoted by
$\mathcal{T}_\nabla M(x) \subset T_xM $.
\begin{lem} For every torsion free gauge structure $(M,\nabla)$ the first fundamental equation $FE^*(\nabla)$ coincides with the second fundamental equation $FE^{**}(\nabla) \clubsuit$
\end{lem}
\textbf{Remark.}\\
In view of the expression of their principal symbol the fundamental equations $FE^*(\nabla)$ and $FE^{**}(\nabla)$ are of finite type, \cite{Kumpera-Spencer}, \cite{Guillemin}, \cite{Singer-Sternberg}. Consequently the real vector spaces $J_\nabla$ and $J(\nabla)$ are finite dimensional. 
We go to implement the formalism of Cartan-Lie group as in Singer-Sternberg \cite{Singer-Sternberg}. Among the significant properties of the first fundamental equation the following is helpful for our purposes.
\begin{thm} For every torsion free Koszul connection $\nabla$ the sheaf $\mathcal{J}_\nabla$ is a LAS (Lie Algebras Sheaf) in the sense of Singer-Sternberg \cite{Singer-Sternberg} $\clubsuit$
\end{thm}
\textbf{Hint}. Theorem is based on the fact that $\mathcal{J}_\nabla$ is defined an equation.\\
\textbf{Reminder}: A Lie Algebra Sheaf (LAS in short) $\mathcal{L}$ is completely defined by its formal version 
$$ J^{\infty}(\mathcal{L}) = \cup_{x \in M} \left\{j^\infty_x \mathcal{L}\right\},$$
Loosely speaking a vector field $X$ is a section of $\mathcal{L}$ if ans only if $j^\infty_x(X)\in J^\infty(\mathcal{L})$  $\forall x \in M$. If this condition fails $\mathcal{L}$ is called a Weakly Lie Algebras Sheaf, WLAS in short \cite{Singer-Stenberg}. Thus a vector field $X$ belongs to $J_\nabla$ if and only if  
$$j^\infty_x(X)\in \mathcal{J}^\infty_\nabla(x)\quad \forall x\in M \clubsuit $$
Before continuing we recall that the differential operators $D^\nabla$ and $D_\nabla$ are involutive of type two. According to the formalism of Kuranishi-Spencer the fundamenatal equations $FE^*(\nabla)$ and $FE^{**}(\nabla)$ are pointwise formally integrable. If they are transitive then they are analytically integrable. 
\subsection{A sketch of the global analysis in symmetric gauge structures}
Without the statement of the contrary we go to focus on symmetric gauge structures $(M,\nabla)$, viz $\nabla$ is a symmetric Koszul connections. 
\begin{lem} In every symmetric gauge structure $(M,\nabla)$ the first fundamental equation coincides with the second fundamental equation  $\clubsuit$
\end{lem}
\textbf{HINT}. Since $\nabla$ is symmetric use the formula
$$L_X\nabla = \iota_xR^\nabla + \nabla^2X \clubsuit$$
\begin{prop} \cite{Singer-Stenberg} For every section $\xi$ of a LAS $\mathcal{L}$ the local flows $\phi_\xi(t)$ are $\mathcal{L}$-preserving $\clubsuit$
\end{prop}
\subsubsection{The formalism of Palais}
The formalism of R. Palais is another version of the theory of Cartan-Lie groups. This formalism links the theory of Cartan-Lie groups and the the theory of infinitesimal dynamics of abstract Lie groups. That is what we call Localizations of Abstract Lie Groups in smooth manifolds. We go to provide the complete definition. We consider the Lie algebra $\mathcal{G}$ of a finite dimensional abstract Lie group $G$.
\begin{defn} An infinitesimal dynamic of $G$ in a manifold $M$ is a Lie algebra homomorphism of $\mathcal{G}$ in the Lie algebra of vector fields $\mathcal{x}(M)$, namely
$$\mathcal{G}\ni \xi\rightarrow \tilde{\xi}\in \mathcal{X}(M) \clubsuit$$
\end{defn}
More generally let $\mathcal{G}$ be a finite dimensional Lie algebra acting in a differentiable manifold $M$ via a Lie algebra homomorphism
$$\mathcal{G}\ni \xi\rightarrow \tilde{\xi}\in \mathcal{X}(M).$$
Let $G$ be the simply connected Lie group the Lie algebra of which is $\mathcal{G}.$ The germs of local flows of all $\tilde{\xi}$ generate a pseudogroup of germs of local diffeomorphisms of $M$. That pseudogroup is denoted by $\Gamma_{\mathcal{G}}$. 
\begin{defn} The pseudogroup $\Gamma_\mathcal{G}$ is called a localization of $G$ in $M \clubsuit$
\end{defn}
\begin{defn} \cite{Kumpera-Spencer}, \cite{Malgrange} A pseudogroup of local diffeomorphisms $\Gamma$ is called a Cartan-Lie group if it is 
the complete solution to a Lie equation $\clubsuit$
\end{defn}
The readers interested in global analysis of Lie equations are referred to \cite{Singer-Stenberg}, \cite{Kumpera-Spencer}, \cite{Malgrange}, \cite{Goldschmidt} $\clubsuit$\\
We consider a gauge structure $(M,\nabla)$ is and a finite dimensional Lie algebra $\mathcal{G}.$
\begin{defn} An action of $\mathcal{G}$ in  gauge structure $(M,\nabla)$ is a Lie algebra homomorphism
$$\mathcal{G}\ni \xi\rightarrow \tilde{\xi}\in \mathcal{X}(M)$$
subject to the requirement
$$ L_{\tilde{\xi}}\nabla = 0 \quad \forall \xi\in \mathcal{G} \clubsuit$$
\end{defn} 
\begin{thm} We assume that $(M,\nabla)$ is geodesically complete. We consider a simply connected Lie group $G$ and its Lie algebra $\mathcal{G}$. Every action of $\mathcal{G}$ in $(M,\nabla)$ is the infinitesimal counterpat of a global action of $G$ in $(M,\nabla) \clubsuit $
\end{thm}
\textbf{Hint} By \cite{Kobayashi} every infinitesimal transformation of $(M,\nabla)$ generates a global flow. The conclusion is based on the  formulation of the Lie theory after Richard Palais \cite{Palais}. 
\subsubsection{The application to $\mathcal{J}_\nabla$.}
Now we go involve the fundamental equation $FE^*(\nabla)$ in performing the formalism of Cartan-Lie groups . The purpose is to discuss the global analysis of $FE*(\nabla)$.\\
\begin{thm} For every locally flat manifold $(M,\nabla)$ the first fundamental equation $FE*(\nabla)$ is a transitive and smoothly integrable Lie equation $\clubsuit$
\end{thm}
\textbf{Demonstration.}\\
Let $U$ be the domain of a system of local affine coordinate functions of $(M,\nabla)$, namely $\left\{x_1,..,x_m\right\}$. Let $X$ be a vector field. We put
$$\partial_j = \frac{\partial}{\partial x_j},$$
$$ X = \sum_j X^j \partial_j.$$
Therefore the equation
$$\nabla^2 X = 0$$
is equivalent to the system of partial differential equations
$$\left\{ \frac{\partial^2X^k}{\partial_i\partial_j} = 0\quad \forall i, j, k \in \left\{1,..,m\right\}\right\}.$$
The system above shows that the scalar components of solutions to $FE^*(\nabla)$, namely $X^k$ depend affinely on the coordinate functions $\left\{x_1,..,x_m\right\}$. We consider a Cauchy problem
$$\frac{\partial^2X^k}{\partial_i\partial_j} = 0,$$
$$X^k(x_0) = v^k,$$
$$\frac{\partial X^k(x_0)}{\partial_j} = a^k_j.$$
One easily sees that every Cauchy problem has a unique maximal solution.\\
The theorem is proved $\clubsuit$\\
The theorem above shows that the LAS $\mathcal{J}_\nabla$ is generated by an effective action of a finite dimensional Lie algebra $g_\nabla$. The simply connected Lie group whose Lie algebra is $g_\nabla$ is denoted by $G_\nabla$ \\
The Cartan-Lie group $\Gamma_{g_\nabla}$ is denoted by $\Gamma_\nabla$. It is called the Cartan-Lie group of $(M,\nabla)$.\\
\textbf{Reminder.}\\
For every symmetric gauge structure $(M,\nabla)$ one has
$$ L_X\nabla = \iota_XR^\nabla + \nabla^2X.$$
Therefore if $(M,\nabla)$ is a locally flat structure then 
$$L_X\nabla = 0$$
if and only if
$$\nabla^2X = 0.$$
Thus the Cartan-Lie group $\Gamma_\nabla$ is the transitive Lie pseudogroup of local affine transformations of
$(M,\nabla)$. The LAS of $\Gamma_\nabla$ is the commutator Lie Algebras Sheaf the the Associative Algebras Sheaf $\mathcal{J}_\nabla,\nabla)$. There $\Gamma_\nabla$ admits a bi-invariant affine structure $(\Gamma_\nabla,\tilde{\nabla})$. For convenience we recall that
$$(\gamma_*\nabla)_XY = \gamma_*(\nabla_{(\gamma^{-1}_*(X)}\gamma^{-1}_*(Y)\quad \forall \gamma \in \Gamma_\nabla.$$
The following statements are straightforward consequence of the theorem above.
\begin{thm} The vector space of infinitesimal transformations of a locally flat manifold $(M,\nabla)$ is an associative algebra $J_\nabla$. Further the connection $\nabla$ is induced by the multiplication of $J_\nabla \clubsuit $
\end{thm}
\begin{cor}
The group of affne transformations of a geodesically complete locally flat manifold is a bi-invariant affine Lie group acting transitively.
\end{cor}
For furture usefulness we go to rephrase this theorem. 
\begin{thm} Every locally flat manifold $(M,\nabla)$ is a homogeneous under the action of its Cartan-Lie group $\Gamma_\nabla$. Further if $(M,\nabla)$ is geodesically complete then $(M,\nabla)$ is an affine quotient of a finite dimensional simply connected bi-invariant affine Lie group $(G_\nabla,\nabla^0)$ modulo a connected bi-invariant affine Lie subgroup $(H_\nabla,\nabla^0) \subset (G_\nabla,\nabla^0)$,viz
$$(M,\nabla) = \frac{(G_\nabla,\nabla^0)}{(H_\nabla,\nabla^0)}\clubsuit $$
\end{thm}
\textbf{Hint}.\\
The simply connected abstract Lie group $G_\nabla$ is endowed with a bi-invariant locally flat structure $(G_\nabla,\nabla^0)$. The localization of $G_\nabla$ is the Cartan-Lie group of $(M,\nabla)$. Thus this Cartan-Lie group has a bi-invariant locally flat structure $(\Gamma_\nabla,\tilde{\nabla}^0)$. Furthermore we have already pointed out that the locally flat structure $(M,\nabla)$ is induced by the bi-invariant locally flat structure
$$(\Gamma_\nabla,\tilde{\nabla}^0).$$
We go to more explain this. We consider the effective action
$$g_\nabla \ni u \rightarrow \tilde{u}\in J_\nabla$$
Then $\forall u,v \in g_\nabla$ the vector field
$$\nabla_{\tilde{u}}\tilde{v}$$
is the image in $J_\nabla$ of
$$ \tilde{\nabla}^0_uv \in g_\nabla. $$
This means that the application
$$g_\nabla\ni u \rightarrow \tilde{u}\in J_\nabla$$
is an algebra isomorphism of the associative algebra $(g_\nabla,\tilde{\nabla})$ onto the algebra $(J_\nabla,\nabla)$. SO $(M,\nabla)$ is a homogeneous space of the bi-invariant affine Cartan-Lie group $$(\Gamma_\nabla,\tilde{\nabla}^0).$$
\section{THE FIRST FUNDAMENTAL EQUATION AND NEW INVARIANTS OF THE LOCALLY FLAT GEOMETRY}
\textbf{Reminder.}\\
We have used a locally flat manifold $(M,\nabla)$ for introducing the following new objects.\\
(1.1): The bi-invariant affine Cartan-Lie group $(\Gamma_\nabla,\tilde{\nabla}^0)$.\\
(1.2): The simply connected bi-invariant affine (abstract) Lie group $(G_\nabla,\nabla^0)$.\\
(1.3): The localization of $(G_\nabla,\nabla^0)$
$$(G_\nabla,\nabla^0)\rightarrow (\Gamma_\nabla,\tilde{\nabla}^0) \clubsuit $$
We go to repeat this process by replacing $(M,\nabla)$ by $(G_\nabla,\nabla^0)$. Thereby we obtain a new transitive Associative Algebras Sheaf in $G_\nabla^0$, namely $\mathcal{J}_{\nabla^0}$ which yields new data\\
(2.1): The bi-invariant affine Cartan-Lie group $(\Gamma_{\nabla^0},\tilde{\nabla}^1)$.\\
(2.2): The simply connected bi-invariant affine (abstract) Lie group $(G_{\nabla^0},\nabla^1).$\\
(2.3): The localization of $(G_{\nabla^0},\nabla^1)$
$$ (G_{\nabla^0},\nabla^1)\rightarrow (\Gamma_{\nabla^0},\tilde{\nabla}^1) \clubsuit $$
\subsection{Inductive-projective systems of bi-invariant affine Lie groups}
Given a locally flat gauge structure $(M,nabla)$ we just introduced the following data\\
(i) the gauge structure $(G_{\nabla^0},\nabla^1)$ is homogeneous under the action of $(\Gamma_{\nabla^0},\tilde{\nabla}^1)$,\\
(ii) $(J_\nabla,\nabla)$ is a subalgebra of the associative algebra $(J_{\nabla^0},\nabla^1)$. Thereby there is a canonical 
affine Lie group immersion 
$$ I^0_1:=\quad (G_\nabla,\nabla^0) \rightarrow (G_{\nabla^0}, \nabla^1).$$
If the gauge structure $(G_{\nabla},\nabla^0$ is complete then we get a canonical projection
$$\pi^1_0:=\quad (G_{\nabla^0},\nabla^1)\rightarrow (G_\nabla,\nabla^0)  $$
We go to repeat this construction through the level $q$. There we get new data\\
(q.1): the bi-invariant affine Cartan-Lie group $(\Gamma_{\nabla^q},\tilde{\nabla}^{q+1})$,\\
(q.2): the simply connected bi-invariant affine (abstract) Lie group $(G_{\nabla^q},\nabla^{q+1}),$\\
(q.3): the localization of $(G_{\nabla^q},\nabla^{q+1})$
$$(G_{\nabla^q},\nabla^{q+1})\rightarrow (\Gamma_{\nabla^q},\tilde{\nabla}^{q+1}). $$
\textbf{Fact} The associative algebra $J_{\nabla^q}$ is subalgebra of the associative algebre $J_{\nabla^{q+1}}$. Therefore we have the canonical immersion 
$$I^q_{q+1}:\quad (G_{\nabla^q},\nabla^{q+1})\rightarrow (G_{\nabla^{q+1}},\nabla^{q+2}). $$
Let assume that the gauge structure $(G_{\nabla^q},\nabla^{q+1})$ is geodesically complete then  the associative algebra 
$J_{\nabla^{q+1}}$ is the infinitesimal counterpart of a transitive action of $(G_{\nabla^{q+1}},\nabla^{q+2})$ in $(G_{\nabla^q},\nabla^{q+1})$. Thereby we also get the canonical projection
$$\pi^{q+1}_q: (G_{\nabla^{q+1}},\nabla^{q+2})\rightarrow (G_{\nabla^q},\nabla^{q+1})$$ 
 Henceforth we assume that all of the bi-invariant gauge structures $ (G_{\nabla^q};\nabla^{q+1})$ are geodesically complete. We consider two non negative integers $p,q$ with
$$ p < q.$$
We define the projection
$$ \pi^q_p:(G_{\nabla^q},\nabla^{q+1})\rightarrow (G_{\nabla^p},\nabla^{p+1})$$
by setting
$$ \pi^q_p = \pi^{p+1}_p\circ...\circ \pi^q_{q-1}.$$
However we the canonical immersion
$$ (G_{\nabla^p},\nabla^{p+1})\rightarrow (G_{\nabla^q},\nabla^{q+1})$$
by setting
$$ I^p_q = I^{q-1}_q\circ..\circ I^p_{p+1}.$$
Let $p, q, r $ be positive integers such that 
$$p < q < r. $$
The canonical maps we just introduced have the following properties\\
$$ \pi^q_p\circ\pi^r_q = \pi^r_p,$$
$$  I^p_q\circ I^q_r = I^p_r,$$
$$ \pi^q_q = I^q_q = 1.$$
Those canonical projection form a \textbf{semi-projective system} of finite dimensional simply connected bi-invariant affine Lie groups 
$$ \mathcal{PS}_\nabla:\quad = \left\{\pi^q_p:\quad(G_{\nabla^q},\nabla^{q+1})\rightarrow(G_{\nabla^p},\nabla^{p+1})\right\}.$$
\textbf{Warning}. The canonical projections are not group homomorphisms. That is the reason we say \textbf{semi-projective system} of bi-invariant affine Lie groups.\\   
By replacing the canonical projections by the canonical immersions we obtain a \textbf{sem-inductive system} of simply connected finite dimensional bi-invariant affine Lie groups
$$\mathcal{IS}^\nabla:\quad = \left\{I^p_q:\quad(G_{\nabla^p},\nabla^{p+1})\rightarrow(G_{\nabla^q},\nabla^{q+1}) \right\}.$$
\textbf{Warning}. The canonical immersion may not be one to one. That is the reason we say \text{semi-inductive system} of bi-invariant affine Lie groups.
From the infinitesimal viewpoint we obtain the inductive system of associative algebras 
$$\left\{I^p_q:\quad(J_{\nabla^p},\nabla^{p+1})\rightarrow (J_{\nabla^q},\nabla^{q+1})\right\}$$
\begin{defn} Starting with a locally flat manifold $(M,\nabla)$ we have constructed the pairs 
$\left\{(M,\nabla),\mathcal{IS}^\nabla\right\}$, and $\left\{(M,\nabla),\mathcal{PS}_\nabla\right\}.$ 
They are called a \text{semi-inductive system} of simply connected bi-invariant affine Lie groups and \textbf{semi projective system} of simply connected bi-invariant affine Lie groups.
\end{defn}
\subsection{Initial objects and Finala objects: the structural theorem of the locally flat geometry}
\textbf{Reminder.}\\
For convenient we adopt the following notation,
$$ (G_q,\nabla^q) = (G_{\nabla^q},\nabla^{q+1})$$
Given a locally flat manifold $(M,\nabla)$ we have introduced\\
(1) the semi inductive system of simply connected bi-invariant affine Lie groups $\left\{(M,\nabla),\mathcal{IS}^\nabla\right\}$,\\
(2) the semi-projective system of simply connected bi-invariant affine Lie groups $\left\{(M,\nabla,\mathcal{PS}_\nabla\right\}$,\\
(3) the localizations of those systems  systems of bi-invariant affine Cartan-Lie groups
$$ \left\{ (\Gamma_{q+1},\tilde{\nabla}^{q+1})\times (G_q,\nabla^q)\rightarrow (G_q,\nabla^q) \right\} $$
Thereby we get new insights of the locally flat geometry. Those links between the locally flat geometry and the theory of Lie is highlighted by a structural theorem.
\begin{thm} The structural theorem.\\
(a) Every finite dimensional geodesic complete locally flat manifold is the Final Object of an optimal semi-projective system of finite dimensional simply connected bi-invariant affine Lie groups.\\ 
(b) Every finite dimensional affine Lie group is the Initial Object of an optimal semi-inductive system of finite dimensional simply connected bi-invariant affine Lie groups.\\
(c) Every finite dimensional affine Lie group is the Final Object of an optimal semi-projective system of finite dimensional simply connected bi-invariant affine Lie groups $\clubsuit$
\end{thm}
\textbf{Reminder}\\
The systems of bi-invariant affine Cartan-Lie groups will be used for exploring the moduli space of finite dimensional simply connected locally flat manifolds. \\
For non specialists we recall that a transitive Cartan-Lie group is the same thing as a transtive Lie Algebras Sheaf (LAS) and that a transitive bi-invariant affine Cartan-Lie group is the same thing as a transitive Associative Algebras Sheaf (AAS). Further a AAS $\mathcal{J}$ is uniquely defined by $J^\infty(\mathcal{J})$. Loosely speaking if a vectr field $X$ satistfies 
$$j^\infty(X) \in J^\infty(\mathcal{J})$$
then $X$ is a section of $\mathcal{J}$. We observe that the infinitesimal version of (b) is useful for exploring the topology of finite dimensional Lie groups. We go to state some algebraic counterparts of the structural theorem.
\begin{thm} The structural theorem versus Algebras Sheaf.\\
(1) In a finite dimensional manifold every finite dimensional transitive Koszul-Vinberg Algebras Sheaf is the Initial Object of an optimal semi-inductive system of transitive Associative Algebras Sheaves.\\ 
(2) In a finite dimensional Lie group every finite dimensional transitive left invariant Koszul-Vinberg Algenbras Sheaf is a Final Object of a semi-projective system of finite dimensional Associative Algebras Sheaves $\clubsuit$
\end{thm}
\section{THE FIRST FUNDAMENATAL EQUATION OF A GAUGE STRUCTURE AND THE FOLIATIONS, continued}
We go to deal with the differential equation $FE^*(\nabla)$. The Associative Algebras Sheaf $(\mathcal{J}_\nabla,\nabla)$ is a geometric invariant of the gauge structure $(M,\nabla)$. In this section we go to focus on the subset of Koszul connections whose first fundamental equation is smoothly integrable. 
\subsection{The special gauge structures}
\begin{defn} A gauge structure $(M,\nabla)$ is called special if
 $$J_\nabla \neq 0\clubsuit$$
\end{defn}
\textbf{Examples}.\\
(1) Every locally flat structure $(M,\nabla)$ is a special gauge structure.\\
(2) For every positive integer $m$ the open unit ball $B^m$ equipped with the Euclidean Levi-Civita connection is a spacial gauge structure.\\
(3) In the real analytic category all of the finite dimensional gauge structures are special. In particuliar metric connections of Kaehlerian manifolds are specail.\\
The pair $(\mathcal{J}_\nabla, [-,-]_\nabla )$ is a sheaf of Lie algebras whose bracket is defined by
$$[X,Y]_\nabla = \nabla_XY - \nabla_YX. $$
\textbf{Warning. The sheaf of Lie algebras $(\mathcal{J}_\nabla, [-,-]_\nabla)$ is a LAS if and only if $(M,\nabla)$ is torsion free.}\\
Whatever the sheaf of Lie algebras $(\mathcal{J}_\nabla, [-,-]_\nabla)$ yields the family of finite dimensional bi-invariant affine Lie group 
$\left\{(G_q,\nabla^q)\right\}$. Under the asumption that all of the simply connected bi-invariant affine Lie groups $(G_q,\nabla^q)$ are geodesically complete we get the systems\\
$$\mathcal{PS}_\nabla: =\quad\left\{(G_q,\nabla^q)\rightarrow(G_p,\nabla^p)\right\},$$
$$\mathcal{IS}^\nabla: =\quad\left\{(G_p,\nabla^p)\rightarrow(G_q,\nabla^q)\right\},$$
and the corresponding projective system of bi-invariant affine Cartan-Lie groups
$$ \left\{(\Gamma_{q+1},\tilde{\nabla}^{q+1}) \times (G_q,\nabla^q) \rightarrow(G_q,\nabla^q)\right\}.$$
We recall that the framework of our structural theorems is the category of finite dimensional special gauge structures. However if $(M,\nabla)$ is non symmetric the links between those systems of bi-invariant affine Lie groups and the differential topology of $M$ are unclear.
\subsection{The topological nature of the semi-inductive systems of bi-invariant affine Lie groups}
\textbf{Warning: geometric data which are invariant by the gauge group $\mathit{G}(M)$ are called gauge invariants of the manifold $M$; those which are invariant by the group of diffeomorphisms $Diff(M)$ are called geometric invariants of the manifold $M$}.\\
Gauge invariants of a gauge structure $(M,\nabla)$ are data which are invariant by $J_{\nabla\nabla}$. Geometric invariants of $(M,\nabla)$ are data which are invariant by the infinitesimal transformations of $(M,\nabla)$. For instance the semi-projective system
$$\mathcal{PS}_\nabla: =\quad \left\{(G_q,\nabla^q)\rightarrow(G_p,\nabla^p) \right\}$$
and the semi-inductive system
$$\mathcal{IS}^\nabla: =\quad\left\{(G_p,\nabla^p)\rightarrow(G_q,\nabla^q) \right\}$$
are geometry invariants of $(M,\nabla)$, Loosely speaking every $\nabla$-preserving diffeomorphism of $M$ preserves those systems. Really it is easy to see that every $(\phi\in Diff(M,\nabla)$ gives rise to an automorphism $\phi_q$ of the affine Lie group $(G_q,\nabla^q)$.\\
We know that in a symmetric gauge structure $(M,\nabla)$ the following assertions are equivalent\\
(1): The Lie equation $FE^{*}(\nabla)$ is transitive,\\
(2): The couple $(M,\nabla)$ is a locally flat manifold.
We keep the notation used in the preceding sections. Assume that $(M,\nabla)$ is a special symmetric gauge structure and assume that the LAS 
$\mathcal{J}_\nabla$ is regular, viz the dimension of $\mathcal{T}_\nabla M(x) \subset T_xM$ does not depend on $x \in M$. Then $J_\nabla$ defines a foliations whose leaves are locally flat manifolds. Every leaf of this foliation is the Final Object of the semi-projective system simply connected bi-invariant affinely flat Lie groups $$\left\{(G_q,\nabla^{q+1}), \pi^q_p | p \leq q\right\}. $$ This topological asppect show that our machineries yiels a solution to $EXF(\mathcal{S}: FL(M))$. We plan going into this perspective in next sections. However go to conclude this ection by
\begin{thm} Every regular gauge structure $(M,\nabla)$ support the locally flat foliation defined by $J_\nabla$ every leaf of which is the Initial object of a semi-inductive system of finite dimensional simply connected bi-invariant affine Lie groups $\left\{(G_q,\nabla^q)\right\}$. \\
Further if all of the gauge structures $(G_q,\nabla^q)$ are geodesically complete then every leaf of $J_\nabla$ is the Final object of semi-projective system of finite dimensional simply connected bi-invariant affine groups $\clubsuit$ 
\end{thm}
The theorem just stated is of great interest in the information geometry 
\section{NEW NUMERICAL INVARIANTS OF GAUGE STRUCTURES}
The moduli space of gauge structures in $M$ is the orbit space 
$$\frac{\mathcal{MG}(M)}{\mathit{G}(M)}.$$
The moduli space of Koszul connections in $M$ is the orbit space
$$ \frac{\mathcal{LC}(M)}{Diff(M)}.$$
In this section we go to use the fundamental equation $FE^*(\nabla)$ for introducing a $\mathbb{Z}$-valued function $r^b$ which is defined in the moduli space 
$$ \mathbb{LC}(M)= \frac{\mathcal{LC}(M)}{Diff(M)}.$$
In this sections we go to focus on two concerns.\\
(i) The first concern is to investigate the links between the function $r^b$ to be introduced and the open problems we have raised in the introduction, namely the problems $EX(\mathcal{S})$, $EXF(\mathcal{S})$ and $DL(\mathcal{S})$.\\
(ii) Another concern is to well understand the instrinsic nature of the function $r^b$. In particular we are interested in its link with the differential topological. The objects of the differential topology we are interested in are intrinsic geometries of foliations, viz the problem  $EXF(\mathcal{S})$.\\
Let $(M,\nabla)$ be a special gauge structure and let $(\mathcal{J}_\nabla, \nabla)$ be its Associative Algebras Sheaf. For convenience we recall the inclsion maps
$$\mathcal{LF}(M)\subset \mathcal{SLC}(M)\subset \mathcal{LC}(M).$$
\subsection{Numerical invariants of the first fundamental equation}
In the category of gauge structures in a manifold $M$ the group of diffeoemorphisms act by the group of 1-jets 
$$J^1(Diff(M)) = \left\{ (\phi,\phi_*)| \phi \in Diff(M)\right\},$$ 
here $\phi^*$ is the differential of the diffeomorphism $\phi$. The gauge group $\mathit{G}(M)$ is a normal subgroup of $J^1(Diff(M))$. The group $J^1(Diff(M))$ acts in $\mathcal{MG}(M)$. Really $\frac{\mathcal{LC}(M)}{Diff(M)}$ stands for 
$\frac{\mathcal{MG}(M)}{J^1(Diff(M))}$. The vector space of sections of the vector bundle $Hom(TM\otimes TM,TM)$ is denoted by $T^2_1(M)$. The vector space  of sections of $Hom(S^2(TM),TM)$ is denoted by $S^2_1(M)$
 \textbf{Warning: Depending on frameworks and needs the same object may be denoted differently}. We recall some structures we go to focus on.\\
\subsubsection{Affine Geometry}
$A1: \mathcal{LC}(M)$ is an affine space whose vector space of translations is $T^2_1(M)$,\\
$A2: \mathcal{SLC}(M)$ is an affine subspace whose vector space of translations is $S^2_1(M)$.\\
\subsubsection{Combinatory}
$GR_1: \mathcal{GR}(LC)$ stands for the graph $\left\{\mathcal{MG},T^2_1(M)\right\}$.\\
$GR_2: \mathcal{GR}(SLC)$ stands for the graph $\left\{\mathcal{SLC}(M), S^2_1(M)\right\}$\\
$GR_3: \mathcal{GR}(LF)$ stands for the graph $\left\{\mathcal{LF}, \mathbb{MC}(LF)\right\}$. \\
Here $\mathbb{MC}(LF)$ stands for the zeros of Maurer-Cartan polynomial of locally flat structures \cite{Nguiffo Boyom(3)}. In fact let $\nabla, \nabla^* \in \mathcal{LF}(M)$ and let 
$$ \tau = \nabla^* - \nabla \in S^2_1(M)$$
Then $\tau$ is a zero of Maurer-Cartan polynomial of both $\nabla$ and $\nabla^*$.\\
\subsubsection{Moduli spaces}
The moduli space we are interested in\\
$M_0$: the moduli space of gauge structures
$$\frac{\mathcal{MG}(M)}{\mathit{G}(M)},$$
$M_1$: the moduli space of Koszul connections
$$\mathbb{LC}(M) = \frac{\mathcal{LC}(M)}{Diff(M)},$$
$M_2$: the moduli space of symmetric Koszul connections
$$\mathbb{SLC}(M) = \frac{\mathcal{SLC}(M)}{Diff(M)},$$
$M_3$: the moduli space locally flat Koszul connections
$$\mathbb{LF}(M) = \frac{\mathcal{LF}(M)}{Diff(M)}.$$
Let $\phi\in Diff(M)$, $\psi \in \mathit{G}(M)$ and $\nabla \in \mathcal{LC}(M)$. Given a vector field $X$ the covariant derivatives $[\phi_*(\nabla)]_X$ and $[\psi(\nabla)]_X$ are defined by
$$ [\phi_*(\nabla)]_X =
\phi^{-1}_*\circ\nabla_{\phi^{-1}_*(X)} \circ \phi^{-1}_* ,$$
$$ [\psi(\nabla)]_X = \psi\circ \nabla_X \circ \psi^{-1}.$$
In both $\frac{\mathcal{MG}(M)}{\mathit{G}(M)}$ and $\frac{\mathcal{LC}(M)}{Diff(M)}$ the orbit of $\nabla)$ is denoted by $[\nabla]$.
\begin{defn} The function
$$\mathbb{LC}(M)\ni [\nabla]\rightarrow r^b([\nabla])\in \mathbb{Z}$$
is defined by
$$ r^b([\nabla]) = \min_{\left\{x \in M\right\}}\left\{dim[T_\nabla M(x)\right\} \clubsuit $$
\end{defn}
If $\mathcal{T}$ is a $Diff(M)$-invariant subset of $\mathcal{LC}(M)$ then we restrict the function $r^b$ to
$$\mathbb{T} = \frac{\mathcal{T}}{Diff(M)}.$$
We go to point emphasize the roles played by some restrictions of the function $r^b$. Those undreamed roles help to address the problems $EX(\mathcal{S})$, $EXF(\mathcal{S})$ and $DL(\mathcal{S})$.\\
In an $m$-dimensional manifold $M$ we have $$ 0 \leq r^b([\nabla])\leq m.$$
\textbf{A comment} A finite dimensional manifold has many numerical geometric invariants such as its dimension, its Betti numbers, its Euler characteric. We go to use the function $r^b$ for enriching every finite dimensional manifold with new numerical invariants
\begin{defn} In the moduli space $\mathbb{LC}(M)$ we use the function $r^b$ for introducing the following non negative integers \\
(1): $r^b_\tau(M) = \min_{\left\{\nabla \in \mathbb{LC}(M)\right\}}\left\{dim(M) -
r^b([\nabla])\right\},$\\
(2): $r^b(M) = \min_{\left\{\nabla\in \mathbb{SLC}(M)\right\}}\left\{dim(M) - r^b([\nabla])\right\}\clubsuit$
\end{defn}
We have the inequality
$$ r^b_\tau(M)\leq r^b(M).$$
Both $r^b(M)$ and $r^b_\tau(M)$ are global geometric invariants of $M$. We go to investigate the nature of those new invariants. We focus on insights which are related to the challenges we go to face, namely $EX(\mathcal{S})$, $EXF(\mathcal{S})$ and $DL(\mathcal{S})$.
\subsection{New insights}
\textbf{(I) The theory of obstructions}. As said before we aim at involving $r^b$ in studying $EX(\mathcal{S})$.\\ 
\textbf{(II) Combinatory}. We aim at involving $r^b$ in studying $DL(\mathcal{S})$.\\
\textbf{III The differential topology}. We aim at involving $r^b$ in studying $EXF(\mathcal{S})$.\\
\subsection{Some relative gauge invariants}
A geometric structure $\mathcal{S}$ may be well handled if it admits a unique specific Koszul connection $\nabla_\mathcal{S}$. The uniqueness of $\nabla_\mathcal{S}$ may be source of crucial geometric invariants of $\mathcal{S}$. That idea goes back to Elie Cartan \cite{Cartan(1)}, \cite{Cartan(2)}, \cite{Cartan(3)}.\\
(1) In the naïve Riemannian geometry the torsion free metric connection is unique. The curvature tensor (of the Levi-Civita conection) and its avatars are source of the major part of the classical Riemannian geometry.\\
(2) In the bi-lagrangian geometry the connection of Hess is the alternative to the Levi-Civita connection in the Riemannian geometry. It is the unique torsion free connection preserving the bi-lagrangian structure, see \cite{Nguiffo Boyom(7)} and references therein.\\
There are few known geometric structures with a unique restricted connection. As announced we go to use the function $r^b$ for introducing new relative numerical geometric invariants which intersting nature. 
\subsection{The notion of Hessian defect}
For our purpose we shall be recalling some well known definitions.
\begin{defn}
A Hessian manifold is a triple $(M,g,D)$ formed of a Riemannian structure $(M,g)$ which is linked with a locally flat structure $(M,D)$ by
$$(D_Xg)(Y,Z) - (D_Yg)(X,Z) = 0\clubsuit $$
\end{defn}
The set of Hessian structures in a manifold $M$ is denoted by $\mathcal{H}es(M)$. 
From the definition of Hessian structure arise two subsequent problems.\\
$EX(\mathcal{S}: Hes(M,g))$: the question whether there exist Hessian structures in a Riemannian manifold $(M,g)$, ( S-I Amari, A.K. Guts)\\ 
$EX(\mathcal{S}: Hes(M,D))$: the question whether there exist Hessian structures in a locally flat manifold $(M,D)$, J-L Koszul, J. Vey.\\
In the same contexts arise the problems $EXF(\mathcal{S})$ and $DL(\mathcal{S}).$
\begin{defn} The solution to $DL(\mathcal{S}: Hes(M,g))$ is called the Riemmannian Hessian defect.\\
The solution to $DL(\mathcal{S}: Hes(M,D))$ is called the affine Hessian defect $\clubsuit$
\end{defn}
To study those Hessian defects we perform an arsenal formed of the function $r^b$ and the functor of Amari.
\subsection{Te functor of Amari}
\begin{defn} The dualistic relation of Amari is the functor
$$\mathcal{R}ie(M)\times \mathcal{LC}(M)\ni (g,\nabla)\rightarrow \nabla^*\in \mathcal{LC}(M).$$
Here the Koszul connection $\nabla^*)$ is defined  by
$$  g(\nabla^*_XZ,Y) = Xg(Y,Z) - g(\nabla_XY,Z) \clubsuit$$
\end{defn}
We use the dualistic relation for introducing two new functor.\\
\subsubsection{The first functor}
We deal with the restriction
$$\left\{g\right\}\times \mathcal{LF}(M)\ni (g,\nabla)\rightarrow \nabla^*\in \mathcal{LC}(M).$$
Here $\nabla^*$ is defined as above, viz
$$g(\nabla^*_XY,Z) = Xg(Y,Z) - g(Y,\nabla_XZ).$$
\subsubsection{The second functor}
We deal with the restriction
$$\mathcal{R}ie(M)\times\left\{\nabla\right\} \ni (g,\nabla)\rightarrow \nabla^g \in \mathcal{LC}(M).$$
Here $\nabla^g$ is defined by
$$g(\nabla^g_XY,Z) = g(Y,\nabla_XZ) - g(Y,\nabla_XZ).$$
The functors we just defined are called functors of Amari. Our arsenal is composed of the functors  of Amari and the function $r^b$.\\
The first functor is used for introducing the following numerical constants\\
$(1):\quad r^b(M,g) = \min_{\left\{\nabla]\in\mathbb{SLC}(M)\right\}}\left\{dim(M) - r^b([\nabla^*])\right\}$,\\
$(2):\quad r^B(M) = \min_{\left\{g\in \mathcal{R}ie(M)\right\}}\left\{r^b(M,g)\right\}.$\\
\begin{defn} Given a Riemannian manifold $(M,g)$\\
$(1):$ The integer $r^b(M,g)$ is called the relative Riemannian Hessian defect of $(M,g)$.\\
$(2):$ The integer $r^B(M)$ is called the absolute Riemanncian Hessian defect of $M$ $\clubsuit$
\end{defn}
The second functor is used for introducing the following numerical invariants\\
$(3):\quad r^b(M,D) = \min_{\left\{g\in \mathcal{R}ie(M)\right\}}\left\{dim(M) - r^b(D^g)\right\}.$\\
$(4):\quad r^{B^*}(M) = \min_{\left\{D\in \mathcal{LF}(M)\right\}}\left\{r^b(M,D)\right\}.$\\
\begin{defn} Given a locally flat manifold $(M,D)$\\
$(3*):$ The integer $r^b(M,D)$ is called the relative affine Hessian defect of $(M,D)$.\\
$(4*):$ The integer $r^{B*}(M)$ is called the absolute affine Hessian defect of $M$ $\clubsuit$.
\end{defn}
\section{THE FIRST FUNDAMENTAL EQUATION AND THE INFORMATION GEOMETRY}
The information geometry is an outstanding domain of application of the Hessian geometry, \cite{Armstrong-Amari}, \cite{Armstrong-Amari}, \cite{Nguiffo Boyom(6)}, \cite{Nguiffo Boyom-Wolak(3)}. Among the importants problems in the theory of statistical models is the question explicitly raised by Murray-Rice \cite{Murray-Rice}. \textbf{ When is exponential a statistical model?} That is the information geometry versus of the question raised by Alexander K. Guts in the context of Physics \cite{Guts}. The same question is linked with the geoemtric theory of heat \cite{Barbaresco}. We just introduced a numerical invariant named Riemannian Hessian defect. The Riemannian Hessian defect is useful for studying this problem of complexity of statistical models, viz the question  whether a statistical model is isomorphic to an exponential model. In this section we restrict the attention to the case of regular statistical models. The Fisher information of a regular model is a Riemannian metric. By \cite{Nguiffo Boyom(6)} the complexity of regular statistical models is nothinng but the  problem $DL(\mathcal{S}: He(M,g))$. Before focusing on the impacts of $r^b$ on the information geometry we go to recall some basic notions. Let $(\Xi,\Omega)$ be a measurable set. We assume that $\Xi$ is homogeneous under the natural action of the group of measurable transformations of $(\Xi,\Omega)$, namely
$$\Gamma = Aut(\Xi,\Omega).$$
\begin{defn} ( see the appendix to this paper.) A $m$-dimensional statistical model for $(\Xi,\Omega)$ is a quintuple
$$ \mathbb{M} = [\mathcal{E},\pi,M,D,p].$$
Here\\
(1): $(M,D)$ is a locally flat manifold.\\
(2): The triple $(\mathcal{E},M,\mathbb{R}^m)$ is endowed with an action of the group $\Gamma$
$$\Gamma\times \mathcal{E}\times M\times \mathbb{R}^m \ni (\gamma,e,x,t)\rightarrow (\gamma.e,\gamma.,\gamma.t)\in \mathcal{E}\times M\times \mathbb{R}^m.$$
(3) The projection $$\pi: \mathcal{E}\ni e\rightarrow \pi(e) \in M$$
is a $\Gamma$-equivariant,viz 
$$\pi(\gamma.e) = \gamma.\pi(e),$$
further this fibration $\pi$ is a locally trivial bundle whose fibers are isomorphic to $\Xi$
(4): The real valued function $p$ is defined in $\mathcal{E}$, it is horizontally differentiable and its restriction to every $\pi$-fiber is a probability density.
(5): The horizontal differentiation is denoted by $\frac{\partial}{\partial \theta}$, the integration along the $\pi$-fiber is denoted by $\int_\Xi$ and we have
$$\frac{\partial}{\partial \theta}\circ \int_\Xi = \int_\Xi\circ \frac{\partial}{\partial \theta}\clubsuit $$
\end{defn}
For more details on the whole theory the readers are referred to \cite{Nguiffo Boyom(6)}
We consider the Hessian differential operator $D^2$. According to condition (5) the operator
$$D^2\circ\int_\Xi$$ is well defined.\\
\begin{defn} A statistical model
$$\mathbb{M} = [\mathcal{E},\pi,M,D,p]$$
is called an exponential family if there exist two functions
$$\mathcal{E}\in \rightarrow a(e)\in\mathbb{R},$$
$$ M \in x \rightarrow \psi(x)\in \mathbb{R}$$
subject to the following requirement.\\
$(1): p(e) = exp(a(e) - \psi(\pi(e))$,\\
$(2): D^2[\int_{\Xi} a(e)] = 0\clubsuit$
\end{defn}
\subsection{The Fisher information}
Among central ingredients of the Information Geometry is the Fisher information. We recall its definition
\begin{defn} The Fisher information of a model
$$\mathbb{M} = [\mathcal{E},\pi,M,D,p]$$
is defined by
$$g_\mathbb{M}(x) = - \int_{\mathcal{E}_x} p(e)D^2 log(p(e))$$
The model $\mathbb{M}$ is called regular if the Fisher information is positive definite
$\clubsuit$
\end{defn}
\begin{defn} Given a regular statistical model $\mathbb{M}$ the invariant $r^b(M,g_\mathbb{M})$ is called the Riemannian exponential-defect of $\mathbb{M} \clubsuit$
\end{defn}
\subsection{The connections of Chentsov}
We go to define the Christoffel symbols of the family of Koszul connections $\left\{\nabla^\alpha; \alpha\in\mathbb{R}\right\}.$ \\
Let $(U,\Phi\times\phi $ be a local chart of the fibration of $[\mathcal{E},\pi,M,D]$.\\
For every $e \in \mathcal{E}_U$ we set
$$\Phi(e) = (\theta(e),\xi(e)),$$
$$\phi(\pi(e)) = \theta(e).$$
Here
$$\theta(e) = (\theta_1,..,\theta_m) \in \mathbb{R}^m.$$
At every $x \in U$ we put
 $$[\nabla^\alpha_{\frac{\partial}{\partial\theta i}}\frac{\partial}{\partial \theta j}](x) = \sum_k[\Gamma^k_{ij:\alpha}](x)\frac{\partial}{\partial k}.$$
Here
$$[\Gamma^k_{ij:\alpha}](x) = \int_\mathcal{E}x p(e)[\frac{\partial^2log(p(e))}{\partial \theta_i\partial \theta_j} + \frac{1+\alpha}{2}\frac{\partial log(p(e))}{\partial \theta_i}\frac{\partial log(p(e))}{\partial \theta_j}][\frac{\partial log(p(e))}{\partial \theta_k}]]d\xi(e)$$
The definition of the family $\nabla^\alpha$ agrees with affine coordinate change.\\
We consider the family of fundamental equations 
$$\left\{FE^*(\nabla^\alpha)| \alpha \in \mathbb{R}\right\}$$ 
Thereby we obtain the family of the functions  
$$r^b(\alpha) = r^b(\nabla^\alpha).$$
\begin{defn} The Chentsov exponential defect of a model $\mathbb{M}$ is defined by
$$r^b(\mathbb{M}) = \min_{\alpha \in \mathbb{R}}\left\{dim(M) - r^b(\alpha)\right\}\clubsuit$$
\end{defn}
Actually we put
$$ r^b(\mathbb{M}) \leq r^b(M,g)$$
Those integers are statistical invariants of the measurable set $(\Xi,\Omega)$.
\subsection{The moduli space of Cartan-Lie groups of simply connected locally flat manifolds}
We are interested in the problem of local equivalence of Lie pseudogroups. The pioneering works on the object are those of Sophus Lie and  Elie Cartan \cite{Singer-Sternberg} and references therein. The major modern step has been the new formalism of Spencer-Kodaira-Kuranishi \cite{Spencer}, \cite{Singer-Sternberg}, \cite{Kuranishi}, \cite{Kodaira}. See also \cite{Goldschmidt}, \cite{Kumpera-Spencer}, \cite{Malgrange}. \\
In this paper we use a contribution due to V. Guillemin \cite{Guillemin}. After \cite{Guillemin} the problem of equivalence of transitive Cartan Lie groups of finite type is reduced to the problem of the formal equivalence of their LAS. Loosely speaking, for those groups formaly equivalence implies differentiable equivalence. \\
We aim at performing the formalism of Guillemin in studying the moduli space of Cartan-Lie groups of locally flat manifolds. We go to show that in our context the study is easy because the Cartan-Lie groups we go to deal with are localizations of finite dimensional simply connected (abstract) Lie groups.
Another challenge we go to face is the problem of moduli space of finite dimensional geodesically complete locally flat manifolds. For this purpose we aim at involving the semi-inductive systems of bi-invariant affine Cartan-Lie groups that have been introduced, namely
 $$\mathcal{CLG}^\nabla: = \left\{(\Gamma_{p+1},\tilde{\nabla}^{p+1}) \times(G_p,\nabla^p)\rightarrow (G_p,\nabla^p) \right\}.$$
Those systems are the localizations of the semi-inductive systems of finite dimensional simply connected bi-invariant affinely flat Lie groups
$$\mathcal{IS}: =\quad \left\{(G_p,\nabla^p)\rightarrow (G_q,\nabla^q)\right\}.$$
\textbf{Reminder}.\\
In the preceding sections the moduli space of locally flat structures in $M$ has been denoted by
$$ \mathbb{LF}(M) = \frac{\mathcal{LF}(M)}{Diff(M)}.$$
We go to implement the invariants we have introduced. Our approach is similar to an approach already used by E.B. Vinberg \cite{Vinberg(1)}, see also Y. Matsushima \cite{Matsushima}
\subsubsection{Notation-Definition, continued}
Before stating some important new theorems we introduce a new algebraic version of the Cartan-Lie Theory of homogeneous spaces. We deal the infintesimal theory of finite dimensional Cartan-Lie groups which is nothingelse than the theory of LAS of finite type \cite{Palais}, \cite{Guillemin}, \cite{Guillemin-Sternberg}, \cite{Singer-Sternberg}.\\
We have already seen that every finite dimensional locally flat manifold $(M,\nabla)$ is the Final Object of a semi-projective system of finite dimensional simply connected bi-invariant affine Cartan-Lie groups. Subsequently every finite dimensional geodesically complete locally flat manifold is the Final Object of a semi-projective system of finite dimensional bi-invariant affine Lie groups.
\subsubsection{The notion of effective pair}
To study the problem of the moduli space we go to introduce some new ingredients.\\
\begin{defn} Let $\mathcal{J}$ be an associative algebra.\\
$(1):$ A right ideal $\mathcal{I}$ is called simple if it does not contain any non zero two-sided ideal of $\mathcal{J}$\\
$(2):$ Consider a pair $\left\{\mathcal{I} < \mathcal{J}\right\}$ where $\mathcal{I}$ is a simple right ideal of an associative algebra $\mathcal{J}$; that pair is called a pair of associative algebras $\clubsuit$.
\end{defn}
Let $\left\{\mathcal{I} < \mathcal{J}\right\}$ be a finite dimensional pair of associative algebras. Let $n$ be the dimension of $\mathcal{I}$. The grassmannian of $n$-dimensional vector subspaces in $\mathcal{J}$ is denoted by $G_n(\mathcal{J})$. Let $G_J$ be the simply connected abstract Lie group whose Lie algebra is the commutator Lie algebra of $\mathcal{J}$ and let
$$[\mathcal{I}] \subset G_n(\mathcal{J})$$
be the orbit of $\mathcal{I}$ under the adjoint representation of $G_J$ in $\mathcal{J}$. We go to use the notion of effective pair.
\begin{defn} Given a pair of associative algebras $\left\{\mathcal{I} < \mathcal{J}\right\}$ the pair
$\left\{[\mathcal{I}] < G_n(\mathcal{J})\right\}$ is called an effective pair of associative algebras $\clubsuit$
\end{defn}
\begin{defn} We keep the notation just settled.\\
(1): The category whose objects are effective pairs $\left\{[\mathcal{I}] < G_n(\mathcal{J})\right\}$ is denoted by
$$\left\{[\mathcal{I}] < G_n(\mathcal{J})\right\}(ASSA)$$
(2): A morphism of $\left\{[\mathcal{I}] < G_n(\mathcal{J}\right\})$ in $\left\{[\mathcal{I}^\star] < G_n(\mathcal{J}^\star\right\})$ is an algebra isomorphism $\phi$ of the algebra $\mathcal{J}$ onto $\mathcal{J}^\star$ with  \\
$\phi([\mathcal{I}]) = [\mathcal{I}^\star]]$ $\clubsuit$\\
\end{defn}
\textbf{Reminder.}\\
(1) Let $M$ be a finite dimensional simply connected locally flat manifold. The moduli space of Cartan-Lie groups of $M$ is denoted by
$$[\mathcal{CLG}(LF)(M)].$$
(2): The moduli space of the category of effective pairs is denoted by $[\mathcal{EP}(ASSA)]$.
\section{THE FIRST FUNDAMENTAL EQUATION AND SOME OPEN PROBLEMS. SOME NEW THEOREMS EXTH}
\textbf{Reminder}. For convenience we go to reset problems already raised. In the category of finite dimensional differentiable manifolds we have raised some fundamental open problems.\\
The problem $EX(\mathcal{S})$ whose framework is the global analysis on manifolds.\\
The problem $EXF(\mathcal{S})$ whose framework is the differential topology.\\
The problem $DL(\mathcal{S})$ whose framework has a combinatorial nature.\\
The problem of complexity of statistical models for measurable sets whose framework is the Information Geometry.\\
We have listed some geometric structures $\mathcal{S}$ whose existence
is an open problem up to nowadays.\\
The existence of locally flat structures in finite dimensional manifolds.\\
The existence of left invariant locally flat structures in Lie groups.\\
The existence of Hessian structures in Riemannian manifolds.\\ 
The existence of Hessian structures in locally flat manifolds.\\
The existence of symplectic structures in finite dimensional manifolds.\\
The existence of left invariant symplectic structures in finite dimensional Lie groups.\\
The existence of bi-invariant Riemannian metrics in finite dimensional Lie groups.\\
\textbf{An obervation. Actually there is not hint of thinking that the problems just listed might be controled by one same invariant}.\\
Up to nowadays there does not exist any characteristic obstruction to those problems. A characteristic obstructions are invariants providing necessary and sufficient conditions.\\
We go to link the function $r^b$ with those problems. Those links highlight the novative nature of the fundamental differential equation $FE^*(\nabla)$. The theorems which are stated in the next subsections will be demonstrated in a later section.
\subsection{The first fundamental equation and the locally flat geometry}
For convenience we recall that there does not exist any characteristic obstruction to $EX(\mathcal{S}: LF(M))$. Recently the search of such an invariant seemed hopeless \cite{Medina-Saldarriaga-Giraldo}.
\begin{thm} EXTH1: In a finite dimensional differentiable manifold $M$ the following assertions are equivalent,\\
\label{flatgeometry}
(1): $r^b(M) = 0,$\\
(2): the manifold $M$ admits locally flat structures $\clubsuit$
\end{thm}
\subsection{The first fundamental equation and the Hessian geoemtry in Riemanian manifolds: answer an old question of S-I Amari and A.K. Guts}
We have recalled an old question of Alexander K. Guts and S-I. Amari. We have observed that the same request is implicit in the works of W. Kaup, J-L Koszul and of J. Vey. (The notion of hyperbolicity in the locally flat geometry).
\begin{thm}EXTH2: In a finite dimensional Riemannian manifold $(M,g)$ the following assertions are equivalent\\
\label{hessiangeometry1}
(1): $r^b(M,g) = 0,$\\
(2): the Riemannian manifold $(M,g)$ admits Hessian structures $\clubsuit$
\end{thm}
\textbf{A comment}\\
We recall that the assertion (2) has a homological versus: $M$ admits a locally flat structure $(M,D)$ having the metric tensor $g$ as a 2-cocycle of its KV complex.
\subsection{The first fundamental equation and the Hesian geometry in locally flat manifolds}
\begin{thm} EXTH3: In a finite dimensional locally flat manifold $(M,D)$ the following statements are equivalent,\\
\label{hessiangeometry2}
(1): $r^b(M,D) = 0,$\\
(2): the locally flat manifold $(M,D)$ admits Hessian structures $\clubsuit$
\end{thm}
\textbf{A comment}\\
Actually assertion (2) means that the manifold $M$ admits a Riemannian structure $(M,g)$ whose metric tensor $g$ is a 2-cocycle of the KV complex of $(M,D)$. This theorem is a first step toward the geometry of Koszul in $(M,D)$. A combination of EXTH2 and EXTH3 last yields the following corollary.
\begin{thm} EXTH4: In a finite dimensional differentiable manifold $M$ the following assertions are equivalent,
\label{hessiangeometry3}
(1): $r^B(M) = 0,$\\
(2): $r^{B*} = 0,$\\
(3): the manifold $M$ admits a Hessian structures $\clubsuit$
\end{thm}
\subsection{The first fundamenatal equation and the information geometry}
Let $(\Xi,\Omega)$ be a measurable set. We consider a statistical model for $(\Xi,\Omega)$, namely 
$$\mathbb{M} = [\mathcal{E},\pi,M,D,p].$$
Its Fisher information is denoted by $g$.\\
(i): When the Fisher information $g$ is definite we get the Riemannian Hessian defect  $r^b(M,g)$ and the exponential defect $r^b(\mathbb{M})$.\\
(ii): When the model $\mathbb{M}$ is singular we use the family of its $\alpha$-connection for introducing its exponential defect $r^b(\mathbb{M})$.\\
Therefore we state
\begin{thm}EXTH5: In a regular statistical model $\mathbb{M}$ whose Fisher information is denoted by $g$ the following assertions are equivalent,\\
\label{informationgeometry1}
(1): $r^b(M,g) = 0,$\\
(2): $r^b(\mathbb{M}) = 0$,\\
(3: $\mathbb{M}$ is an exponential family $\clubsuit$
\end{thm}
\begin{thm} EXTH6: In a singular statistical model $\mathbb{M}$ the following assertions are equivalent,\\
\label{informationgeometry2}
(1): $r^b(\mathbb{M}) = 0,$\\
(2): $\mathbb{M}$ is an exponential family $\clubsuit$
\end{thm}
The six Theorems just stated and some avatars of them will be demonstrated in a next section.\\
We recall that our purposes includ the moduli space of bi-invariant affine Cartan-Lie groups.  
\begin{thm} MODS1: The moduli space of Cartan-Lie groups of finite dimensional locally flat manifolds, namely $\mathbb{CLG}(LF)$ is isomoprhic to the moduli space of finite dimensional effective pairs of associative algebras, namely $\mathbb{EP}(ASS)\clubsuit$
\end{thm}
\textbf{Miscellaneous}. In \cite{Medina-Saldarriaga-Giraldo} the authors study the projectively flat structures in finite dimensional Lie groups. This purpose is not discussed in this paper. However the function $r^b$ might be linked with the projectively flat geometry in finite dimensional Lie groups. The concern would be the links of the Associative Algebras Sheaf $(\mathcal{J}_\nabla, \nabla)$ with the Associative Algebras Sheaf $(\mathcal{J}_{\nabla^*},\nabla^*)$ when $\nabla$ is projectively equivalent to $\nabla^*$.\\
\subsection{The hyperbolicity after Kaup, Koszul and Vey}
Before continuing we summarize a construction due to J-L Koszul \cite{Koszul(1)}. It will be involved in studying a variant of the conjecture of Markus.\\
We consider an $m$-dimensional locally flat manifold $(M,\nabla)$. Its universal covering is denoted by
$$(\tilde{M},\tilde{\nabla})\rightarrow (M,\nabla).$$
Up to homeomorphism the topological space $\tilde{M}$ is the space of homotopy class of origin fixed paths
$$\left\{(0,[0,1)\rightarrow (x^*,M)\right\}$$
The topology of $\tilde{M}$ is the compact-open topology.\\
For $s \in [0,1]$ the parallel transport along $c$ of $T_{x^*}M$ in $T_{c(s)}M$ is denoted by
$\tau(s)$. Koszul has introduced the continuous map
$$ \tilde{M} \ni [c]\rightarrow Q_\nabla([c])\in T_{x^*}M]. $$
which is defined by
$$Q_\nabla([c]) = \int^1_0 \tau^{-1}(s)[\frac{dc(s)}{ds}] ds .$$
\begin{defn} We keep the notation just posed.
$(1)$: A locally flat manifold $(M,\nabla)$ is called geometrically complete if the map $Q_\nabla$ is a diffeomorphism of $\tilde{M}$ onto $T_{x*}M$.\\
$(2):$ A locally flat manifold $(M,\nabla)$ is called hyperbolic if the map $Q_\nabla$ is a diffeomorphism of $\tilde{M}$ onto a convex domain not containing any straight line $\clubsuit$
\end{defn}
\textbf{Reminder.}\\
For more details about the notion of hyperbolicity the readers are referred to \cite{Katsumi},  \cite{Kaup}, \cite{Koszul(2)}, \cite{Vey(1)}.\\
For our purpose the KV cohomology (to be introduced) is a unifying framework for investigating the following topics.\\
(i): The geometry of Koszul.\\
(ii: The information geometry.\\
(iii): The theory of Hessian Special Dynamical Systems, namely the data
$$\left\{(\Gamma^0,\tilde{\nabla}^0,g^0),(M,\nabla,g)\right\}\clubsuit$$
Let $(M,\nabla,g)$ be a Hessian manifold. It gives rise to the transitive Associative Algebras Sheaf $(\mathcal{J}_\nabla, \nabla)$. The sections of $\mathcal{J}_\nabla$ are solutions to the first fundamental equation $FE*(\nabla)$. Loosely speaking those sections are zeros of the Hessian differential operator $D^2$. The restriction to $J_\nabla$ of the metric tensor $g$ is denoted by $g^0$. This restriction defines a $C^\infty(M)$-valued cohomology class
$$[g^0] \in H^2_{KV}(J_\nabla, C^\infty(M)).$$
We also recall that the commutator Lie Algebra of  $(\mathcal{J}_\nabla, \nabla^0)$ is the Lie Algebra Sheaf (LAS) of the Cartan-Lie group $\Gamma_0.$
\begin{defn} The triple $(\Gamma_0,\tilde{\nabla}^0,g^0)$ is the Hessian Cartan-Lie group of affine transformation of the Hessian manifold $(M,\nabla,g)$ $\clubsuit$
\end{defn}
The machinery we just discussed will helps to relate the geometric completeness of $(M,\nabla)$ and the geometric completeness of the affine Lie group $(G_0,\nabla^0)$.\\
\section{THE MODULI SPACES OF LOCALLY FLAT STRUCTURES, CANONICAL LINEAR REPRESENTATIONS OF FUNDAMENTAL GROUPS}
\subsection{A long digression}
Since the pioneering works of Elie Cartan the geometry of bounded domains has been a major research subject of many mathematicians among whom are Koszul, Vinberg and their schools. The majors addresses of those researches may be linked with the problem of the moduli space of simply connected locally flat manifolds. The case of hyperbolic locally flat manifolds is linked with the theory of Siegel domains, see {Gindikin-Piatecii-Shapiro-Vinberg} \cite{GPSV}, Vey.
Both J-L Koszul and E.B. Vinberg have implemented abstract algebra machinery to handle the geometry of homogeneous domains \cite{Koszul(4)}, \cite{Vinberg(1)}, \cite{Dorfmeister}, \cite{Dorfmeister-Nakajima}
For instance E.B. Vinberg introduced and used the notion of $T-algebra$ for studying the moduli space of homogeneous convex cones \cite{Vinberg(2)}. J-L Koszul introduced and used the notion of $j-algebra$ for studying the homogeneous Kaehlerian manifolds \cite{GPSV}. J. Dorfmeister and K. Nakajima involved the theory of j-algebra of Koszul in demonstrating the fundamental conjecture of Gindikin-Pyateccii-Shapiro-Vinberg. In this paper we implement the fundamental equation $FE^*(\nabla)$ in demonstrating an alternative to the fundamental conjecture. Our fibration is topologically finer and our demonstation is shorter than the demonstration of Dorfmeister -Nakajima. In the preceding sections we have recalled some relevant categories of geometrical structures and the problems of their moduli spaces. For our next purposes we recall some notions and notation.\\
$1st:$ $\mathcal{M}od(LF)$ is the moduli space of finite dimensional simply connected locally flat manifolds.\\
$2nd:$ $\mathcal{EP}(ASSA)$ is the moduli space of effective pairs of finite dimensional associative algebras.\\
$3rd:$ The category of finite dimensional Hessian manifolds is denoted by $\mathcal{HES}$.\\
Besides the categories just recalled are
$4th:$ The category of Hessian Cartan-Lie groups.\\
$5th:$ the category of semi-projective systems of bi-invariant affine Cartan-Lie groups.\\
$6th:$ The category of semi-inductive systems of bi-invariant affine (abstract) Lie groups $\clubsuit$\\
A challenge is the well understanding of links between those categories and the links between their moduli spaces.\\ The global analysis of the differential equation $FE^*(\nabla)$ is helpful in studying the problems $EX(\mathcal{S})$ some major geometric structures such as the symplecti cstructures. This ends the digression $\clubsuit$\\
\subsection{The first fundamental equation and the canonical linear representations of the fundamental groups}
Our aim is to introduce a familly of linear representations of $\pi_1(M)$. Those linear representations are encoded by the first fundamental equation of gauge structures in $M$. They are used for addressing the moduli space of geometrically complete locally flat structures in finite dimensional manifolds. Every gauge structure $(M,\nabla)$ is involved via its first fundamental equation $FE^*(\nabla)$.\\
Let $(M,\nabla)$ be a gauge structure whose universal covering is denoted by $(\tilde{M},\tilde{\nabla})$. Let $D_\nabla$ and $D_{\tilde{\nabla}}$ be the Hessian differential operators of $(M,\nabla)$ and $(\tilde{M},\tilde{\nabla})$ respectively. The covering projection is denoted by $\pi$. The action of the fundamental group $\pi_1(M)$ in $\tilde{M}$ is $\tilde{D}$ preserving. Therefore the differential of every $\gamma \in \pi_1(M)$ is an isomorphism of the Associative Algebra Sheaf $(\mathcal{J}_{\tilde{D}}, \tilde{D})$. So we obtain a canonical linear action of the group $\pi_1(M)$ in the AAS 
$(\mathcal{J}_{\tilde{D}, \tilde{D}})$, namely 
$$\pi_1(M) \times \mathcal{J}_{\tilde{D}} \ni (\gamma,\xi)\rightarrow \gamma_\star(\xi) \in \mathcal{J}_{\tilde{D}}.$$ 
At one side the equality
$$\pi_\star(\tilde{D}) = D$$
means that given two $\pi_1(M)$-invariant vector fields $X, Y \in \mathcal{X}(\tilde{M})$ one has
$$ \pi_*(\tilde{D}_XY) = D_{\pi_*(X)}\pi_*(Y).$$
Consequently there exists a unique linear of action 
$$ \pi_1(M) \times (\mathcal{J}_D, D) \rightarrow (\mathcal{J}_D, D) $$
such that the differential $\pi_\star$ is an equivariant local isomorphism of Associative Algebras Sheaf. We keep the notation used in the preceding sections. So we consider the associative algebras $(J_D,D)$ and $(J_{\tilde{D}},\tilde{D})$. Therefore we obtain a canonical linear action of $\pi_1(M)$ in the associative algebra $(J_D,D)$, namely
$$\pi_1(M) \times (J_D,D)\rightarrow (J_D,D)$$
Since $(G_D,D^0)$ is simply connected the canonical action above is the differential at the unit element of a unique continuous action in the bi-invariant affine Lie group $(G_D, D^0)$, namely
$$ \pi_1(M)\times (G_D,D^0)\rightarrow (G_D,D^0).$$ 
We emphasize this discussion by the following statement
\begin{prop} A gauge structure $(M,\nabla)$ defines a canonical homomorphism of the fundamental group $\pi_1(M)$ in the transformations group of the simply connected bi-invariant affine Lie group
$(G_\nabla,\nabla) \clubsuit$
\end{prop}
This canonical dynamics will be involved in studying the moduli space of geodesically complete locally flat manifolds. The canonical representation of the fundamental group has a natural extension in the semi-inductive systems and in the semi-projective systems that we have constructed in preceding sections.\\
We go to state a generic theorem
\begin{thm} The fundamental group of a differentible manifold have a canonical representation in the automorphisms group of an object of the following categories,\\
$Cat.1.1$: semi-inductive systems of finite dimensional bi-invariant affine Cartan-Lie groups,\\
$Cat.1.2$: semi-projective systems of finite dimensional bi-invariant Cartan-lie group,\\
$Cat.2.1$: semi-inductive systems of finite dimensional simply connected bi-invariant affine (abstract) Lie groups,\\
$Cat.2.2$: semi-projective systems of finite dimensional simply connected bi-invariant (abstract) affine Lie groups,\\
$Cat.3.1$: inductive systems of finite dimensional bi-invariant affine Lie groupoids,\\
$Cat.3.2$: projective systems of finite dimensional bi-invariant affine Lie groupoids $\clubsuit$
\end{thm}
\textbf{Warning}\\
Currently some litteratures dealing with Lie algebroids aim in extending to Lie algebroids the third theorem of Lie for finite dimensional Lie algebras. They use the term "Lie groupoid". In  this paper the term "Lie groupoid" is used in its historical menaning. Lie groupoids are objects of the global analysis on manifolds, (cf MSC2010). From the duality viewpoint 
$$ Equation\leftrightarrow Solution$$ 
Lie groupoids are (truncated) formal Solutions to Lie Equations, \cite{Singer-Sternberg}, \cite{Kumpera-Spencer} \cite{KU-RO}, \cite{Guillemin-Sternberg}, \cite{Kuranishi}, \cite{Kodaira}, \cite{Spencer}.
\subsection{The moduli space of geometrically complete locally flat manifolds}
Let $(M,\nabla)$ be a geodesically complete locally flat manifold whose universal covering is denoted by $(\tilde{M},\tilde{\nabla})$. We go to recall a technical construction. To every $v \in T_xM$ one assigns the Cauchy problem
$$ \left\{\nabla_{\frac{dc(t)}{dt}}\frac{dc(t)}{dt} = 0, c(0) = x, \frac{dc(t)}{dt}(0) = v.\right\}$$
The unique maximal solution to this Cauchy problem is defined in an open interval
$$ |t| < |t_0| \subset \mathbb{R}.$$  
The maximal solution is denoted by $c_{\left\{x,v\right\}}(t)$. If $ 1 < |t_0|$ then $(x,v)$ belongs to the domain of the exponential map and
$$Exp_\nabla(v) = c_{\left\{x,v\right\}}(1).$$
The group of affine transformations of $(M,\nabla)$ is denoted by $Aff(M,\nabla)$. 
After \cite{Koszul(1)} the domain of $Exp_\nabla$ is a universal covering of $M$. We go to deal with this construction of a universal convering of $M$. Therefore every  $\phi \in Aff(M,\nabla)$ gives rise to an action 
$$\left\{\phi,\phi_\star\right\} \in Aff(\tilde{M},\tilde{\nabla});$$ 
actually this action of $Aff(M,\nabla)$ in $(\tilde{M},\tilde{\nabla})$ commuts with the action of the fundamental group $\pi_1(M)$ in $(\tilde{M},\tilde{\nabla})$. To make simple we say that the $Aff(M,\nabla)$ acts in the triple $\left\{\pi_1(M),\tilde{M},\tilde{\nabla}\right\}$.
In the next we aim at involving two dynamics in addressing the problem of mduli space of locally flat manifolds.\\ 
The first dynamic is the gauge dynamic
$$\left\{\pi_1(M)\times(\tilde{M},\tilde{\nabla})\rightarrow (\tilde{M}, \tilde{\nabla})\right\}.$$ 
The second dynamic is the canonical representation of $\pi_1(M)$ in the bi-invariant affine Lie group $(G_0,\nabla^0)$
$$\left\{\pi_1(M)\times (G_0,\nabla^0)\rightarrow (G_0,\nabla^0)\right\}.$$
\textbf{Reminder}. By localization the second dynamic leads to the canonical representation of $\pi_1(M)$ in the bi-invariant affine Cartan-Lie group $(\Gamma_{\nabla^0},\tilde{\nabla}^0)$
$$\left\{\pi_1(M)\times(\Gamma_{\nabla^0},\tilde{\nabla}^0)\rightarrow (\Gamma_{\nabla^0},\tilde{\nabla}^0)\right\}.$$
A step toward our aim is the following proposition.
\begin{prop} Consider $m$-dimensional geodesically complete locally flat manifolds $(M,\nabla)$ and $(M^*,\nabla^*)$ and their gauge dynamics 
$$\left\{\pi_1(M),(\tilde{M},\tilde{\nabla}^*)\right\},$$ 
$$\left\{\pi_1(M^*),(\tilde{M}^*,\tilde{\nabla}^*)\right\}.$$ 
Then the following assertions are equivalent\\
(1): The locally flat manifold $(M,\nabla)$ is isomorphic to $(M^*,\nabla^*)$,\\
(2): The gauge dynamic 
$$\left\{\pi_1(M)\times(\tilde{M},\tilde{\nabla})\right\}$$ 
is isomorphic to the gauge dynamic 
$$\left\{\pi_1(M^*)\times (\tilde{M}^*,\tilde{\nabla}^*)\right\}\clubsuit$$
\end{prop}
\textbf{Reminder.\\
An isomorphism of 
$$\left\{\pi_1(M),(\tilde{M},\tilde{\nabla})\right\}\rightarrow \left\{\pi_1(M^*),(\tilde{M}^*,\tilde{\nabla}^*)\right\}$$  
is a pair $(\Phi,\psi)$ formed of an affine diffeomorphism
$$\Phi: (\tilde{M},\tilde{\nabla})\rightarrow (\tilde{M}^*,\tilde{\nabla}^*)$$
and a group isomorphism  
$$\psi:\pi_1(M)\rightarrow \pi_1(M^*)$$
Such that
$$\Phi\circ \gamma = \psi(\gamma)\circ \Phi \quad \forall \gamma \in \pi_1(M),$$
$$\Phi_\star(\tilde{\nabla}) = \tilde{\nabla}^*.$$
Lossely speaking $\Phi$ is $\pi_1(M),\pi_1(M^*)$-equivariant}.\\
\textbf{The idea of proof of Proposition}.\\
We go to prove that (2) implies (1).\\
Assume that there exists an affine diffeomorphism
$$\phi:(M,\nabla)\rightarrow (M^*,\nabla^*)$$
such that
$$\phi\circ \pi = \pi^*\circ\Phi.$$
Here $\pi$ and $\pi^*$ are (universal) covering maps. Of course one has
$$\phi_\star(\nabla) = \nabla^*.$$
The last property of $\phi$ yields
$$\phi_\star(\tilde{\nabla}) = \tilde{\nabla}^*.$$
In final the differential $\phi$ is an isomorphism of $(\tilde{M},\tilde{\nabla})$ on 
$(\tilde{M}^*,\tilde{\nabla}^*)$ which AGREES the actions of the fundamental groups 
$$\left\{\pi_1(M)\times (\tilde{M},\tilde{\nabla})\rightarrow ( \tilde{M},\tilde{\nabla})\right\},$$
$$\left\{\pi_1(M^*)\times (\tilde{M}^*,\tilde{\nabla}^*)\rightarrow (\tilde{M}^*,\tilde{\nabla}^*  )\right\}.$$
We take into account the equality
$$ \phi(\pi_1(M)) = \pi_1(M^*),$$
therefore this AGREEMENT is but the assertion (1).\\
It is a folklore that (1) implies (2)$\clubsuit$\\
\textbf{Reminder.}\\
The pair $\left\{\pi_1(M),(J_{\tilde{\nabla}},\tilde{\nabla})\right\}$ stands for the canonical representation of $\pi_1(M)$
$$\pi_1(M)\times J_{\tilde{\nabla}} \ni (\gamma,\xi)\rightarrow \gamma_{\star}(\xi)\in J_{\tilde{\nabla}} \clubsuit $$
Every $\gamma_\star$ is an automorphism of the associative algebra 
$(J_{\tilde{\nabla}},\tilde{\nabla})$
Now let $(G_0,\nabla^0)$ be the simply connected bi-invariant affine Lie group whose Lie algebra is 
$$(J_{\tilde{\nabla}}, [-,-]_{\tilde{\nabla}}).$$ 
The canonical representation 
$$\left\{\pi_1(M),(J_{\tilde{\nabla}},\tilde{\nabla})\right\}$$ 
is the infinitesimal counterpart of a well defined dynamic 
$$\left\{\pi_1(M)\times(G_0,\nabla^0)\rightarrow (G_0,\nabla^0)\right\}$$.\\
We go to deal with the problem of moduli space.\\
The following theorem is hepful 
\begin{thm} [\cite{Kobayaschi}, Theorem 6.16]. Assume that a gauge structure $(M,\nabla)$ is geodesically complete. Then every infinitesimal transformation of $(M,\nabla)$ $X$ is complete, viz $X$ generates a one parameter subgroup of affine transformations 
$$\left\{\phi_X(t)\right\}\subset Aff(M,\nabla) \clubsuit$$
\end{thm}
 We go to use the following corrolary of the theorem of Kobayashi, see also \cite{Hano-Morimoto}
\begin{thm} Let $(M,\nabla)$ be a finite dimensional geodesically complete locally flat manifold. The Lie group of affine transformations of $(M,\nabla)$ admits a structure of bi-invariant affine Lie goup whose universal covering is the simply connected bi-invariant affine Lie group $(G_\nabla, \nabla^0)$. Further the Koszul connectiong $\nabla$ is induced by the locally effective transitive affine dynammic 
$$(G_\nabla,\nabla^0)\times(M,\nabla)\rightarrow (M,\nabla)\clubsuit$$
\end{thm}
\textbf{Rimder}: Let $H \subset G_\nabla$ be the isotropy subgroup at $x\in \tilde{M}.$ 
Let $\mathcal{H} \subset J_{\tilde{\nabla}}$ be the Lie algebra of $H$. Then $\mathcal{H}$ is a simple right ideal of $J_{\tilde{\nabla}},\tilde{\nabla}$. Up to isomorphism 
$(\tilde{M},\tilde{\nabla})$ is defined by the effective pair
$$\left\{[\mathcal{H}] < J_{\tilde{\nabla}}\right\}$$
This effective pair is preserved by the canonical representation of $\pi_1(M).$
So we have the canonical linear dynamic
$$ \left\{\pi_1(M)\times \left\{[\mathcal{H}] < J_{\tilde{\nabla}}\right\}\right\}$$
Another step toward our aim is the following statement.
\begin{thm} We consider geodesically complete locally flat manifolds $(M,\nabla)$ and $(M^*,\nabla^*)$ and the corresponding canonical representations 
$$\left\{\pi_1(M)\times(J_{\tilde{\nabla}},\tilde{\nabla})\right\}$$ 
and 
$$\left\{\pi_1(M^*)\times(J_{\tilde{\nabla}^*},\tilde{\nabla}^*) \right\}.$$ The following assertions are equivalent\\
$$(1):\quad\left\{\pi_1(M)\times (\tilde{M},\tilde{\nabla})\right\}$$ 
is isomorphic to 
$$\left\{\pi_1(M^*)\times (\tilde{M}^*,\tilde{\nabla}^*)\right\},$$
$$(2):\quad \left\{\pi_1(M)\times(J_{\tilde{\nabla}},\tilde{\nabla})\right\}$$ 
is isomorphic to 
$$\left\{\pi_1(M^*)\times(J_{\tilde{\nabla}^*},\tilde{\nabla}^*)\right\}\clubsuit$$
\end{thm}
\textbf{Demonstration}.\\
It is easy to see that (1) implies (2).\\
We go to prove that (2) implies (1).\\
The covering map
$$ \pi: (\tilde{M},\tilde{\nabla})\rightarrow (M,\nabla)$$
yields the local isomorphism of the Associative Algebras Sheaf 
$$ \pi_\star: (\mathcal{J}_{\tilde{\nabla}},\tilde{\nabla})\rightarrow (\mathcal{J}_\nabla, \nabla).$$
By \cite{Kobayaschi} the associative algebra $(J_{\tilde{\nabla}}, \tilde{\nabla})$ is derived from a well defined transitive locally effective action
$$(G_\nabla,\nabla^0)\times (\tilde{M},\tilde{\nabla})\rightarrow (\tilde{M},\tilde{\nabla}).$$
Suppose one has a linear isomorphism
$$ \psi\times\Phi: \left\{\pi_1(M)\times \left\{[\mathcal{H}] < J_{\tilde{\nabla}}\right\}\right\}\rightarrow \left\{\pi_1(M^*)\times \left\{[\mathcal{H}^*] < J_{\tilde{\nabla}^*}\right\}\right\}$$
Since the Lie groups $G_\nabla$ and $G_{\nabla*}$ are simply connected we have the following consequences:\\
(1*): there exists a unique isomorphism of bi-invariant affine Lie groups whose differential at the unit element is $\Phi$, namely
$$ \phi: (G_\nabla,\nabla^0)\rightarrow (G_{\nabla*},\nabla^{*0}),$$
(2*): the action of $\pi_1(M)$ in $J_{\tilde{\nabla}},\tilde{\nabla})$ is the differential of a well defined action of $\pi_1(M)$ in $(G_\nabla,\nabla^0)$.
(3*): the same argument yields a well defined action of $\pi_1(M^*)$ in $(G_{\nabla{*0}},\nabla^{*0})$.\\
(4*): We consider the effective pairs 
$$\left\{[H] < G_\nabla\right\},$$  
$$\left\{[H^*] < G^{*0} < G_{\nabla*}\right\};$$ 
they are defined by 
$$\left\{[\mathcal{H}] < J_{\tilde{\nabla}}\right\},$$  
$$\left\{[\mathcal{H}^*] < J_{\tilde{\nabla}}^*\right\}.$$ 
Then we conclude that $\psi\times\Phi$ is the linear counterpart of a well defined isomorphism 
$$\left\{\pi_1(M)\times [H] < G_\nabla\right\}\rightarrow \left\{\pi_1(M^*)\times [H^*] < G_{\nabla^{*0}}\right\}$$
The theorem is demonstrated $\clubsuit$\\
\subsection{An algebraic model for complete locally flat geometry}
We keep the notation used in the last subsections. The theorems proved in the later subsection allow the following identifications\\
$$ \left\{\pi_1(M)\times \left\{[\mathcal{H}] < J_{\tilde{\nabla}}\right\}\right\}\leftrightarrow \left\{\pi_1(M)\times \left\{[H] < G_\nabla\right\}\right\}, $$
$$ \left\{\pi_1(M)\times \left\{\tilde{M},\tilde{\nabla}\right\}\right\}\leftrightarrow\left\{\pi_1(M)\times\left\{[H] < G_\nabla,\nabla^0\right\}\right\},$$
Actually the simply connected bi-invariant affine Lie group $(G_\nabla,\nabla^0)$ suports the following dynnamics\\
(1): the canonical representation of $\pi_1(M)$
$$ \pi_1(M)\times (G_\nabla,\nabla^0)\rightarrow (G_\nabla,\nabla^0),$$ 
(2): the left translations by the isotropy subsgroup $H$ 
$$H\times(G_\nabla,\nabla^0)\rightarrow (G_\nabla,\nabla^0).$$
Those dynamics agree each withother, further their double orbit space is the complete locally flat manifold $(M,\nabla)$, thereby we get
$$ (M,\nabla) = H\backslash \left\{G_\nabla,\nabla^0\right\}/\pi_1(M).$$
Consequently up to isomorphism the geodesically complete locally flat manifold $(M,\nabla)$ is the same think as the canonical representation $$\left\{\pi_1(M),\left\{[\mathcal{H}] < J_{\tilde{\nabla}}\right\}\right\}$$
In other words the locally effective pair of associative algebra
$$\left\{[\mathcal{H}] < J_{\tilde{\nabla}}\right\}$$
is the same thing as the canonical representaion of the fundamental group
$$\pi_1(M)\times \left\{[\mathcal{H}] < J_{\tilde{\nabla}}\right\}\rightarrow\left\{[\mathcal{H}] < J_{\tilde{\nabla}}\right\}.$$
Actually the canonical linear dynamic is effective. In final we pose
\begin{defn}
The effective linear dynamic $\left\{\pi_1(M)\times\left\{[\mathcal{H}] < J_{\tilde{\nabla}}\right\}\right\}$ is called the algebraic model of $(M,\nabla)$
\end{defn}
The identifications we have made lead to the following statement
\begin{thm}
The moduli space of finite dimensional geodesically complete locally manfolds is equivalent to the moduli space of effective pair of finite dimensional associative algebras $\clubsuit$
\end{thm}
\section{A SKETCH OF THE KV ALGEBRAIC TOPOLOGY OF LOCALLY FLAT MANIFOLDS}
To handle some intrinsic meaning of both the fundamental equation $FE^*(\nabla)$ and the function $r^b$ we shall be using more or less explicitly the KV homological. This short section is devoted to briefly present the theory of KV homology in the locally flat geometry. Let $(M,D)$ be a locally flat manifold. The vector space of differentiable vector felds $\mathcal{X}(M)$ is endowed with the abstract algebra structure whose multiplication is defined by
$$X.Y = D_XY.$$
The pair
 $$\mathcal{A} = (\mathcal{X}(M),D)$$ 
is called the Koszul-Vinberg algebra of $(M,D)$, KV algebra in short. The vector space $\mathcal{X}(M)$ is canonically a two-sided module of the algebra $\mathcal{A}$. The vector space $C^\infty(M)$ is a left module of $\mathcal{A}$. The left action is defined as
$$X.f = df(X).$$
We mentione that the source of our interest in $FE^*(D)$ was a conjecture of Muray Gerstenhaber. This conjecture claims that \textbf{every restricted theory of deformation generates its proper theory of cohomology}. Regarding this conjecture in the category of KV algebra the pioneering attempt was made by Albert Nijenhuis \cite{Nijenhuis}. Our concern was to demonstrate the conjecture of Gerstenhaber in the locally flat geometry \cite{Koszul(1)}. Details and applications in the information geometry can be found in \cite{Gerstenhaber}, \cite{Nguiffo Boyom(2)}, \cite{Nguiffo Boyom(6)}, \cite{Nijenhuis}, \cite{Koszul(1)}, \cite{Barbaresco}.
 \subsection{The vector valued cohomology}
We keep the notion used in the preceding sections. Our concern is the $\mathcal{X}(M)$-valued KV complex of $\mathcal{A}$. That complex is the $\mathbb{Z}$-graded vector space
$$C(\mathcal{A}) = \oplus_q C^q(\mathcal{A}).$$
Here the homogeneous subspaces $C^q(\mathcal{A})$ are defined by 
$$C^q(\mathcal{A}) = 0\quad \forall q < 0,$$
$$ C^0(\mathcal{A}) = J_D,$$
$$ C^q(\mathcal{A}) = Hom_\mathbb{R}(\mathcal{A}^{\otimes q},\mathcal{A}).$$
We go to use the brute coboundary operator $\delta$ as in \cite{Nguiffo Boyom(3)}, see also \cite{Nguiffo Boyom(6)} for alternative definition of the coboundary operator. This linear brute coboundary operator of degree $+1$ is defined as it follows
$$(1):\quad \delta\xi(X) = -X.\xi + \xi.X \quad \forall \xi \in C^0(\mathcal{A}),$$
$$(2):\quad \delta f(\xi) = \sum^q_1(-1)^i[X_i.f(\partial_i\xi) +
(f(\partial^2_{i,q+1}\xi\otimes X_i)).X_{q+1} - f(X_i.\partial_i\xi)]. $$
Here we have 
$$\xi = X_1\otimes..\otimes X_{q+1}\in \mathcal{A}^{\otimes q+1},$$
$$\partial_i\xi = X_1\otimes..\hat{X}_i\otimes..\otimes X_{q+1},$$
$$\partial^2_{i,q+1}\xi = X_1\otimes..\hat{X}_i..\otimes X_q \otimes\hat{X}_{q+1},$$
$$ Y.\xi = \sum^{q+1}_1 X_1\otimes..\otimes Y.X_j\otimes..\otimes X_{q+1}.$$
The symbol ${\hat{X}_i}$ indicates that $X_i$ is missing. The pair $(C(\mathcal{A}),\delta)$ is a cochain complex. Its $q^{th}$ cohomology space is denoted by $ H^q_{KV}(\mathcal{A})$. There are two other definitions of the coboundary operator \cite{Nguiffo Boyom(6)}. Those alternative definitions yield two cochain complexes which are quasi isomorphuic to that the complex we just defined.
\subsection{The scalar valued cohomology}
The scalar KV complex is the $\mathbb{Z}$-graded vector space
$$C(\mathcal{A},\mathbb{R}) = \oplus_q C^q(\mathcal{A},C^\infty(M)).$$
The homogeneous vector subspaces are defined by \\
$$C^q(\mathcal{A},\mathbb{R}) = 0\quad \forall q  < 0,$$
$$C^0(\mathcal{A},\mathbb{R}) = \left\{f \in C^\infty(M)|\quad D^2(f) = O,\right\}$$
$$C^q(\mathcal{A},\mathbb{R}) = Hom_\mathbb{R}(\mathbb{A}^{\otimes q}, C^\infty(M)).$$\\
The coboundary operator is defined as it follows\\
$$\delta f = - f \quad\forall f\in C^0(\mathcal{A},\mathbb{R}),$$
$$\delta f(\xi) = \sum^q_1(-1)^i[X_i.f(\partial_i\xi) - f(X_i.\partial_i\xi)].$$
The $q^{th}$ cohomology space of $C(\mathcal{A},\mathbb{R})$ is denoted by
$$H^q_{KV}(\mathcal{A},\mathbb{R}).$$
Every Riemannian metric tensor in $M$ $g$ is a 2-cochain of $C(\mathcal{A},\mathbb{R})$\\
\begin{thm} Let $g$ be a metric tensor in a locally flat manifold $(M,D)$. The following statements are equivalent.\\
(1): $(M,g,D)$ is a Hessian manifold.\\
(2): $\delta g = 0.$\\
(3): Every point has a neighborhood $U$ supporting a local smooth function $h$ such
$$ g = D^2 h \clubsuit $$
\end{thm}
This theorem is a partial rephrasing of Proposition 2.1 as in  \cite{Shima(2)}.\\
We rephrase Theorem 3 as in \cite{Koszul(1)}\\
\begin{thm} In a compact positive Hessian manifold $(M,g,D)$ the following assetions are equivalent\\
(1) $[g] = O \in H_{KV}(A,\mathbb{R})$,\\
(2) is hypebolic $\clubsuit$
\end{thm}
\section{NAIVE HOMOLOGICAL ALGEBRA, continued}
In this short section we go to emphasize perspectives for new developments. We fix a finite dimensional manifold $M$. The SET of Koszul connexions in $M$ is denoted by $\mathcal{LC}(M)$. As seen in preceding sections we have emphasized some important mathematical structures in the set 
$\mathcal{LC}(M)$.\\
(a) The set $\mathcal{LC}(M)$ is an infinite dimensional affine space. This structure is linked with the theory of analysis-topology deformation of linear connections, \cite{Koszul(1)}, \cite{Raghunathan}.\\
(b) $\mathcal{LC}(M)$ may be regarded as a convex subset of the vector space of Differential operators of $TM\otimes TM$ in $TM$.\\
(c) $\mathcal{LC}(M)$ may be regarded as the set of vertices of the connected graph $\mathcal{GR}(LC)$. By the functor of Amari every Riemannian a tensor metric $g$ is a symmetry of $\mathcal{GR}(LC)$. Thereby $\mathcal{R}ie(M)$ generates a group $\mathit{COX}(M)$ which acts effectively in the graph $\mathcal{GR}(M)$. So the metric tensors have a dynamical nature. We have decided not to address this perspectve in the present paper.
\subsection{Some remarkable homological functors in gauge structures}
We go to emphasize the following categories\\ 
C.1 : the category of gauge structures
$$\mathcal{MG}(M) = \left\{(M,\nabla) | \nabla \in \mathcal{LC}(M)\right\}$$
C.2 : the category Associative Algebras Sheaves on $M$ which is denoted by
$$\mathcal{J}(M)$$
C.3 : the category of finite dimensional simply connected bi-invariant affine Lie groups denoted by
$$\mathcal{BALG}$$
C.4 : the category of inductive systems of finite dimensional bi-invariant affine Lie groups.\\
C.5 : the category of projective systems of finite dimensional bi-invariant affine Lie groups.\\
We have introduced the Hessian functor
$$\nabla \rightarrow D_\nabla = \nabla^2.$$  
We have used the Hessian operator for defining the functors\\
$$\mathcal{MG}(M)\ni (M,\nabla)\rightarrow \mathcal{J}_\nabla \in \mathcal{J}(M),$$
$$\mathcal{MG}(M)\ni (M,\nabla)\rightarrow (J_\nabla,\nabla),$$
$$\mathcal{MG}(M)\ni (M,\nabla)\rightarrow (G_\nabla,\nabla^0),$$
$$\mathcal{MG}(M)\ni (M,\nabla)\rightarrow \left\{(G_q,\nabla^q):| \Pi^q_p; I^p_q | p\leq q\right\}$$
Those functors might impact the differential topology of the underlying manifold $M$. They also emphasize the richness of the category $\mathcal{MG}(M)$.\\
\textbf{Remind}\\
We have defined the following functor
$$\frac{\mathcal{MG}(M)}{\mathit{G}(M)} \ni [\nabla]\rightarrow r^b([\nabla])\in\mathcal{Z}.$$
This last functor $r^b$ is widely used in this paper. 
\section{NEW (co)HOMOLOGICAL FUNCTORS, continued}
The purpose of this section is to introduce new homological functor defined in the category of gauge structures. Those new functors are global geometric invariants of finite dimensional manifolds. Those functors are not studied in this paper. However they deserve the attention.\\
To a gauge structure $(M,\nabla)$ we have assigned the differential operators
$$D^\nabla(X) = L_X\nabla - \iota_XR^\nabla,$$
$$D_\nabla(X) = \nabla^2X.$$
They are involutive 2nd order differential operators of type $2$. The sheaves of germs of their zeros are denoted by $\mathcal{J}(\nabla)$ and $\mathcal{J}_\nabla$. \\
We go to focus on the Associative Algebras Sheaf
$$\left\{\mathcal{J}_\nabla,\nabla\right\}.$$
Behind the global analysis of $FE^*(\nabla)$ are many classical homological complexes.
\subsection{A sketch of the algebraic topology}
A symmetric gauge structure $(M,\nabla)$ defines the AAS $\left\{\mathcal{J}_\nabla,\nabla\right\}$ and the LAS $\left\{\mathcal{J}_\nabla, [-,-]_\nabla\right\}$. 
The sheaf of germs of vector fields in $M$ is a sheaf of is RIGHT modules of the AAS 
($\mathcal{J}_\nabla,\nabla)$. Thereby we get the following cohomology Sheaves. 
\subsubsection{The Hochschild cohomology sheaf of the AAS $(\mathcal{J}_\nabla,\nabla)$} 
$$HH^\star([\mathcal{J}_\nabla,\nabla],\mathcal{X}(M)),$$
\subsubsection{The KV cohomology sheaf of the KVAS $(\mathcal{J}_\nabla,\nabla)$}
$$H^\star_{KV}([\mathcal{J}_\nabla,\nabla],\mathcal{X}(M))$$
\subsubsection{The Chevalley-Eilenberg cohomology of the sheaf of Lie algebras $(\mathcal{J}_\nabla,[-,-]_\nabla)$}
$$H^\star_{CE}([J^\star_\nabla, [-,-]_\nabla],\mathcal{X}(M))$$
The homological complexes just pointed out are sources new geometric invariants of gauge structures such as spectral sequences \cite{McCleary}. Arise the question how those new invariants impact the topology of underlying manifolds

\subsection{More sophisticated}
\textbf{All of those cohomology sheaves are gauge in variants of the gauge structure $(M,\nabla)$}.\\
We have introduced two new affine geometry invariants.\\
(1): The inductive system and the projective system of bi-invariant affine Cartan-Lie groups
$$\left\{ (\Gamma_q,\tilde{\nabla}^q)\pi^q_p; I^p_q| p < q\right\}$$
(2) the inductive system and the projective system of bi-invariant abstract affine Lie groups
 $$\left\{(G_p,\nabla^p), I^p_q |p < q\right\} \clubsuit$$
According to the problem $EXF(\mathcal{S})$ those gauge invariants impact the differential topology of $M$ whenever the gauge structure $(M,\nabla)$ is symmetric. Those prespective are not not studied in this work.
\subsection{The first fundamental equation and deformations of the Poisson bracket of vectors fields}
Henceforth the concerns are transitive non symmetric gauge structures.\\
\textbf{Reminder}. The transitivity of $(M,\nabla)$  means that
$$r^b(\nabla) = dim(M).$$
Therefore $\mathcal{J}_\nabla$ is a transitive Associative Algbera Sheaf.\\
Thus we face a new structure of Lie algebra in the vector space of vector fields $\mathcal{X}(M)$, namely $g_\nabla$. The bracket of $g_\nabla$ defined by 
 $$[X,Y]_\nabla = \nabla_XY - \nabla_YX$$
 $\mathcal{X}(M)$.
This new structure of Lie algebra is a deformation of the Lie algebra $g$ whose bracket is the Poisson bracket
$$[X,Y] = X\circ Y - Y\circ X.$$
Here $X$ and $Y$ are linear endomorphisms of the vector space $C^\infty(M)$ and
$$ (X\circ Y)(f) = X(Y(f)).$$
The theory of deformations of the Lie algebras $g$ is controlled by the Chevalley-Eilenberg cohomology 
$$H_{CE}(g,g) = \frac{ker(d)}{im(d)}.$$
Here $d$ is the coboundary operator of Chevalley-Eilenberg. The skew symmetric tensor
$$ B(X,Y) = X\circ Y - Y\circ X - [X,Y]_\nabla $$
is also called a deformation of $g$. Thereby it a zero of the polynomial of Maurer-Cartan of both brackets $[-,-]$ and $[-,-]_\nabla$. So for all triple of vector fields $(X, Y, Z$ one has
$$dB(X,Y,Z) + J_B(X,Y,Z) = 0.$$
Here $J_B$ is the Jacobi anomaly function of the abstract algebra $\left\{\mathcal{X}(M),B\right\}$, 
\cite{Nguiffo Boyom(6)}. The anomaly function $J_B$ is defined by 
$$J_B(X,Y,Z) = \Sigma_{cyclic} B(X,B(Y,Z)).$$
\subsubsection{Application to finite dimensional Lie group}
The address in the last subsubsection emphasizes the nature of transitive left invariant gauge structures in finite dimensional simply connected Lie groups. Consider a simply connected Lie group $G$ the Lie algebra of which denoted by $\mathcal{G}$. Let $(G,\nabla)$ be a transitive left invariant gauge structure. We recall that $(G,\nabla)$ is called transitive if the associative algebra $(J_\nabla,\nabla)$ is transitive in $G$. The subset transitive left invariant gauge structures is denoted by 
$$tlc(\mathcal{G})\subset lc(\mathcal{G})$$
Actually one has
$$ lf(\mathcal{G}) = tlc(\mathcal{G})\cup slc(\mathcal{G}).$$
We go to interpret the subset
$$    defor(\mathcal{G}) = tlc(\mathcal{G})\setminus lf(\mathcal{G}).$$
A deformation of the Lie algebra $\mathcal{G}$ is a deformation of the simply connected Lie group $G$. According to the last subsubsection we go to state
\begin{prop} Every $(G,\nabla) \in defor(\mathcal{G})$ is a deformation of the simply connected Lie group $G \clubsuit$
\end{prop}
Readers interessed in the theory of deformation of algebraic structures are referred to \cite{Piper}, \cite{Nijenhuis-Richardson}.
\section{INTRINSIC MEANINGS OF THE FIRST FUNDAMENTAL EQUATION}
Really through the space of solutions to $FE^*(\nabla)$ we are intersted in intrinsic meanings of the function $r^b$. In preceding sections we have used the function $r^b$ for introducing characteristic obstructions to $EX(\mathcal{S})$. In this section we go to face a new challenge. We need to well understand the intrinsic nature of $r^b$.\\
Firstly we are interested in the links of $r^b$ with the combinatory topology, this means the links with the problem $DL(\mathcal{S})$ in the graphs we have introduced, viz $\mathcal{GR}(LC)$, $\mathcal{GR}(SLC)$ and $\mathcal{GR}(LF)$.\\
Secondly we are interested in the links of $r^b$ with the differential topology, this means the links with the problem $EXF(\mathcal{S})$ in the category of finite dimensional manifolds $\clubsuit$\\
Before proceding we recall that the problems $DL(\mathcal{S})$ and problem $EXF(\mathcal{S})$ are the search of solutions to the following questions.\\
\textbf{(1*): How far from admitting $\mathcal{S}$ is a given manifold $M$?} \\
\textbf{(2*): How far from admitting a non discrete $\mathcal{S}$-foliation is a given manifold $M$?} \\
Odd dimensional manifolds do not admit symplectic structure. However some odd dimensional manifolds are foliated by symplectic manifolds. An instance is $\mathbb{R}\times\mathbb{C}.$
\subsection{Definition-Notation: the special dynamical systems}
We go to deal with differentiable dynamical systems, viz the differentiable dynamics of finite dimensional abstract Lie groups.\\
\subsubsection{Geodesically complete gauge structures}
\begin{defn} Let $\mathbb{CMG}(M)$ be the category whose objects are geodesically complete gauge structures $(M,\nabla)$ \\
We consider two objects of $\mathcal{CMG}(M)$ $(M,\nabla)$ and $(M^*,\nabla^*)$. A morphism of $(M,\nabla)$ in $(M^*,\nabla^*)$ is an embedding $\Psi$ of $M$ in $M^*$ satisfying the following identity
$$ \Psi_*(\nabla_XY) = \nabla^*{\Psi_*(X)}\Psi_*(Y)\clubsuit$$
\end{defn}
\textbf{Reminder}. For convenience we recall some important notions.\\ 
The category whose objects are finite dimensional simply connected bi-invariant abstract affine Lie groups is denoted by $\mathbb{BLF}$ $\clubsuit$
An object of $\mathbb{BLF}$ is a couple $(G,\nabla)$ formed of a finite dimensional simply connected abstract Lie group $G$ and a two-sided invariant locally flat Koszul connection $\nabla$\\
Up to group isomorphisms a couple $(G,\nabla)$ is the same thing as a finite dimensional associative algebra $(\mathcal{G},\nabla)$ whose commutator Lie algebra $(\mathcal{G},[-,-]_\nabla)$ is the Lie algebra of the Lie group $G$. So an isomorphism of $(G,\nabla)$ in $(G^*,\nabla^*)$ is the same thing as a linear isomorphism of the associative algebra $(\mathcal{G},\nabla)$ in the associative algebra $(\mathcal{G}^*,\nabla^*)$.\\
Now we consider a complete locally flat gauge structure $(M,\nabla)$. Here the two fundamental equations $FE^*(\nabla)$ and $FE^{**}(\nabla)$ have the same sheaf of solutions, viz 
$$ \mathcal{J}(\nabla) = \mathcal{J}_\nabla. $$
 Therefore according to \cite{Kobayaschi} every section of the sheaf $\mathcal{J}_\nabla$ is a complete vector field. According to Palais  \cite{Palais} the sheaf $\mathcal{J}_\nabla$ is the linear counterpart of a locally effective affine action of the finite dimensional simply connected abstract bi-invariant affine Lie group $(G_\nabla,\nabla^0)$ whose associative algebra is the algebra $(J_\nabla,\nabla)$.\\
\textbf{Reminder. For convenience we recall the definition we are dealing with} 
\begin{defn} A left action of a affine Lie group $(G^*,\nabla^*)$ in a gauge structure $(M,\nabla)$ is formed of\\
(1) a locally effective left action
$$ G^*\times M\ni (\gamma,x)\rightarrow \gamma.x \in M,$$
(2) Let $\xi$ be a left invariant vector field in $G^*$  and let $\xi^* $ be the fundamental vector field in $M$ defined by
$$\xi^*(x) = \frac{d(exp(t\xi)(x))}{dt}(0),$$
then given left invariant vector fields $\xi_1, \xi_2$ one has
 $$[\nabla^*_{\xi_1}\xi_2]^{*} = \nabla_{\xi^*_1} \xi^*_2 $$
 (3) An orbit of smallest dimension is called a minimal orbit. That smallest dimension (of orbits) is called the index of the action.
\end{defn}
\begin{defn} A left action \\
$\left\{(G^*,\nabla^*)\times(M,\nabla)\rightarrow (M,\nabla) \right\}$ \\
is called a Special Dynamical System if the action of $G^*$ in $M$ is locally effective, viz the normal subgroup of $G^*$ acting trivially in $M$ is a discrete subgroup of $G^*$.\\
\end{defn}
\subsection{The existence of special dynamical systems}
The concerns of this subsection are extensions of the problems of affine embeddings of real Lie groups as in \cite{Nguiffo Boyom(9)}
The question is whether a gauge manifold $(M,\nabla)$ admits a special dynamic system. The problem is the search of a finite dimensionl affine Lie group $(G^*,\nabla^*)$ with a left dynamic in $(M,\nabla)$ in such a way that $\left\{(G^*,\nabla^*);(M,\nabla)\right\}$ be a Special Dynamical System.\\
\begin{thm} Every compact special symmetric gauge structure $(M,\nabla)$ admits a special dynamical system.
\end{thm}
\textbf{Hint}. Let $(G_\nabla,\nabla^0)$ be the simply connected bi-invariant affine Lie group whose associative algebra is $(J_\nabla,\nabla)$. Since $M$ is compact $J_\nabla$ is the infinitesimal counterpart of a locally effective action of $G_\nabla$ in $M$. It is easy to see that $(G_\nabla,\nabla^0)$ acts locally effectively in $(M,\nabla) \clubsuit$
The optimal nature of the function $r^b$ is highlighted by the following statement.
\begin{thm} Let $(M,\nabla)$ be a compact symmetric gauge structure.\\
 There is no special dynamical system 
$$\left\{(G^*,\nabla^*)\times(M,\nabla)\rightarrow (M,\nabla)\right\}$$
whose index is smaller than $r^b(\nabla)$.
\end{thm}
According to Theorem EXTH1 the numerical invariant $r^b(M)$ is a \textbf{characteristic obstruction} to the existence of locally flat structure in $M$. \\
\textbf{Assume $\nabla$ to be a vertex of the graph $\mathcal{GR}(SLC)$ then the integer $r^b(\nabla)$ measures the \textbf{distance-like} from the vertex $\nabla$ to the sub-graph $\mathcal{GR}(LF)\subset \mathcal{G}(SLC)$}
\textbf{From the Lie Theory point of view the theorem \textbf{emphasizes the dynamical nature} of the invariant $r^b(M)$}.\\
We summarize this discussing as it follows.
\begin{thm} In a compact symmetric gauge structure $(M,\nabla)$\\
$(1): r^b(\nabla)$ measures how far from admitting a transitive Special Dynamical System is $(M,\nabla).$\\
$(2): r^b(M)$ measures how far from admitting any transitive \textbf{gauge special dynamical system} is the manifold $M$ $\clubsuit$
\end{thm}
\textbf{Assertion (2) expresses that $r^b(M)$ is an obstruction to the existence of transitive special gauge dynamical systems, viz transitve actions of finite dimensional bi-invariant affine Lie groups}.\\
The following corollary is an alternative of the existence theorem EXTH1.
\begin{thm} A compact symmetric gauge structure $(M,\nabla)$ admits a transitive special dynamical system if and only if $(M,\nabla)$ is a locally flat manifold $\clubsuit$
\end{thm}
\textbf{A comment}\\
The two theorems just stated walk in the category of non compact gauge structures. The readers must only replace affine Lie group by affine Cartan-Lie Groups or by affine Lie Groupoid.
\subsection{The differential topological nature of the first fundamental equation, continued}
In this subsection we pursue facing the problem $EXF(\mathcal{S}: FL(M))$. In the analytic category the scope of $EXF(\mathcal{S})$ includes singular foliations, \cite{Zeghib}.\\In the non analytic category the theorem of Frobenius cannot be performed, however there exist completely integrable singular differential systems. Notable instances are submanifolds, foliations of Poisson manifolds by sympectic orbits of Hamiltonian vector fields, foliations of Jacobi manifolds by symplectic orbits or contact orbits of Hamiltonian vector fields. Sadely those instances are rare. In \cite{Nguiffo Boyom(6)} we have used methods of the information geometry for studying Riemannian foliations and symplectic foliations in differentiable manifolds. Some reference for positive Riemannian foliations are \cite{Molino}, \cite{Reinhardt}, \cite{Moerdijk-Mrcun}; there the theory of Riemannian foliation deals with the restricted transverse geometry and the restricted transverse topology of foliations. That theory says nothing when signatures of metric tensors are positive. Therefore the methods of the information geometry is more general. In this section the concerns are intrinsic geometry and intrinsic topology of "geodesic" foliations. That is the meaning of the challenge $EXF(\mathcal{S})$.\\
Before pursuing we go to settle some subtle links between the differential topology and the gauge geometry.
\begin{defn} Consider a finite dimensional manifold $M$ endowed with a differentiable dynamical system 
$$G\times M\rightarrow M,$$
(1) a gauge foliation in $M$ is a pair $\left\{\mathcal{F},\nabla_{\mathcal{F}}\right\}$ formed of a foliation $\mathcal{F}$ and a differentiable family of partial Koszul connections 
$\nabla(\mathcal{F})$ which are defined along the leaves of $\mathcal{F}$\\
(2) a foliation $\mathcal{F}$ in a gauge structure $(M,\nabla)$ is called auto-parallel if its leaves are $\nabla$-auto-parallel $\clubsuit$
\end{defn}
\textbf{Coments. An auto-parallel foliation in a gauge structure is a gauge foliation but the converse is wrong; a geodesic foliation is an auto-parallel foliation but the converse is wrong. Every foliation in a positive Riemannian manifold is canonically a gauge foliation}\\
\textbf{Warning} In the definition above foliations may be singular foliations.
\begin{lem} If $(G,\nabla)$ is a finite dimensional bi-invariant affine Lie group then every locally effective differentiable dynamic of $G$ defines a locally flat gauge foliation.
\end{lem}
\textbf{Hint}. Let $\mathcal{G}$ be the Lie algebra of $G$. A locally effective dynamic
$$G\times M\rightarrow M$$
yields a Lie algebra isomorphism of $\mathcal{G}$ in the Lie algebra of vector fields $\mathcal{M}$
$$ \mathcal{G}\ni a \rightarrow a^\star \in \mathcal{X}(M)$$
 Therefore along the orbits of $G$ the partial Koszul connection $\nabla^*$ is defined by
 $$\nabla^*_{a^*}b^* = (\nabla_ab)^*.$$
 The pair $\left\{Orbit(G),\nabla^*\right\}$ is a locally flat gauge foliation in $M$ $\clubsuit$
For our forthcoming purposes we emphasize the following proposition. 
\begin{lem} Let $M$ be a compact manifold satisfying the conditions 
$$ 0 < r^b(M) < dim(M).$$
Then $M$ has the following properties\\
$M$ supports a locally flat gauge foliation $\left\{\mathcal{F},\nabla(F)\right\}$ whose leaves are orbits of the dynamic of a bi-invariant affine Lie group $(G^*,\nabla^*)$. Further the codimension of minimal orbits is $r^b(M)$ $\clubsuit$
\end{lem}
\textbf{Hint}.\\
In the graph $\mathcal{GR}(SLC)$ there exists a vertex $\nabla$ such that
$$r^b(\nabla) = r^b(M).$$
Therefore the gauge structure $(M,\nabla)$ is special. Let $(G^*,\nabla^*)$ be the simply connected bi-invariant affine Lie group the associative algebra of which is $J_(\nabla,\nabla)$.\\
Since $M$ is compact the Lie Algebras Sheaf (LAS) $(\mathcal{J}_\nabla, \nabla)$ is the Lie algebra of a Cartan Lie group Lie $\Gamma$ which is but the localization of the locally effective action of $G^*$ in $M$.
The orbits of $G^*$ define a locally flat gauge foliation $\mathcal{F}$ whose partial Koszul connection is defined by the bi-invariant affine structure $(G^*,\nabla^*)$. A relevant corollary of the lemma we just proved is 
\begin{thm}
(1) In every compact manifold with $0 < r^b(M) < dim(M)$ the problem $EXF(\mathcal{S}: FL(M))$ admits solutions.\\
(2) Suppose a differentiale manifold $M$ to be homogeneous under a proper action of a locally compact topological group $G$. Then the problem $EXF(\mathcal{S}:LF)$ admits regular soltions in $M \clubsuit$
\end{thm}
\textbf{Hint}\\
The claim (1) is a straightforward corollary of the lemma above.\\
Regarding the claim (2) one could involve the slice theorem of Koszul. The simpler is that the proper action of $G$ yields the existence of $G$-invariant positive Riemannian structure $(M,g)$. Therefore let $\nabla$ be the Levi-Civita connection of $(M,g)$. The AAS $(\mathcal{J}_\nabla,\nabla)$ is $G$-homogeneous, thereby the orbits of $\mathcal{J}_\nabla$ form a regular foliation. Those orbits are orbits of the effective infinitesimal action of the simply connected bi-invariant affine Lie group the associative algebra of which is $(J_\nabla,\nabla)$. The claim (2) is proved $\clubsuit$
\subsection{The fundamental gauge equation and regular geodesic foliations} 
We use the notation already posed in preceding section. We choose a triple 
 $$(g,\nabla,\phi) \in \mathcal{R}ie(M)\times\mathcal{SLC}\times J_{\nabla\nabla^g}.$$
Then we consider the unique pair 
$$(\Phi,\Phi^*)\in J_{\nabla\nabla^g}\times J_{\nabla\nabla^g}$$
which is defined by
$$g(\Phi(X),Y) = \frac{1}{2}[g(\phi(X),Y) + g(X,\phi(Y))],$$
$$g(\Phi^*(X),Y) = \frac{1}{2}[g(\Phi^*(X),Y) - g(X,\phi(Y))].$$
\begin{defn} An solution $\phi\in J_{\nabla\nabla^g}$ is called simple if the foliations $Ker(\Phi^*)$ is simple, viz the quotient topology space $\frac{M}{Ker(\Phi^*)}$ is a differentiable manifold.
\end{defn}
(1) The Riemannian foliation $Ker(\Phi)$ is $\nabla$-totally geodesic.\\
(2) The symplectic foliation $(Ker(\Phi^*)$ is a $\nabla$-totally geodesic.\\
Thus in a Riemannian manifold $(M,g)$ every solution to $FE(\nabla\nabla^g)$, namely $\phi J_{\nabla\nabla^g}$ may be identified with a pair of $\nabla$-totaly geodesic foliations. This is another impact of the differential operator $D^{\nabla\nabla^g}$ on the differential topology of $M$. Those quantitative invariants suggest to investigate the impacts of the arrow\\
$\left\{Riemannian\quad Geometry\times Gauge\quad Geommetry\rightarrow Differential\quad Topology\right\}.$\\
In particular a tensor metric metric in $M$ is a functor of the category $\mathcal{MG}(M)$ in the category of pairs of totally geodesic foliations in objects of $\mathcal{MG}(M)$.  
If $(M,g,g^*)$ is a statistical manifold then for every $\phi \in J_{\nabla\nabla^*}$,\\
(3) the pair $\left\{Ker(\Phi),Ker(\Phi^*\right\}$ is a pair of $\nabla$-totally geodesic foliations,\\
(4) the pair $\left\{im(\Phi),im(\Phi^*)\right\}$ is a pair of $\nabla^*$-totally geodesic foliations.\\
(5) We already pointed out the $g$-othogonal decomposition
$$TM = Ker(\Phi) \oplus im(\Phi),$$
$$TM = Ker(\Phi^*) \oplus im(\Phi^*).$$
\begin{thm}\cite{Nguiffo Boyom(6)}
If $\phi \in J_{\nabla\nabla^*}$ is simple then the quotient map
$$\mathcal{K}er(\Phi^*)\rightarrow M \rightarrow \frac{M}{\mathcal{K}er(\Phi^*)}$$
is a Riemannian submersion on a symplectic manifold $\clubsuit$
\end{thm}
\section{THE FIRST FUNDAMENTAL EQUATION AND OPEN PROBLEMS : DEMONSTRATIONS OF THEOREMS EXTH}
In a preceding section we have stated existence theorems (EXTH). This section is partially devoted to demonstrate those theorems.
\subsection{DEMONSTRATION OF THEOREM \ref{flatgeometry}.}
\textbf{Reminder EXTH1: The aim is to prove that (1) $r^b(M) = 0$ if and only if (2) $M$ admits  locally flat structures}\\
\textbf{(1) implies (2)}.\\
According to the assertion (1) one has
$$r^b(M) = 0.$$
There exists a torsion free Koszul connection $\nabla$ such that
$$r^b(\nabla) = dim(M)$$
Therefore the Associative Algebras Sheaf $(\mathcal{J}_\nabla, \nabla)$ is transitive. Thus let $X,Y, Z$ be sections of $\mathcal{J}_\nabla$, the KV anomaly of $\mathcal{J}_\nabla,\nabla)$ $KV(X,Y,Z)$ vanishes identically.\\
At another side one has the following identity
$$KV(X,Y,Z) = R^\nabla(X,Y).Z \quad \forall X, Y, Z \in J_\nabla$$
Since the application
$$X\otimes Y\otimes Z\rightarrow R^\nabla(X,Y).Z $$
is $C^\infty(M)$-multi-linear the transitivity of $\mathcal{J}_\nabla $ implies the flatness of the torsion free Koszul connection $\nabla $. Thus $(M,\nabla)$ is a locally flat manifold.\\
\textbf{(2) implies (1).}\\
Let $(M,\nabla)$ be a locally flat structure in a manifold $M$.\\
The useful tool one needs is the relationship between the differential equations $FE^*(\nabla)$ and $FE^{**}(\nabla)$. Let $X$ be a vector in $M$. Let $L_X$ be the Lie derivative in the direction $X$ and $\iota_X$ is the inner product by $X$. \\
Every torsion free connection $\nabla^\star$ satisfies the identity
$$L_X\nabla^\star = \iota_XR^{\nabla^\star} - D_{\nabla^\star}(X)$$
We apply this identity to a locally flat manifold $(M,\nabla)$ to get the following identity
$$L_X\nabla = - D_\nabla(X)$$
Thereby $\mathcal{J}_\nabla$ is the sheaf of infinitesimal transformations of $(M,\nabla)$
Thus if $\nabla$ is a locally flat, viz 
$$T^\nabla = 0$$ and
$$R^\nabla = 0 $$
then the Lie equation $FE^{**}(\nabla)$ is the equation of the Lie algebra of infinitesimal transformations $aut(M,\nabla)$. Consider a system of local affine coordinate functions of 
$(M,\nabla)$
$$x: = (x_1,..,x_m).$$ 
let
$$ X(x) = \sum_j X_j(x)\frac{\partial}{\partial x_j}.$$
The SPDE which is derived from $FE^{**}(\nabla)$ is reduced to the linear system
$$FE^{**}_{ij:q}(\nabla):\quad \frac{\partial^2 X_q}{\partial i\partial j} = 0.$$
Thus a local solution $X$ to $FE^{**}(\nabla)$ has the form
$$X(x) = \sum_q[\sum_{j}A_{[q,j]}x_j + c_q]\frac{\partial}{\partial x_q}$$
Therefore the local solutions depend affinely of the point $x$.\\
This proves that at every point $x \in M$ one has the following equality
$$ (***):\quad dim(T_\nabla M(x)) = dim(M)$$
Thereby we get
$$ r^b(M) = 0.$$
This ends the demonstration of Theorem \ref{flatgeometry} $\clubsuit$\\
\textbf{Some corollaries of EXTH1}.\\
The demonstration of Theorem \ref{flatgeometry} is useful for hightlighting some intersting facts. We focus on the global analysis of the $FE^*(\nabla)$, \cite{Kumpera-Spencer} \\
\begin{cor} Let $(M,\nabla)$ be a locally flat manifold. Then the fundamental equation $FE*(\nabla)$ is a transitive Lie equation $\clubsuit$
\end{cor}
\textbf{An idea of proof}. Let $x: = (x_1,..,x_m)$ be a system of local affine coordinate functions of $(M,\nabla)$. We put
$$X = \sum_k X^k\frac{\partial}{\partial x_k}.$$
Then the first fundamental equation $ FE*(\nabla)$ is equivalent to the following of SPDE
$$\Theta^k_{ij} = 0.$$
Here
$$\Theta^k_{ij} = \frac{\partial^2X^k}{\partial x_i\partial x_j}$$
Therefore all of the functions $X^k(x_1,..,x_m)$ depend affinely on $\left\{x_1,..,x_m\right\}$,viz
$$ X(x) = [A].x + [c] \clubsuit$$
The classical machinery of R. Palais leads to the following statement
\begin{cor}. Let $(M,\nabla)$ be a locally flat manifold.\\
(1) The Associative Algebras Sheaf $\mathcal{J}_\nabla$ is the infinitesimal counterpart of the action of a bi-invariant affine Cartan-Lie $(\Gamma,\tilde{\nabla})$. \\
(2) The bi-invariant affine Cartan-Lie group $(\Gamma,\tilde{\nabla})$ is the localization of the simply connected bi-invariant affine Lie group $(G^0,\nabla^0)$ whose associative algebra is isomorphic to $(J_\nabla,\nabla)$.\\
(3) If $(M,\nabla)$ is geodesically complete then the sheaf $\mathcal{J}_\nabla$ is the infinitesimal counterpart of a locally free action of the simply connected bi-invariant affine Lie group $(G^0,\nabla^0)$ $\clubsuit$
\end{cor}
\subsection{DEMONSTRATION OF THEOREM \ref{hessiangeometry1}}
\textbf{Reminder EXTH2: The goal is to prove that $r^b(M,g) = 0$ if and only if the $m$-dimensional Riemannian manifold $(M,g)$ admits Hessian structures $(M,g,\nabla)$.}\\
\textbf{We go to prove that (1) implies (2)}
For our purpose we recall that every Riemannian manifold $(M,g)$ defines the Amari functor
$$\mathcal{LF}(M)\ni \nabla \rightarrow \nabla^g\in \mathcal{LC}(M)$$
where the Koszul connection $\nabla^g$ is 
$$ g(\nabla^g_XZ,Y) = Xg(Y,Z) - g(\nabla_XY,Z).$$
We assume that
$$ r^b(M,g) = 0.$$
There exists $\nabla \in \mathcal{LF}(M)$ such that
$$r^b(\nabla^g) = dim(M)$$
We aim in showing that $\nabla^g \in \mathcal{LF}(M)$. Both $\nabla$ and $\nabla^g$ are flat, viz 
$$R^\nabla = 0,$$
$$R^{\nabla^g} = 0.$$
Therefore it remains to show that the connection $\nabla^g$ is symmetric.\\
The sheaf $\mathcal{J}_{\nabla^g}$ is transitive. We choose a local base of sections of $\mathcal{J}_\nabla^g$
$$\left\{X_1,..,X_m\right\}$$
The flatness of $\nabla^g$ is equivalent to the identity
$$\nabla^g_{X_i}(\nabla^g_{X_j}X_k) - \nabla^g_{X_j}(\nabla^g_{X_i}X_k) = \nabla^g_{[X_i,X_j]}X_k $$
Since $\mathcal{J}_{\nabla^g}$ is an Associative Algebras Sheaf we have the identity
$$ \nabla^g_{X_i}(\nabla^g_{X_j}X_k) - \nabla^g{X_j}(\nabla^g_{X_i}X_k) = \nabla^g_{[X_i,X_j]_{\nabla^g}}X_k $$
We recall that both $R^{\nabla^g}(X_i,X_j)X_k$ and $KV_{\nabla^g}(X_i,X_j,X_k)$ are
$C^\infty(M)$-multi-linear. Therefore one gets the identity
$$ \nabla^g_{[X,Y] - [X,Y]_{\nabla^g}} = O.$$
A consequence of the last identity is that $\nabla^\star(g)$ is symmetric.\\
Thus the quadruple $(M,g,\nabla,\nabla^g)$ is a dually flat pair in the sense of \cite{Amari-Nagaoka}. Consequently both $(M,g,\nabla^*)$ and $(M,g,\nabla^*(g))$ are Hessian manifolds.\\
\textbf{We go to prove that (2) implies (1)}.\\
We assume that there exists  $\nabla \in \mathcal{LF}(M)$ such that $(M,g,\nabla)$ is a Hessian manifold. Thereby $(M,\nabla)$ is a locally flat manifold . Its KV algebra is denoted by
$\mathcal{A}$. Thus our starting assumption means that the metric tensor $g$ is a scalar 2-cocycle of the KV complex  $C_{KV}(\mathcal{A},\mathbb{R})$. This last assumption implies that $\nabla^g$ is torsion free. Therefore the couple $(M,\nabla^g)$ is a locally flat manifold. Thereby
$$r^b(\nabla^g) = dim(M)$$
In final we get
$$r^b(M,g) = 0.$$
Theorem EXTH2 is demonstrated $\clubsuit$
\subsection{DEMONSTRATION OF THEOREM \ref{hessiangeometry2}.}
\textbf{Reminder EXTH3: We are to show that $r^b(M,\nabla)$ vanishes if and if $(M,\nabla)$ admits  Hessian structures.}\\
\textbf{We go to prove that (1) implies (2)}\\
Owning $(M,\nabla)$ We use the functor of Amari for defining the map of $\mathcal{R}ie(M)$ in 
$\mathcal{LC}(M)$ 
$$\mathcal{R}ie(M)\ni g\rightarrow \nabla^g\in \mathcal{LC}(M),$$
here $\nabla^g$ is defined by
$$ g(\nabla^g_XZ,Y) = Xg(Y,Z) - g(\nabla_XY,Z).$$
Since 
$$r^b(M,\nabla) = 0$$
there exists $g^* \in \mathcal{R}ie(M)$ such that
$$r^b(\nabla^{g^*}) = dim(M).$$
Therefore the Associative Algebras Sheaf $\mathcal{J}_{\nabla^{g^*}}$ is transitive.
We put
$$K(X,Y) = \nabla^{g^*}_X\circ\nabla^{g^*}_Y - \nabla^{g^*}_Y\circ\nabla^{g^*}_X$$
and we involve the calculations  made in the demonstration of Theorem EXTH1, namely
$K.1:\quad K(X,Y) = \nabla^{g^*}_{[X,Y]},$\\
$K.2:\quad K(X,Y) = \nabla^{g^*}_{\nabla^{g^*}_XY - \nabla^{g^*}_YX}.$\\
From the properties K.1 and K.2  we deduce that $\nabla^{g^*}$ is torsion free. Thereby $(M,g,\nabla,\nabla(g))$ is a dually flat pair. In final both $(M,g^*,\nabla)$ and $(M,g^*,\nabla^{g^*})$ are Hessian manifolds.\\
\textbf{We go to prove that (2) implies (1)}.\\
We assume that a locally flat manifold $(M,\nabla)$ admits a Hessian structure
$$(M,g,\nabla).$$
Then we define the Koszul connection $\nabla^*$ by
$$g(\nabla^*_XZ,Y) = X(g(Y,Z) - g(\nabla_XY,Z).$$
Let us consider the KV algebra of $(M,\nabla)$
 $$\mathcal{A} = (\mathcal{X}(M),\nabla).$$
At one side the metric tensor $g$ is a 2-cocycle of the KV complex $C(\mathcal{A},\mathbb{R})$.\\
At another side we know that the identity
$$\delta_{KV}g = 0$$
holds if and only if $\nabla^*$ is torsion free. Since $\nabla^*$ is flat $(M,\nabla^*)$ is a locally flat manifold. According to Theorem \ref{flatgeometry} one has
$$r^b(\nabla^*) = dim(M)$$
The last conclusion yields the equality
$$r^b(M,\nabla) = 0.$$
Theorem \ref{hessiangeometry2} is demonstrated $\clubsuit $\\
\subsection{DEMONSTRATION OF THEOREM \ref{hessiangeometry3}.}
\textbf{Digression. Actually it is easy to see that Theorem \ref{hessiangeometry3} namely EXTH4, may be regarded either as a corollary of Theorem \ref{hessiangeometry1} or as a corollary of Theorem \ref{hessiangeometry2}.}\\
At one side $ r^B(M) = 0$ is equivalent to the existence of a Riemannian metric tensor
$$g\in \mathcal{R}ie(M)$$
which is a solution to 
$$r^b(M,g) = 0.$$
At another side $ r^{B^*} = 0$ is equivalent the existence of a locally flat structure 
$$\nabla \in \mathcal{LF}(M)$$ 
which is a solution to
$$r^b(M,\nabla) = 0.$$
Those observations yield Theorem 9.4 $\clubsuit$
\subsection{DEMONSTRATION OF THEOREM \ref{informationgeometry1}.}
\textbf{Reminder EXTH5: Let $\mathbb{M} = [\mathcal{E},\pi,M,D,p]$ be a statistical model for a measurable set $(\Xi,\Omega)$. The Fisher information of $\mathbb{M}$ is denoted by $g$.\\
We recall that the complexity problem of $\mathbb{M}$ is the search of invariants informing about how far from being an exponential family is the model $\mathbb{M}$. Our purpose is to demonstrate that this problem is linked with the function $r^b$ through the combitorial problem 
$DL(\mathcal{S}) \clubsuit $}\\
\textbf{We go to prove that (1) implies (2)}.\\
The starting assumption is that $\mathbb{M}$ is a regular model with
$$r^b(M,g) = 0.$$
We apply Theorem \ref{flatgeometry} to conclude that there exists a locally flat structure $(M,\nabla)$ such that $(M,g,\nabla)$ is a Hessian structure. Let $\mathcal{A}$ be the KV algebra of $(M,\nabla)$. Since $(M,g,\nabla)$ is a Hessian manifold $g$ is a 2-cocycle of the scalar complex $C^*(\mathcal{A},\mathbb{R})$.\\
The property $\delta_{KV}g = 0$ is equivalent the following properties,\\
 there exist a random function
$$\mathcal{E}\ni e\rightarrow a(e)\in \mathbb{R},$$
and there exists a differentible function
$$M \ni x\rightarrow \psi(x)\in\mathbb{R}$$
which are subject to the following requirements\\
 $$\nabla^2(\int_{\mathcal{E}_x}a(e)) = 0,$$
 $$p(e) = exp\left\{a(e) - (\psi\circ\pi)(e)\right\}.$$
For more details the readers are referred to \cite{Nguiffo Boyom(6)}.\\
\textbf{We go to prove that (2) implies (1)}.\\
It is well known that the Fisher information $g$ of a regular exponential model is a Hessian Riemannian metric, (see the e-m-flatness as in Amari-Nagaoka \cite{Amari-Nagaoka}). Therefore we are in position to conclude that
$$r^b(M,g) = 0\quad\clubsuit$$
\subsection{DEMONSTRATION OF THEOREM \ref{informationgeometry2}}
\textbf{Reminder EXTH6:
In Theorem \ref{informationgeometry2}  We go to deal with singular models $\mathbb{M}$. The concern is to show that $r^b(\mathbb{M})$ vanishes if and only $\mathbb{M})$ is an exponential family}.\\
The readers have to notice that the framework to address Theorem \ref{informationgeometry2} is the differential topology. Indeed let $g_\mathbb{M}$ be the Fisher information of a singular model $\mathbb{M}$. The distribution $Ker(g_\mathbb{M})$ is in involution, viz the Poisson bracket $[X,Y]$ of two sections of $Ker(g_\mathbb{M})$ is a section of $Ker(g_\mathbb{M})$. Further if $X$ is section of $Ker(g_\mathbb{M})$ then
$$L_Xg_\mathbb{M} = 0.$$
Thus the statistical geometry of $\mathbb{M}$ is transverse to the distribution $Ker(g_\mathbb{M})$.\\
For convenience we suppose to be dealing with a singular analytic model $\mathbb{M}$. Therefore the triple $(M,Ker(g_\mathbb{M}),g_\mathbb{M})$ is stratified Riemannian foliation. The family of $\alpha$-connections $\nabla_\alpha$ agrees with the stratification of $M$ by the analytic singular foliation $(M,Ker(g_\mathbb{M}))$.\\
\textbf{Reminder}\\
We go to use the stratification of $M$ defined by $Ker(g_\mathbb{M})$. Loosely speaking we deal with the filtration of $M$ by closed analytic submanifolds
$$ \subset F_{j+1}\subset F_j\subset F_{j-1}\subset..\subset F_1\subset M.$$
At the level $j$ let $g_j$ be the restriction to $\left\{F_{j-1} - F_j\right\}$ of
$g_\mathbb{M}$. Then the triple
$$\left\{F_{j-1} - F_j, Ker(g_j),g_j)\right\}$$
is regular Riemannian foliation as in \cite{Reinhardt}, \cite{Molino}.\\
By the way the model $\mathbb{M}$ is endowed with the filtration defined by setting\\
$\mathbb{M}_j = [\mathcal{E}_{F_j}, \pi, F_j, p|\mathcal{E}_{F_j}].$\\
Since the family $\nabla^\alpha$ agrees with the stratification Theorem \ref{informationgeometry2} also agrees with the filtration of $[\mathcal{E},\pi,M,p]$ by the family $\mathbb{M}_j$.\\
To pursue it is sufficient to prove Theorem \ref{informationgeometry2} when $(M,Ker(g_\mathbb{M}),g_\mathbb{M})$ is a regular Riemannian foliation. Thereby the assumption
$$r^b(\mathbb{M}) = 0$$
yields the exitence of $\alpha \in \mathbb{R}$ such both $(M,\nabla^\alpha)$ and
$(M,\nabla^{-\alpha})$ are locally flat manifolds. Therefore the triple $(g_\mathbb{M},\nabla^\alpha, \nabla^-{\alpha})$ is transversely a dually flat pair. Consequently the Fisher information $g_\mathbb{M}$ is a KV cocycle of both $(M,\nabla^\alpha)$ and $(M,\nabla^{-\alpha})$.\\
Locally, viz in a local statistical model $ \mathbb{M} = (\Theta,p)$ we get $\delta_{KV}g_\mathbb{M} = 0$ if and only if the probability density $p$ is an exponential family. Thus a global statistical model whose Fisher information is a KV cocycle is locally exponential. Since a (real) exponential function is one-to-one we conclude that
$\delta g_\mathbb{M} = 0$ if and only if $p$ is an exponential family. This end the demonstration of Theorem \ref{informationgeometry2} $\clubsuit$
\section{THE GEOMETRIC COMPLETENESS OF LOCALLY FLAT MANIFOLDS}
(1) An important problem in the affine geometry is the geometric completeness. This probem is whether the developing map of a $m$-dimensional locally flat manifold is a diffeomorphism onto the affine euclidean space $(\mathbb{R}^m,D)$. \\
(2): Another problem is whether a given a virtually polycyclic group can be the fundamental group of a geometrically complete compact locally flat manifold. This is known as a problem of Milnor \cite{Milnor}, \cite{Auslander-Kuranishi}, \cite{Auslander-Markus}. Among references are \cite{FGH}, \cite{Milnor}, \cite{Carriere}.\\
(3): A conjecture of Markus claims the following,\\
\textbf{Every compact locally flat manifold whose linear holonomy group is unimodular is geometrically complete $\clubsuit$}\\
The conjecture of Markus has been (and remains) an outstanding open conjecture in the locally flat geometry. It has been proved in many particular cases, see \cite{Carriere}, \cite{FGH}, \cite{Kim} and references therein.\\
Another outstanding result in the locally flat geometry is 
\begin{thm}(Benzecri) The only compact orientable surface supporting geometrically complete locally flat structures is the flat torus $$ \mathbb{T}^2 = \frac{\mathbb{R}^2}{\mathbb{Z}^2}\clubsuit$$
\end{thm}
We aim at linking the completeness problem with the function $r^b$. The framework for doing that is the notion of Special Dynamical Systems that we have introduced. This approach allows to reformulate some general criterion for a geodesically complete locally flat manifold being geometrically complete. In our approach the geometric completeness of geodesically complete locally flat manifolds is replaced  by the geometric completeness of bi-invariant affine Lie groups. The return of this approach is that the geoemtric completeness of finite dimensional affine Lie groups can be expressed in terms of affine representations of Lie algebras \cite{Milnor}, \cite{Nguiffo Boyom(8)}. More than forty years ago we have been concerned with affine embeddings of finite dimensional real Lie groups. The purposes are affine representations of a $q$-dimensional Lie group having a $q$-dimensional orbit which lies in a $q$-dimensional affine subspace \cite{Nguiffo Boyom(9)}. We go to perform a similar idea in this section. \\
Readers interested in a survey on the completeness problem are refered to W.M. Goldman \cite{Goldman} 
\subsection{The completeness of affine Lie groups and affine representations}
We go to recall the approach by affine representations.\\
Let $\left\{H < G\right\}$ be a pair formed of a connected closed Lie subgroup $H$ of a connected Lie group $G$. We consider the action of $G$ by inner autormorphisms
$$G\times G \ni(g,g^*)\rightarrow gg^*g^{-1}\in G.$$ 
The orbit of $H$ under this inner action 
$$[H] = \left\{gHg^{-1},\quad g \in G\right\}$$
\textbf{Reminder}\\
We recall some notions already introduced. \\
(1) A pair $\left\{[H] < G\right\}$  is called effective if $H$ does not contain any nontrivial normal subgroup of $G$.\\
(2) A pair $\left\{[H] < G\right\}$ is called simply connected if $G$ is simply connected.\\
(3) A pair $\left\{[H] < G\right\}$ is a pair of bi-invariant affine Lie groups if $H$ is a bi-invariant affine subgroup of the bi-invariant affine Lie group $G$.
\begin{prop}
The category of simply connected finite dimensional effective pairs of bi-invariant affine Lie groups, namely
$$\left\{[H] < G\right\}$$
is equivalent to the category of finite dimensional effective pairs of associative algebras
$\left\{[I] < J\right\}$ $\clubsuit$
\end{prop}
\subsubsection{The geometric completeness of abstract algebras}
\begin{defn} A real algebra is a pair $(A,\mu)$ formed of a real vector space $A$ and a bi-linear map
$$A\times A \ni (a,a^*)\rightarrow \mu(a,a^*)\in A.$$
\end{defn}
Given a real algebra $(A,\mu)$ we go to make simple by putting
$$a.a^* = \mu(a,a^*).$$
To every $a^* \in A$ we assign the linear endomorphism $\psi_{a^*}$ defined by
$$\psi_{a^*}(a) = a.a^* + a.$$
\begin{defn} A finite dimensional abstract algebra $A$ is called geometrically complete if all of the linear endomrphisms 
$\psi_{a^*}$ are injective $\clubsuit$
\end{defn}
Actually every subalgebra of a geometrically complete algebra is geometrically complete too. Therefore an effective pair of associative algebras $\left\{[I] < J\right\}$ is called complete if $J$ is complete.\\
Let $(G,\nabla)$ be a connected affine Lie group whose (left invariant) KV algebra is denoted by 
$\mathcal{G}_\nabla$. \\
Actually the KV algebra of $G$ is the real algebra
$$\mathcal{G}_\nabla = (\mathcal{G},\nabla),$$
here $\mathcal{G}$ is the Lie algebra of the Lie group $G$. Note that $1_\mathcal{G}$ stands for the identity endomorphism of the vector space $\mathcal{G}$.\\
\subsubsection{The canonical affine representation of affine Lie groups} 
We define the affine representation of Lie algebra $\mathcal{G}$ in the vector space $\mathcal{G}$ by 
$$ \mathcal{G}\ni a\rightarrow (\nabla,1_\mathcal{G})(a) = (\nabla_a, a) \in Hom_\mathbb{R}(\mathcal{G},\mathcal{G})\times \mathcal{G}.$$
The affine map $(\nabla_a,a)$ is defined by
$$(\nabla_a,a)(a^*) = \nabla_aa^* + a.$$
The universal covering of $(G,\nabla)$ is denoted by $(\tilde{G},\tilde{\nabla})$. The affine representation $(\nabla,1_\mathcal{G})$ is the differential at the unit element of a unique affine representation $(f,q)$ of $\tilde{G}$ in $\mathcal{G}$, namely  
$$ \tilde{G} \ni g\rightarrow (f(g),q(g))\in Gl(\mathcal{G})\times \mathcal{G}.$$
Thus at the level of the Lie algebra of $G$ the representaion  $(\nabla,1_\mathcal{G})$ is the linear counterpart of the affine representaion $(f,q)$. To every $a^* \in \mathcal{G}$ is assigned the orbital  map $\Psi_{a^*}$ 
$$ \tilde{G}\ni g\rightarrow \Psi(g) \in \mathcal{G},$$
here $\Psi$ is defined by
$$\Psi_{a^*}(g) = f(g)(a^*) + q(g)$$
\begin{defn} Keep the notation above.\\
(1) An affine Lie group $(G,\nabla)$ is called geometrically complete if at least one orbital maps $\Psi_{a^*}$ is surjective.\\
(2) An affine Lie group $(G,\nabla)$ is called hyperbolic if all of the $\Psi_{a^*}(\tilde{G}$ are convex domains not containing any straight line $\clubsuit$
\end{defn}
\begin{thm} We keep the notation $(G,\nabla)$,  $\mathcal{G}_\nabla$. The following assertions are equivalent\\
(1) $(G,\nabla)$ is geometrically complete,\\
(2) $(\mathcal{G}_\nabla,\nabla)$ is geometrically complete $\clubsuit$
\end{thm}
\textbf{Hint: By (2) all of the orbits are open $\clubsuit$}\\
We go to match the completeness problem of geodesicaly complete locally flat manifolds.\\
Let $(\tilde{M},\tilde{\nabla})$ be the universal covering of a finite dimensional geodesically complete locally flat manifold $(M,\nabla)$. The covering map is denoted by $\pi$.
Let $(G_\nabla,\nabla^0)$ be the simply connected bi-invariant affine Lie group whose Lie algebra is
$(J_{\tilde{\nabla}},[-,-]_{\tilde{\nabla}})$. \\
The convering map $\pi$ is a local isomorphism of the AAS $\left\{\mathcal{J}_{\tilde{\nabla}},\tilde{\nabla}\right\}$ in the AAS $\left\{\mathcal{J}_{\nabla},\nabla\right\}$. In particular we have
$$(J_\nabla,\nabla) = \pi_\star\left\{J^{\pi_1(M)}_{\tilde{\nabla}},\tilde{\nabla}\right\}.$$
Here $\pi_1(M)$ is the fundamental group of $M$ and 
$$ J^{\pi_1(M)}_{\tilde{\nabla}} = \left\{\xi \in J_{\tilde{\nabla}}\quad|\quad \gamma_\star(\xi) = \xi \quad \forall \gamma \in \pi_1(M) \right\}.$$
We recall that $(G_\nabla,\nabla^0)$ is the simply connected bi-invariant affine Lie group whose associative algebra is $(J_{\tilde{\nabla}},\tilde{\nabla})$. \\
The locally flat manifold $(\tilde{M},\tilde{\nabla})$ is a homogeneous space of $(G_\nabla,\nabla^0)$.\\
\textbf{Reminder.}\\
\textbf{The fundamental group $\pi_1(M)$ is a normal subgroup of $G_\nabla$ $\clubsuit$} \\
We fix a base point $x^*\in \tilde{M}$. The isotropic subgroups at $x^*$ is denoted by
$G_{\nabla x^*} \subset G_\nabla$. Consider the orbital map
$$ G_\nabla \ni g\rightarrow \Psi_{x^*}(g)\in \tilde{M}.$$
This yields the principal bundle  
$$ G_{\nabla x^*}\rightarrow G_\nabla \rightarrow \tilde{M}.$$
Let $I_{x^*}$ be the Lie algebra of $G_{\nabla x^*}$. It is a simple right ideal of the associative algebra $(J_{\tilde{\nabla}},\tilde{\nabla})$. Then we get the effective pair of associative algebras
$$ \left\{[I_{x^*}] < J_{\tilde{\nabla}}\right\}$$
which is the algebraic counterpart of the pair of bi-invariant affine Lie groups  
$$\left\{[G_{\nabla x^*}] < G_\nabla\right\}.$$
We consider the gauge fibration
$$(G_{\nabla x^*},\nabla^0)\rightarrow (G_\nabla,\nabla^0) \rightarrow (\tilde{M},\tilde{\nabla}).$$
Since both $G_\nabla$ and $\tilde{M}$ are simply connected the Lie subgroup $G_{\nabla x^*}$ is simply connected  as well. The gauge fibration agrees with the developping maps   
$$\Delta_{x^*}: G_{\nabla x^*} \rightarrow  \mathbb{R}^n,$$
$$\Delta: G_\nabla\rightarrow \mathbb{R}^{m+n},$$
$$ \tilde{\Delta}: \tilde{M}\rightarrow \mathbb{R}^m.$$
We observe that the developping maps have many subtle properties. In particular $\tilde{\Delta}$ is a $\left\{\pi_1(M),H(\pi_1(M))\right\}$- equivariant local affine isomorphism in $\mathbb{R}^m$, \cite{Carriere}, \cite{FGH}.\\
Here $H(\pi_1(M))$ is the affine holonomy group of $(M,\nabla)$; it is the image of $\pi_1(M)$ under the holonomy affine representation. For more detail see D. Fried, W Goldman and M. Hirsh,\cite{FGH}.
Since $(M,\nabla)$ is geodesically complete the group $Diff(M,\nabla)$ is a bi-invariant affine Lie group whose universal covering is $(G_\nabla,\nabla^0)$. 
\textbf{Remind. If $M$ is compact then the conjecture of Markus claims that for $(M,\nabla)$ being complete it is sufficient that its linear component of $H(\pi_1(M))$ be unimodular.} \\
To a pointed geodesically complete locally flat manifold $\left\{x^* \in M,\nabla\right\}$ we assign the triple  
$$ \tau(\nabla): = \left\{(G_{\nabla x^*},\nabla^0), (G_\nabla,\nabla^0), (\tilde{M},\tilde{\nabla})\right\}.$$
\begin{thm} (The geometric completeness theorem). In the triple 
$$\tau(\nabla) = \left\{(G_{\nabla x^*},\nabla^0); (G_\nabla,\nabla^0); (\tilde{M},\tilde{\nabla})\right\}$$ 
the completeness of a pair yields the completeness of the third member $\clubsuit$
\end{thm}
\section{ THE INTRINSIC NATURE OF THE FUNCTION $r^b$, continued.}
To motivate the concerns of this sections we go to overview a few purposes discussed in preceding sections.
\subsection{Cominatorial distancelike}
In preceding sections we have involved the function $r^b$ in studying the problem $DL(\mathcal{S})$.\\
At another side we have been concerned with the graph $\mathcal{GR}(LC)$ containing the remarkable subgraphs\\
$(1):$ the subgraph $\mathcal{GR}(LF)$ whose vertices are locally flat connections,\\
$(2):$ The subgraph $\mathcal{GR}(SLC)$ whose vertices are torsion free connections,\\
$(3):$ The subgraph $\mathcal{GR}(HES)$ whose vertices are locally flat connections admitting  Hessian Riemannian metric tensors.\\
$(4):$ The subgraph $\mathcal{GR}(Rie)$ whose vertices are metric connections.\\
\textbf{Reminder: a Koszul connection in $M$ $\nabla$ is called a metric connection if there exists a Riemannian structure $(M,g)$ such that $\nabla_Xg = 0$ for all $X\in \mathcal{X}(M)$}.\\
One can use the restricted holonomy group for discussing the question whether a Koszul connection is metric;( see also the classical Theorem of Ambrose-Singer). In a forthgoing program a concern is to investigate links of metrics with the dynamic 
 $$ \mathit{COX}(M)\times \mathcal{MG}(M)\rightarrow \mathcal{MG}(M).$$
 Here $\mathcal{COX}(M)$ is the group of isomorphisms of $\mathcal{MG}(M)$ generated by the isometries $$(M,\nabla)\rightarrow (M,\nabla^g, g\in \mathcal{R}ie(M).$$
\textbf{We have already widely discussed some challenges $DL(\mathcal{S})$ in $\mathcal{GR}(LC)$;  we aim at enriching those purposes}.\\ 
\begin{prop} We consider the restriction to $\mathcal{GR}(SLC)$ of the function $r^b$. \\
(1) For every vertex $\nabla \in \mathcal{GR}(SLC)$ the integer $dim(M) - r^b(\nabla)$ measures how far from the subgraph $\mathcal{GR}(LF)$ is the vertex $\nabla$,\\
(2) For every $\nabla \in \mathcal{GR}(LF)$ the integer $r^b(M,\nabla)$ measures how far from the subgraph $\mathcal{GR}(HES)$ is the vertex $\nabla$ $\clubsuit$
\end{prop}
\textbf{Hint}. The proposition is a straight corollary of Theorems EXTH.\\
The prposition emphasizes the distancelike nature of some restrictions of the function $r^b$. \\
From the naïve classical analysis perspective one poses 
$$[r^b|_{\mathcal{GR}(SLC)}]^{-1}(0) = \mathcal{GR}(LF),$$
here $r^b|_{\mathcal{GR}(SLC)}$ is the restriction to $\mathcal{GR}(SLC)$ of the function $r^b$.
\subsection{The fundamental equation $FE^*(\nabla)$ and characteristic obstructions}
We go to focus on the restriction 
$$dim(M) - r^b|_{\mathcal{GR}(SLC)}.$$
The lower bound of this restriction is the numerical invariant $r^b(M).$ \\
By Theorems EXTH the lower bound of $dim(M) - r^b|_{\mathcal{GR}(SLC)}$ is interpreted as a characteristic obstruction. That is the deep meaning of EXTH1.  
So we get the following equivalence 
$$ r^b(M) > 0\leftrightarrow \mathcal{GR}(LF) = \emptyset.$$
\textbf{A few examples}\\
\textbf{Example.1} For every euclidean flat torus $\mathbb{T}^m$ 
$$r^b(\mathbb{T}^m) = 0.$$
\textbf{Example.2} If $(M,\nabla)$ is a special compact orientable surface satisfying
$$r^b(S) > 0 $$
then the Euler characteristic of $S$ is zero.
\section{THE FUNDAMENTAL EQUATION $FE^*(\nabla)$ AND THE RIEMANNIAON GEOMETRY}
\textbf{A digression.}\\
Let $(M,g)$ be a Riemannian manifold. The convex set of Riemannian conections is denoted by $\mathcal{LC}(M,g)$. The Levi-Civita connection $\nabla$ is the unique torsion free element of 
$\mathcal{LC}(M,g)$. We go to sketch some ideas proving that non symmetric metric connections deserve the attention.\\
Given $\nabla \in \mathcal{MC}(M,g)$ the differential operators $D^\nabla$ and $D_\nabla$ coincide if and only if $\nabla$ is the Levi-Civita connection. Otherwise for non symmetric metri conection $nabla$ solutions to $FE^*(\nabla)$ differ from solutions to $FE^{**}(\nabla)$. In the next subsection we go to focus on the differential equation $FE^*(\nabla)$.\\
\subsection{New geometric invariants of Riemannian manifolds}
We deal with a generic non symmetric metric connection $\nabla) \in \mathcal{MC}(M,g)$. Our concern is the operator $D_\nabla$ which products an inductive system of finite dimensional Associative Algebras Sheaf and a projective system of finite dimensional Associative Algebras Sheaf 
$$\left\{(\mathcal{J}_{\nabla^q},\nabla^q) |\quad : I^p_q; \Pi^q_p , p\leq q\right\}$$
The vector spaces of sections of those sheaves yield an inductive system of finite dimensional associative algebras and a projective system of finite dimensional associative algebras
$$\left\{(J_{\nabla^q},\nabla^q) |\quad: I^p_q;\Pi^q_p, p\leq q \right\}$$
Of course those systems are the algebraic versus of an inductive system (respectively of a projective system) of finite dimensional simply connected bi-invariant affine Lie groups.\\
At the algebraic versus the Initial Object of that inductive system (respectively the Final Object of that projective system) is the Associative Algebras Sheaf $(\mathcal{J}_\nabla,\nabla)$.\\
In the Lie theory versus the Initial Object (respectively the Final Object) of those systems is the simply connected bi-invariant affine Lie group $(G_\nabla,\nabla^0)$ the associative algebra of which is $(J_\nabla,\nabla)$.\\ 
We denote the Sheaf of differentiable vector fields by $\mathcal{T}M$; it is a right module of the AAS $\mathcal{J}_\nabla$. The right action is defined by 
$$ X.\xi = \nabla_X\xi \quad \forall (X,\xi)\in \mathcal{X}(M)\times J_\nabla.$$
Thereby $\mathcal{T}M$ is a right module of both the KVAS $(\mathcal{J}_\nabla,\nabla)$ and of $(\mathcal{J}_\nabla, [-,-]_\nabla)$.\\
As already seen we face three sheaves of cochain complexes\\
(1): the Hochschild complex of the AAS $(\mathcal{J}_\nabla,\nabla)$ 
$$C_H(\mathcal{J}_\nabla,\mathcal{T}M)$$,\\
(2): the Chevalley-Eillenberg complex of $(\mathcal{J}_\nabla, [-,-]_\nabla)$ yields the sheaf
$$C_{CE}(\mathcal{J}_\nabla,\mathcal{T}M),$$
(3): the KV complex of the KVAS $(\mathcal{J}_{\nabla},\nabla)$
$$ C_{\tau}(\mathcal{J}_\nabla,\mathcal{T}M).$$
\textbf{The triple formed of the complexes justed described is denoted by $\mathcal{TR}_\nabla(M,g)$}.
\begin{defn} We put
$$\mathcal{TR}(M,g) = \left\{\mathcal{TR}_\nabla |\quad \nabla \in \mathcal{MC}(M,g)\right\}$$
\end{defn}
\textbf{$\mathcal{TR}(M,g)$ is a geometric invariant of the Riemannin manifold $(M,g)$}
\subsection{Deformation of metric connections} 
A pioneering work on the deformations of metric connections is \cite{Raghunathan}. The methods S. Raghunathan are essentially topological-analytic. The subgrah $\mathcal{GR}(M,g)$ of metric connections is convex. Therefore the classical analysis and the topology may provide tools for exploring neighbourhoods of a metric connection.\\
From the viewpoint of the algebra structure $(\mathcal{X}(M),\nabla)$ one would be interested in an analogous of the theory of deformation of algebraic structures as in Piper \cite{P   iper}. At the moment every Koszul connection $\nabla$ give rise to three Algebras Sheaves we just decribed in the last subsection, namely \\
(a) the AAS $(\mathcal{J}_\nabla, \nabla)$ whose deformation is controlled by the Hochschild cohomology,\\
(b) the KVAS $(\mathcal{J}_\nabla,\nabla)$ whose deformation is controlled by the KV cohomology,\\
(c) the sheaf $(\mathcal{J}_\nabla,[-,-]_\nabla)$ whose deformation is controlled by the Chevalley-Eilenberg cohomology. Thereby depending on concerns and on needs one could study the deformation of $J_\nabla$ as an associative algebra, as a Koszul-Vinberg algebra or as a Lie algebra. We decide not to discuss that general matter in this paper. However we go to focus on the case of transitive AAS. We recall that a Koszul connection $\nabla$ is called transitive if its Associative Algebras Sheaf $\mathcal{J}_\nabla$ is transitive. So we have 
$$ dim(\mathcal{T}_\nabla M(x)) = dim (M) \quad \forall x \in M.$$
Therefore $$span_{C^\infty(M)}(J_\nabla) = \mathcal{X}(M)$$
We have ready seen that every transitive connection yields a deformation of the Poisson bracket of $\mathcal{X}(M).$\\
Obviously transitive Riemannian Levi-Civita connections  are flat. If $(M,g)$ is positive then the  flatness of the connection of Levi-Civita implies that up to isometries one has 
$$(M,g) = (\frac{\mathbb{T}^{b_1(M)}}{K}\times\mathbb{R}^{m-b_1(M)},g_0).$$
Here $g_0$ is the Euclidean (flat) metric of $\mathbb{R}^m$, $b_1(M)$ is the first Betti number of $M$, $K$ is a finite group of isometries of $g_0$, \cite{Wolf}.\\
 \subsection{The differential topology of geodesically complete special Riemannian manifolds, continued}
Let $(M,g,\nabla)$ be a triple formed of a geodesically complete Riemannian manifold $(M,g)$ and its Levi-Civita connection $\nabla$. We can put
$$\mathcal{T}_\nabla M(x) = \frac{J_\nabla}{I_x(M)},$$
here $I_x(M)$ is the ideal of functions that vanish at $x \in M$.\\
We focus on special $(M,g,\nabla)$ whith the property
$$ dim(\mathcal{T}_\nabla M(x)) = constant$$
It is so if a Riemannian manifold $(M,g)$ is homogeneous. Therefore we put
$$r^b(\nabla) = dim(\mathcal{T}_\nabla M(x)).$$
The orbits of $\mathcal{J}_\nabla$ are leaves of a regular foliation of $(M,g)$, namely 
$\mathcal{F}$. Further those leaves are auto-parallel flat submanifolds of ($M,\nabla)$. 
We go to focus on machineries already constructed, namely  \\
(1): the semi-inductive system of finite dimensional simply connected bi-invariant affine Lie groups
$$ \left\{I^p_q: (G_p,\nabla^p)\rightarrow (G_q,\nabla^q)| \quad p\leq q\right\}, $$
(2): the semi-projective system of simply connected bi-invariant affine Lie groups
$$ \left\{ \pi^q_p: (G_q,\nabla^q)\rightarrow (G_p,\nabla^p)|\quad p\leq q \right\}$$
\textbf{Reminder: digressions}.\\
We recall that $G_q$ is the simply connected Lie group whose Lie algebra is the commutator Lie algebra of the associative algebra $(J_{\nabla^q},\nabla^q)$. The Intial Object of this semi-inductive system of bi-invariant affine Lie groups is the locally flat foliation $\mathcal{F}$. \\
We have introduced this construction in every gauge structure $(M,\nabla)$. When a metric connection $\nabla^* \in \mathcal{LC}(M,g)$ differs from Levi-Civita connection of $(M,g)$ the equation $FE^*(\nabla^*)$ is not a Lie equation in $M$. Therefore the AAS $\mathcal{J}_\nabla$ is not tangent to a foliation. However for every $q > 0$ $FE^*(\nabla^q)$ is a Lie equation in the bi-invariant affine Lie group $(G^{q-1},\nabla^{q-1})$. According to the preceding subssection, we know that a Riemannian manifold $(M,g)$ admits many plenty of geometric invariants which are derived from its non symmetric Riemannian connections. In the case of positive signatures the leave of $\mathcal{F}$ many be degenerate. We do not know whether those invariants impact the topology of $M$. Nevertheless they impact the geometry of the Riemannian manifold $(M,g)$; here is an instance
\begin{thm}
Let $Iso(M,g)$ be the group of isometries of a special Riemannian manifold $(M,g)$. There exists a canonical homomorphism of $Iso(M,g)$ in the group of affine isomorphisms of the semi-inductive (respectively semi-projective) system of simply connected bi-invariant affine Lie groups 
$$\mathcal{IS}: =\quad \left\{I^p_q:\quad (G_p,\nabla^p)\rightarrow(G_q,\nabla^q)\right\},$$ 
$$\mathcal{PS}: =\quad \left\{\pi^q_p:\quad (G_q,\nabla^q)\rightarrow (G_p,\nabla^p)\right\} \clubsuit$$
\end{thm}
\textbf{Reminder. A gauge structure $(M,\nabla)$ is called special if the AAS $\mathcal{J}_\nabla$ is non trivial. There exist non special gauge structres; an instance is the +connection of Cartan in a perfect Lie group} 
\subsection{Special complete positive Riemannian manifolds, continued}
In this subsection we go to focus on the category of geodesically complete positive Riemannian manifolds.\\
We consider a special complete triple $(M,g,\nabla)$ formed of a complete positive Riemannian manifold $(M,g)$ and its Levi-Civita connection $\nabla$. Here the rank of distribution $\mathcal{T}_\nabla M$ is constant. The Associative Algebras Sheaf $\mathcal{J}_\nabla, \nabla)$ is tangent to the auto-parallel locally flat foliation 
$$ \left\{\mathcal{F}_\nabla,\nabla\right\}.$$
Every leaf of $\mathcal{F}_\nabla$ inherits the flat Riemannian structure  
$$\left\{F, g |_F, \nabla\right\}$$ 
Up to finite covering a leaf $F$ is a tube over a twisted Euclidean torus Viz
$$ F = \frac{\mathbb{T}^k}{K}\times \mathbb{R}^{r^b(\nabla)-k}$$
here $K$ is finite group of isometries of the Euclidean space $\mathbb{R}^k$, \cite{Wolf}.\\
At one side every leaf $F$ is the domain of the Special Dynamical System
$$\left\{(G_\nabla,\nabla^0)\times (F,\nabla)\rightarrow (F,\nabla)\right\}.$$
At another side every triple $\left\{F,g|_F,\nabla\right\}$ is obviously a Hessian manifold. The action of $G_\nabla$ in $F$ preserves the gauge structure $(F,\nabla)$ but is not an isometric action in $(F,g|_F)$. Nevertheless we go to use the following definition
\begin{defn} We keep the notation just used.\\
(1): Every quintuple 
$$\left\{(\Gamma_\nabla,\nabla^0),(F,g|_F,\nabla)\right\}$$
is called a Geodesic Flat Hessian Special Dynamical System,(GFHSDS in short).\\
(2): The family
$$\mathcal{F}_M = \cup \left\{(\Gamma_\nabla,\nabla^0),(F,g|_F,\nabla)\right\}$$
is called a GFHSDS-foliation of $(M,g,\nabla)$ $\clubsuit$
\end{defn}
Thus in the category of geodesically complete special positive Riemannian manifolds the function $r^b$ has a topological nature which is linked with $EXF(\mathcal{S})$
\begin{prop} In a complete special positive Riemannian triple $(M,g,\nabla)$ the positive integer  $r^b(\nabla^*)$ is the optimal dimension for GFHSDS-foliation in $(M,g,\nabla)$ $\clubsuit$
\end{prop}
\subsection{GFHSDS-foliation of statistical models}
We aim at performing GFHSDS-foliations in the information geometry and in statistical manifolds. To motivate our purpose we go to recall some preliminary notions.\\
Let $(\Xi,\Omega)$ be a measurable set and let $\Gamma$ be the group of measurable isomorphisms of $\Xi$.
We use the theory of statistical models as in \cite{Nguiffo Boyom(6)}. Consider an $m$-dimensional regular statistical model for $(\Xi,\Omega)$, namely 
$$\mathbb{M} = [\mathcal{E},\pi,M,D,p]$$ 
\begin{defn} A statistical foliation in $\mathbb{M}$ is a $\nabla$- auto-parallel foliation $\mathcal{F}$, here $\nabla$ is the Levi-Civita connection of the Fisher information $\clubsuit$
\end{defn}
\begin{defn}
A statistical model $\mathbb{M}$ is called special-complete if all of its $\alpha$-connections are special and geodesically complete.
\end{defn}
We go to apply the GFHSDS-machinery to regular complete special statistical models. 
\begin{thm} Consider an $m$-dimensional regular special-complete statistical model $\mathbb{M}$. The Levi-Civita connection of the Fisher information is denoted by $\nabla^0$. Set
$$ q = r^b(\nabla^0),$$
 $$\mathbb{M} = [\mathcal{E},\pi,M,D,p].$$
(1) $\mathbb{M}$ admits a $q$-dimensional GFHSDS-foliation
$$\mathcal{F}_\mathbb{M} = [\mathcal{E}_\mathcal{F},\pi,\mathcal{F},p]$$
which is defined by the orbits of the simply connected bi-invariant affine Cartan-Lie group 
$$ \left\{\Gamma_0, \tilde{\nabla}^0\right\}$$ 
the  Associative Algebras Sheaf of which is 
$$\left\{\mathcal{J}_{\nabla^0},\nabla^0\right\}$$ \\
(2) The Cartan-Lie group $\Gamma^0$ is the localization of the differentible action of the simply connected bi-invariant affine Lie group $(G_{\nabla^0},\nabla^0)$ whose associative algebra is 
$$\left\{J_{\nabla^0},\nabla^0\right\}\clubsuit$$
\end{thm}
\textbf{$Hint_1$}\\
(1): For every leaf $F$ of $J_{\nabla^0}$, the model whose base manifold is $F$ is denoted by
$$\mathbb{M}_F = [\mathcal{E}_F,\pi,F,\nabla^0,p].$$
Every triple $(F,g|_F,\nabla^0)$ is a flat Hessian manifold.\\
(2): According to (1) every model $\mathbb{M}_F$ is an exponential family for $(\Xi,\Omega)$. Thus along every leaf $F$, $\mathbb{M}_F$ is an optimal model for $(\Xi,\Omega)$ $\clubsuit$.\\
The remarks above are particularly interesting in the applied information geometry.\\
\textbf{$Hint_2$}\\
Now we go to focus on the family of $\alpha$-connections $\nabla^\alpha$. Since $\mathbb{M}$ is special and complete, to every $\alpha \in \mathbb{R}$ we assign the foliation $\mathcal{F}^\alpha$ which is auto-parallel in $(M,\nabla^\alpha )$. Leaves of $\mathcal{F}^\alpha$ are orbits of the simply connected bi-invariant affine Lie group 
$$\left\{G_\alpha,\nabla^\alpha\right\},$$ 
the associative algebra of this affine Lie group is
$$J_\alpha = (J_{\nabla^\alpha},\nabla^\alpha)$$ 
Actually we get a locally flat foliation defined by the pair $$\left\{\mathcal{F}^\alpha,\nabla^\alpha\right\}.$$
We can apply the process of statistical reduction as in \cite{Amari-Nagaoka}, \cite{Nguiffo Boyom(6)}.
By the reduction theorem since the the following AAS is auto-parallel in $(M,\nabla^\alpha)$
$$\mathcal{J}_\alpha = \left\{\mathcal{J}_{\nabla^\alpha}, \nabla^\alpha\right\}.$$
Further the dual pair $(M,g,\nabla^\alpha,\nabla^{-\alpha})$ endows the sheaf $\mathcal{J}_\alpha$ with a structure of dual pair
$$\left\{\mathcal{J}_\alpha, g|_{\mathcal{J}^\alpha},\nabla^\alpha,\nabla^{-\alpha}\right\}$$
Now we go to point out a significant consequence. Let $F$ be a leaf of $\mathcal{J}_\alpha$. Then the following quadruple is a dually flat pair
$$\left\{F,g|_F,\nabla^\alpha,\nabla^{-\alpha}\right\}.$$
Therefore we consider the foliation $\left\{\mathcal{F}^\alpha,g,\nabla^\alpha\right\}$. It is an auto-parallel Hessian foliation in $(M,\nabla^\alpha)$ $\clubsuit$
\textbf{In view of those investigations we go to state a structure theorem of the statistical models}
\begin{thm} (A structure Theorem)
Consider a special complete regular statistical model $\mathbb{M}$ whose Fisher metric is denoted by $g$.\\
(1) $\mathbb{M}$ admits a one parameter family of Hessian foliations 
$$\left\{\mathcal{F}^\alpha, g|_{F^\alpha},\nabla^\alpha,\quad \alpha \in \mathbb{R}\right\}$$ which converges to a flat Hessian foliation as $\alpha$ tends to zero \\
(2) Further the family $\mathcal{F}^\alpha$ is defined by the Geodesic Flat Hessian Special Dynamical Systems of one parameter family of simply connected bi-invariant affine Lie groups 
$$ \left\{ G_\alpha,\nabla^{\alpha 0} \right\}\clubsuit$$ 
\end{thm}
\begin{cor} A special complete regular statistical model $\mathbb{M}$ admits a one parameter family of foliations by exponential families. This family of foliations converges to a foliation by flat exponetial families as the parameter tends to zero $\clubsuit$
\end{cor}
\subsection{Special statistical manifolds}
\begin{defn}
 A statistical manifold is a triple $(M,g,\nabla)$ where $(M,g)$ is a Riemannian manifold, $\nabla$ is a torsion free Koszul connection satisfying the following identity
$\nabla_Xg(Y,Z) - \nabla_Yg(X,Z) = 0.$
\end{defn}
 The functor of Amari endows every statistical manifold $(M,g,\nabla)$ with a structure of DUAL PAIR $(M,g,\nabla,\nabla^g)$. A family of $\alpha$-connections of $(M,g,\nabla)$ is defined by
 $$\nabla^\alpha = \frac{1+\alpha}{2}\nabla + \frac{1-\alpha}{2}\nabla^g.$$
 All of the connections $\nabla^\alpha$ are symmetric.
\begin{defn}
 A statistical manifold $(M,g,\nabla)$ is called special-complete if all of its $\alpha $-connections are special-complete $\clubsuit$
\end{defn}
The following statement is a straight corollary of our discussion about the Riemannian geometry
\begin{thm} (a): Every special-complete statistical manifold $(M,g,\nabla)$ admits a $\nabla$- auto-parallel foliation $\mathcal{F}$. The triple $\left\{\mathcal{F}, \nabla, g|_\mathcal{F}\right\}$ is a flat Hessian foliations.\\
(b): The leaves of $\mathcal{F}$ are obits of a GFHSDS 
$$\left\{(G_0,\nabla^0)\times(M,\nabla)\rightarrow (M,\nabla)\right\} \clubsuit$$
\end{thm}
The information geometry  highlights the richness of the differential topology of statistical manifolds, see \cite{Nguiffo Boyom(6)}
\section{THE FUNDAMENTAL EQUATION $FE^*(\nabla)$ AND HOMOGENEOUS KAEHLERIAN GEOMETRY}
In this section we sketch the study of links between the fundamental equation $FE^*(\nabla)$ and the differential topology of homogeneous Kaehler manifolds. Our motivation is an old (ex) conjecture of Gindikin-Pyatecci-Shapiro-Vinberg. In \cite{GPSV} Gindiking the authors state the following fundamental conjecture.\\ 
\textbf{The fundamental conjecture: Every homogeneous Kaehlerian manifold $(M,g,J)$ is the total space of a fiber bundle on a bounded domain. The fibers are product of a compact homogeneous Kaehlerian manifolds and a simply connected homogeneous Kaehlerian manifolds $\clubsuit$}\\
This conjecture has been partially demonstrated by Dorfmeister in \cite{Dorfmeister} and totally demonstrated by Dorfmeister and Nakajima in \cite{Dorfmeister-Nakajima}. Their demonstarion is long. They widely implement the theory of j-algebra which had been introduced by Koszul \cite{Koszul(5)}. See also Nakajima K \cite{Nakajima} In this fundamental conjecture the fibration is not a Kaehlerian fibration, it is uniquely analytic.\\
We aim at using the global analysis of the first fundamental equation $FE^*(\nabla)$ for introducing another fibration of homogeomeous Kaehlerian manifolds over homogeneous spaces. The function $r^b$ allows to obtain finer informations about the topology of the fibers and about the geometry of the base space of the fibration.\\  
Our theorem and the former fundamental conjecture are alternative informations on the differential topology of homogeneous Kaehlerian manifolds.  However our demonstaration is shorter than the demonstartion of Dorfmeister and Nakajima.\\ 
\textbf{Reminder}\\
As announced we go to use the function $r^b$ for introducing another type of fibration of homogeneous Kaehlerian manifolds. \\
\textbf{Reminder}.\\
(a): The $m$-dimensional complex flat torus is the quotient space
$$\mathbb{T}^m = \frac{\mathbb{C}^m}{\mathbb{Z}^\frac{m}{2} + i\mathbb{Z}^{\frac{m}{2}}}$$
(b): An $m$-dimensional complex twisted torus is the quotient
$$ \mathcal{T}^m = \frac{T^m}{K},$$
here $K$ is a finite subgroup of $SU(m)$.\\
(c) A $m$-dimensional complex cylinder over a complex $n$-dimensional twisted torus is the direct product 
$$ \tau = \mathcal{T}^n \times \mathbb{C}^{m-n}.$$ 
(d) A foliation $\mathcal{F}\subset M$ is called simple if the leaves space endowed with the quotient topology is a differentiable maniflods; we say that $M$ is simply foliated.\\
Here is an alternative to the former fundamental conjecture of Gindikin-Pyatecii-Shapiro-Vinberg as in \cite{GPSV}  
\begin{thm} Let $G$ be a finite dimensional real Lie group.\\
(i) Every finite dimensional $G$-homogeneous Kaehlerian manifold $(M,g,J)$ admits a simple foliation $\mathcal{F}$ whose leaves are homogeneous  Kaehlerean submanifolds. Further the leaves are biholomorphic to the same complex homogeneous cylinder over a complex twisted torus.\\
(ii) Further the space of leaves manifold is a homogeneous space of the Lie group $G$ $\clubsuit$
\end{thm}
\textbf{Warning}.\\
The theorem we just stated and the former fundamental conjecture look alike. However the demonstartion we go to give is shorter than the demonstration of the former fundamental conjecture given by Dorfmeister in \cite{Dorfmeister} and the demonstration given by Dorfmeister and Nakajima in \cite{Dorfmeister-Nakajima}.\\
\begin{lem} Let $\gamma^*\in G$ and assume that there exist $x^* \in F$ such that $\gamma^*(^*x) \in F$, then
$$ \gamma^* \in G_F \clubsuit$$
\end{lem}
\textbf{Hint} The reason is that an element of $G$ sends a leaf on another leaf.
Therefore the differentiable dynamc system
 $$G_F \times F \rightarrow F $$
is transitive. Every leaf $F$ supports two transitive differentiable dynamical systems,\\
$\Delta.K$: the Kaehlerian dynamic
$$\left\{G_F \times F \rightarrow F\right\}$$
$\Delta.G$: the Geodesic Flat Hessian Dynamic
$$\left\{(G_\nabla,\nabla^0)\times (F,\nabla,g|_F)\rightarrow (F,\nabla,g|_F)\right\}.$$
We recall that the Geodesic Flat Hessian Dynamic does not preserves the Hessian metric tensor $g$.
\begin{lem} 
The subgroup $G_F \subset G$ is closed $\clubsuit$
\end{lem}
\textbf{Hint}. Take $\gamma$ in the adherence of $G_F$ and choose a sequence 
$$\left\{\gamma_n \in G_F\right\}$$  
which converge to $\gamma$ as $n$ tends to $+\infty$, viz
$$\lim_{\left\{n\rightarrow +\infty\right\}}\gamma_n = \gamma.$$
Assume that $\gamma$ does not belong to $G_F$. Then there exists $x^* \in F$ whose image $\gamma(x^*)$ does not belong to F. Therefore 
$$F\cap F_{\gamma(x^*)} = \emptyset.$$
In $F_{\gamma(x^*)}$ consider an open neighourhood of $\gamma(x^*)$, $U \subset F_{\gamma(x^*)}$. There exists a positive integer $N^*$ such that
$$\left\{\gamma_n(x^*) |\quad N^*\leq n\right\} \subset U$$
The existence of $N^*$ implies
$$F\cap F_{\gamma(x^*)} \neq \emptyset $$
This yields a contradiction. The lemma is proved.\\
Subsequent facts.\\
(3) At one side, by Lemma just proved the following homogenous space is Hausdorff 
$$\mathcal{M} = \frac{G}{G_F},$$
(4) at another side we have the equality
$$ \frac{G}{G_F} = \frac{M}{\mathcal{F}_\nabla}.$$
In final the auto-parallel foliation $\mathcal{F}_\nabla \subset ((M,\nabla)$ is simple
By (2) every leaf $F$ carries the structure of flat Riemannian manifold
$$\left\{F,g|_F, J|_F\right\}$$
Let $q$ be the real dimension of $F$ and let $b_1(F)$ be first real Betti number of $F$. Then up to finite covering the real analytic manifold $F$ is isomophic to the tube
$$\mathcal{T}^{b_1(F)}\times \mathbb{R}^{q - b_1(F)},$$ \cite{Wolf}
On one hand the universal covering of $(F,g|_F,\nabla|_F$ is the flat Euclidean space, viz
$$\left\{\tilde{F}, \tilde{g}|_{\tilde{F}}, \tilde{\nabla}|_{\tilde{F}}\right\} = \left\{\mathbb{R}^q, g°,D° \right\}$$
On another hand we have the topological decomposition
$$ \tilde{F} = \mathbb{R}^{b_1(M)}\times \mathbb{R}^{q - b_1(M)}.$$
Now we take into account the action of the affine holonomy group of $(F,\nabla|_F)$, namely $H(\nabla|_F)$. Let $h(\nabla|_F)$ be the linear component of $H(\nabla|_F)$. we have 
$$ h(\nabla|_F) \subset SU(\frac{q}{2}),$$
$$ \mathcal{T}^{b_1(F)}\times \mathbb{R}^{q - b_1(M)} = \frac{\mathbb{R}^q}{H(\nabla|_F)}.$$
We conclude that $b_1(F)$ is even and
$$ H(\nabla|_F)\subset SU(\frac{b_1(M)}{2})\times \mathbb{C}^{\frac{b_1(M)}{2}}.$$ 
In final all of the leaves of $\mathcal{F}_\nabla $ are biholomorphic to a  $\frac{q}{2}$-dimensional complex cylinder over a $\frac{b_1(M)}{2}$-dimensional complex twisted torus.\\
The part (i) of Theorem is proved.\\
The part (ii) of Theorem is nothing else than the identification (4) of this demonstration. Theorem is demonstrated $\clubsuit$\\
\textbf{Warning} The systems 
$$\left\{(G_q,\nabla^q),\quad I^p_q; \Pi^q_p| \quad p \leq q\right\}$$
are systems of finite dimensional simply connected complex Lie groups. 
\subsection{The fundamental equation $FE^*(\nabla)$ and the problem $EXF(\mathcal{S})$ in homogeneous Kaehlerian manifolds, canonical affine representations}
\textbf{The theorem we just demonstrated in the last susbsection links the equation $FE^*(\nabla)$ with some canonical affine dynamics of isotropy groups}.\\
Let $\mathcal{G}_x$ and $\mathcal{G}_{F_x}$ be the Lie algebras of $G_x$ and $G_{F_x}$ respectively. They are Lie algebras of infinitesimal transformations of 
$$\left\{F_x,g|_{F_x},\nabla|_{F_x}\right\},$$
Further $\mathcal{G}_{F_x}$ is a Lie algebra of infinitesimal automorphisms of 
$$\left\{F_x,g|_{F_x}, J|_{F_x}, \nabla|_{F_x}\right\}$$
Since $(F_x, \nabla|_{F_x})$ is a locally flat manifold we get
$$\mathcal{G}_x \subset \mathcal{G}_{F_x} \subset J_{\nabla|_{F_x}}$$
Up to diffeomorphisms the universal coverings are
$$\left\{\tilde{F}_x,g|_{F_x},\tilde{\nabla}|_{F_x} \right\} = \left\{\mathbb{R}^q,g°, D°\right\},$$
$$\left\{\tilde{F}_x,J\right\} = \mathbb{C}^{\frac{q}{2}};$$
$G_\nabla $ is the simply connected complex Lie group whose Lie algebra is $J_\nabla$.
\begin{thm} The differentiable dynamical systems $\Delta.G$ and $\Delta.K$ yield the affine representations\\
$$\Delta.G:\quad G_\nabla \ni \gamma\rightarrow (\rho(\gamma),\tau(\gamma))\in Gl(q,\mathbb{R})\times\mathbb{R}^q,$$
$$\Delta.K:\quad G_{F_x} \rightarrow SU(\frac{q}{2}) \subset Gl(q,\mathbb{R})\times\mathbb{C}^{\frac{q}{2}}.$$
\end{thm}
Before proceeding we observe that
$$q = r^b(\nabla)$$
Then we get the following link of $r^b$ with the Kaehlerian dynamic.
\begin{thm}
Let $G$ be a finite dimensional real Lie group acting transitively in a Kaehlerean manifold $(M,g,J)$ the Levi-civita connection of which is denoted by $\nabla$. \\
(1) $r^b(\nabla)$ is a positive even integer.\\
(2) Every isotropy subgroup $G_x \subset G$ admits a canonical unitary affine representation 
$$G_x \ni \gamma\rightarrow (\lambda(\gamma),\tau(\gamma))\in 
SU(\frac{r^b(\nabla)}{2})\times \mathbb{C}^\frac{r^b(\nabla)}{2}\clubsuit$$   
\end{thm}
\subsection{The first fundamental equation $FE^*(\nabla)$ and flat foliations in geodesically complete Kaehlerian Geometry}. Assume that $(M,g,J)$ is a geodesically complete Kaehlerian manifold whose Levi-Civita connection is denoted by $\nabla$. We assume the gauge structure $(M,\nabla)$ to be spacial. The associative algebra $J_\nabla,\nabla)$ is the infinitesimal versus of a locally effective action of the simply connected bi-invariant complex affine Lie group $G_\nabla$. The orbits of $G_\nabla$ are flat Kaehlerian submanifolds $(F_x,g, J,\nabla)\subset (M,g,J,\nabla)$. Of course one has 
$$ dim(\mathcal{T}_\nabla M(x)) = r^b(\nabla)(x) \in 2\mathbb{Z}\quad \forall x \in M. $$
Minimal leaves are $r^b(\nabla)$-dimensional. We go to summarize
\begin{thm} Every (special) geodesically complete Kaehlerian manifold admits a holomorphic foliation $\mathcal{F}$ whose leaves are flat Kaehlerian submanifolds. Further every leaf of $\mathcal{F}$ is the orbit of a simply connected bi-invariant complex affine Lie group $\clubsuit$
\end{thm}
\textbf{Warning:\\ 
(1) Let $(M,g,J)$ be a Kaehlerian manifold whose Levi-Civita connection is denoted by $\nabla$. The transformation group of $(M,g,J)$ is included in the transformation group $Aff(M,\nabla)$. Those groups have canonical representation in the Lie group of automorphisms of the bi-invariant complex affine Lie group $(G_\nabla,\nabla^0).$\\
(2) The canonical representaion of the fundamental group of $(M,g,J)$ is a homomorphism in the Lie group of automorphisms of the bi-invariant complex affine Lie group $(G_\nabla,\nabla^0$}  
\subsection{Additional numerical invariants of gauge structures}
We go to deal with a special Kaehlerian manifolds. We consider such a manifold $(M,g,J)$; its Levi-Civita connection is denoted by $\nabla$. Then $(\mathcal{J}_\nabla,\nabla)$ is a complex analytic AAS. Thereby one get the inductive system of finite dimensional bi-invariant complex affine Lie groups $\mathcal{IS}^\nabla$ and the projective system of finite dimensional bi-invariant complex affine Lie group $\mathcal{PS}_\nabla$. We have introduced a canonical homomorphisms of the fundamental group $\pi_1(M)$ is the Lie groups of automorphisms of those inductive-projective systems $\left\{ \mathcal{IS}^\nabla-\mathcal{PS}_\nabla\right\}$.
The AAS $(\mathcal{J}_\nabla,\nabla)$ is a geometric invariant of $(M,g,J)$. Therefore the pair
$\left\{\mathcal{IS}^\nabla-\mathcal{PS}_\nabla\right\}$ is a geometri invariant of $(M,g,J)$. More generally let $(M,\nabla)$ be a symmetric special gauge structure. we consider the  $\left\{\mathcal{IS}^\nabla-\mathcal{PS}_\nabla \right\}$. We go to define the $p,q){^th}$ Betti number  $b_{p,q}$ of a manifold $M$
\begin{defn}
$$(1):\quad b_{pq}(\nabla) = b_p(G_{\nabla^q}),$$
here $b_p(G_{\nabla^q})$ is the $p^{th}$ real Betti number of the Lie group $G_{\nabla^q}$, viz 
$$ b_p(G_{\nabla^q} = dim(H_{CE}(J_{\nabla^q},[-,-]_{\nabla^q}),$$ \\
$$(2):\quad b_{p,q} = \min_{\left\{\nabla\in \mathcal{LC}(M)\right\}}\left\{b_{p,q}(\nabla)\right\}\clubsuit$$
\end{defn}
\textbf{Warning. The family $\left\{b_{p,q}(\nabla)\right\}$ is a geometric invariant of $(M,\nabla)$. The family $\left\{ b_{p,q}\right\}$ is a global geometric invariant of the manilfold $M$}.
\section{LINKS OF THE FIRST FUNDAMENTAL EQUATION $FE^*(\nabla)$ WITH Cartan-Lie GROUPS AND WITH ABSTRACT Lie GROUPS}
There are two canonical Differentiable Dynamical Systems of an abstract Lie group in itself.\\
(1) The first dynamic is the action by left translations.\\
(2) The second dynamic is the action by inner automorphisms.\\
In the preceding sections we have been interested in the impacts of the function $r^b$ on the bi-invariant locally flat structures in Cartan-Lie groups. In this section we go to study the problems $EX(\mathcal{S})$, $EXF(\mathcal{S})$ and $DL(\mathcal{S})$ for left invariant structures and for bi-invariant structures in finite dimensional Lie groups.\\
We restrict ourself to structures which have significant impacts on both the theorectical information geometry \cite{Amari}, \cite{Amari-Nagaoka}, \cite{McCullagh}, \cite{Nguiffo Boyom(6)}] and on the applied information geometry. 
\cite{Arnaudon-Barbaresco}, \cite{Barbaresco}, \cite{Pennec}.\\
Manifolds and Lie groups we deal with are finite dimensional. We recall some open problems.\\
$EX(\mathcal{S}: lf(G))$ is the existence of left invariant locally flat structures in a Lie group $G$.\\
$EX(\mathcal{S}: br(G))$ is the existence of bi-invariant Riemannian metrics in a Lie group $G$.\\
Two other challenges we go to face are\\
$EX(\mathcal{S}: Symp(M))$ is the existence of symplectic structures in a differentiable manifold $M$.\\
$EX(\mathcal{S}: symp(G))$ is the existence of left invariant symplectic structures in a Lie group $G$.\\
A significant new insight of the fundamental equations $FE^*(\nabla\nabla^*)$, $FE^*(\nabla)$ and $FE^{**}(\nabla)$ is their links with those open problems. 
\subsection{The fundamental equations $FE(\nabla\nabla^*)$ and  $FE^*(\nabla)$ in abstract Lie groups}
\textbf{The gauge geometry}\\
Let $G$ be a $m$-dimensional Lie group whose Lie algebra is denoted by $\mathcal{G}$. Though the set $\mathcal{LC}(G)$ is not a vector space it contains some interesting finite dimensional vector spaces,\\
(1): $lc(\mathcal{G}$ is the vector space of left invariant Koszul connections in $G$; \\
that is but the $m^3$-dimensional vector space of bilinear maps of $\mathcal{G}$ in itself;\\
(2): $bc(\mathcal{G})$ is the vector space of bi-invariant Koszul connections in $G$; \\
(3): $lf(\mathcal{G})$ is the set of left invariant locally flat Koszul connections in $G$,
(4): $slc(\mathcal{G})$ is the convex set of symmetric Koszul connections in $G$
We have the inclusions
$$lf(\mathcal{G})\subset slc(\mathcal{G})\subset lc(\mathcal{G})\subset \mathcal{LF}(G).$$
\textbf{The Riemannian geometry}\\
We go to focus on Riemannian geometry is finite dimensional Lie groups.\\
$(1^*)$: $lr(\mathcal{G})$ is the real cone of left invariant Riemannian metrics in $G$;\\
$(2^*)$: $br(\mathcal{G})$ is the cone of bi-invariant Riemannian metrics in $G$;\\ 
$(3^*)$: $br_+(\mathcal{G})$ is the convex cone of bi-invariant positive Riemannian metrics in $G$.\\
\textbf{The symplectic geometry}\\
(1**): $\mathcal{SYP}(M)$ is the set of symplectic structures in a manifold $M$,\\
(2**): $symp(\mathcal{G})$ is the set of left invariant symplectic structures in a Lie group $G$.\\
\textbf{Reminders} \\
The vector space lc($\mathcal{G})$ contains three canonical connections named Cartan connections. They are described in terms of the torsion tensor:\\
$ the -connection: \quad \nabla^-_XY = 0 \quad\forall X,Y \in \mathcal{G},$\\
$ the 0-connection: \quad \nabla^0_XY = \frac{1}{2}[X,Y] \quad\forall X,Y \in \mathcal{G},$\\
$ the +connection: \quad \nabla^+_XY = [X,Y]\quad\forall X,Y \in \mathcal{G}.$
Every bi-invariant Riemannian metric admits those Cartan connections as metric connections. 
Consider $\nabla \in lc(\mathcal{G})$. The vector space of symmetric bilinear forms in $\mathcal{G}$ which are $\nabla$-parallel is denoted by $\mathcal{S}^\nabla_2(\mathcal{G})$;\\
the vector space of skew symmetric bilinear forms in $\mathcal{G}$ which are $\nabla$-parallel is denoted by $\Omega^\nabla_2(\mathcal{G})$.\\
\textbf{We go to perform the functor of Amari} \cite{Nguiffo Boyom(6)}\\
Every left invariant Riemannian metric $g$ defines a linear map
$$lc(\mathcal{G})\ni \nabla\rightarrow \nabla^g \in lc(\mathcal{G})$$
where $\nabla^g$ is defined by 
$$ g(\nabla^g_XY,Z)+g(Y,\nabla_XZ) = 0 \quad \forall X,Y,Z \in \mathcal{G}.$$
\textbf{We go to deal with the differential equation $FE(\nabla\nabla^*)$ and with the sheaf $\mathcal{J}_{\nabla\nabla^*}$}\\
We conside the vector subspace 
$$\mathcal{M}(\nabla,\nabla^g) = \left\{\phi \in Hom_{\mathbb{R}}(\mathcal{G},\mathcal{G})|\quad 
\nabla^g_X\circ \phi - \phi\circ \nabla_X = 0 \quad \forall X \in \mathcal{G}\right\}\clubsuit $$
The formula just posed is is the fundamental equation $FE(\nabla\nabla^*)$ versus Lie theory.\\
Thereby one has 
$$ \mathcal{M}(\nabla,\nabla^g) = J_{\nabla\nabla^g}\cap\mathit{g}(|G|),$$
we recall that $\mathit{g}(G)$ is the Lie algebra of infinitesimal gauge transformations of the manifold $|G|$
To every $\phi \in\mathit{g}(|G|)$ are assigned the bilinear forms 
$B(g,\phi)$ and $\omega(g,\phi)$ which are defined by \\
$$B(g,\phi)(X,Y) = \frac{1}{2}[g(\phi(X),Y) + g(X,\phi(Y))],$$
$$\omega(g,\phi)(X,Y) = \frac{1}{2}[g(\phi(X),Y) - g(X,\phi(Y))].$$
\begin{defn} To every pair 
$$\phi \in \mathcal{M}(\nabla,\nabla^g)$$
is assigned a unique pair 
$$(\Phi, \Phi^*)\in g(|G|)^2$$
which is defined by
$$g(\Phi(X),Y) = B(g,\phi)(X,Y),$$
$$g(\Phi^*(X),Y) = \omega(g,\phi)(X,Y).$$
\end{defn}
We go to state the Lie theory versus of a already state proposition.
\begin{prop} For every triple 
$$(g,\nabla,\phi)\in lr(\mathcal{G})\times lc(\mathcal{G})\times\mathcal{M}(\nabla,\nabla^g)$$
the following assertions are equivalent\\  
 $$(1):\quad (B(g,\phi),\omega(g,\phi))\in \mathcal{S}^\nabla_2(\mathcal{G})\times \Omega^\nabla_2(\mathcal{G}),$$
 $$(2):\quad (\Phi,\Phi^*)\in \mathcal{M}(\nabla,\nabla^g).$$
Furthermore the following maps are sujective\\
 $$(3):\quad \mathcal{M}(\nabla,\nabla^g) \ni \phi \rightarrow B(g,\phi) \in \mathcal{S}^\nabla_2(\mathcal{G}),$$
 $$(4):\quad \mathcal{M}(\nabla,\nabla^g) \ni \phi\rightarrow \omega(g,\phi) \in \Omega^\nabla_2(\mathcal{G})\clubsuit$$
\end{prop}
\subsubsection{The left invariant locally flat geometry}
We go to deal with $lc(\mathcal{G})$. When there is no risk of confusion we will make simpler by setting
$$X.Y = \nabla_XY. $$
We consider $\nabla \in lc(\mathcal{G})$. In the manifold $|G|$ we consider the sheaf $\mathcal{J}_\nabla$ whose sections are solutions to the fundamental equation $FE*(\nabla)$. According to the notation used in the preceding sections the vector space of sections of $\mathcal{J}_\nabla$ is denoted by $J_\nabla$. We consider the vector space $\mathcal{G}_\nabla$ which is defined by 
$$\mathcal{G}_\nabla = J_\nabla \cap \mathcal{G}.$$
The pair $(\mathcal{G}_\nabla,\nabla)$ is an associative sub-algebra of the (non associative) algebra $(\mathcal{G},\nabla)$. Henceforth we focus on the sets $lc(\mathcal{G})$, $slc(\mathcal{G})$ and $lf(\mathcal{G})$. We recall the inclusions
$$lf(\mathcal{G})\subset slc(\mathcal{G})\subset lc(\mathcal{G}) \subset \mathcal{LC}(|G|),$$
$$br_+(\mathcal{G})\subset br(\mathcal{G})\subset lr(\mathcal{G})\subset \mathcal{R}ie(|G|).$$
\textbf{We recall the notation already used}\\
$1_G$ is the neutral element of the group $G$,\\
$$\mathcal{T}_\nabla G(1_G) = \frac{J_\nabla}{I_{1_G}}.$$
For $\nabla \in lc(\mathcal{G})$ the vector subspace $J_\nabla \subset \mathcal{X}(G)$ is left invariant in the Lie group $G$. Consequently 
$$rank(\mathcal{T}_\nabla G) = constant.$$ 
\begin{defn}
$$r^b(\nabla) = dim(\mathcal{T}_\nabla G(1_G)) $$
$$ r^b(G) = \min_{\left\{\nabla\in slc(\mathcal{G})\right\}}\left\{dim(G)-r^b(\nabla)\right\}\clubsuit$$
\end{defn}
Here is Theorem EXTH1 versus Lie theory
\begin{thm} In a finite dimensional Lie group $G$ the following assertions are equivalent.\\
$(1): r^b(G) = 0.$\\
$(2):$ The Lie group $G$ admits a left invariant locally flat structures $\clubsuit$
\end{thm}
The Demonstration of the theorem just stated is similar to the demonstration of Theorem EXTH1\\
\textbf{Remark}. We have the inequality
$$ r^b(|G|) \leq r^b(G) \leq dim(G).$$
\textbf{Examples.}\\
Examples. Let $G$ be the universal covering of $SL(2,\mathbb{R})$,
$$r^b(G) > 0,$$
$$r^b(|G|) = 0.$$
In \cite{Benoit} Yves Benoit has used computing methods for constructing an $11$-dimensional nilpotent Lie groups $G$ with 
$$ r^b(G) > 0.$$
If $G$ is a finite dimensional perfect Lie group then
$$r^b(G) > 0.$$
\textbf{The differential topology, continued}.\\
Let $G$ be a $q$-dimensional Lie group with
$$ 0 < r^b(G) < q$$
Therefore $G$ admits a left invariant torsion free connection $\nabla$ with
$$ r^b(\nabla) = q - r^b(G)$$
The associative algebra $(J_\nabla,\nabla)$ is finite dimensional and is invariant by the left translations in $G$. Therefore the distribution
$$\mathcal{D} = C^\infty(G)J_\nabla = span(J_\nabla). $$
defines a $r^b(G)$-co-dimensional left invariant affine foliation.\\
In a Lie group the invariant $r^b$ has the following meaning.
\begin{thm} Every finite dimensional Lie group $G$ with
$$ 0 < r^b(G) < dim(G).$$
admits a left invariant $r^b(G)$-codimensional locally flat foliation $\mathcal{F}$ whose leaves are orbits of a bi-invariant affine Cartan-Lie group. Furthermore $r^b(G)$ is the optimal codimension of such a type of foliation $\clubsuit$
\end{thm}
\subsubsection{The bi-invariant Riemannian geeometry in finite dimensional Lie groups}
\textbf{The problem is $EX(\mathcal{S}: br(\mathcal{G}))$. Our approach involves the +connection of Elie Cartan $\nabla^+$ and the functor of Amari.} Here\\
$$\nabla^+_XY = [X,Y]\quad \forall X, Y \in \mathcal{G}.$$ 
We consider the functor of Amari 
$$lr(\mathcal{G}) \ni g \rightarrow \nabla^{+g}\in lc(\mathcal{G}).$$
Here $\nabla^{+g}$ is defined by
$$ g(\nabla^{+g}_XY,Z) + g(Y,[X,Z]) = 0\quad \forall x, Y, Z \in \mathcal{G}.$$
To every $g \in lr(\mathcal{G})$ we assign the vector space 
$$\mathcal{M}(\nabla^+,\nabla^{+g}) \subset J_{\nabla\nabla^g}.$$
An element $\phi \in \mathcal{M}(\nabla^+,\nabla^{+g})$ is defined by \\
 $$\nabla^{+g}_X\phi(Y) - \phi([X,Y] = 0\quad \forall X, Y \in \mathcal{G}.$$
To $\phi \in \mathcal{M}(\nabla^+,\nabla^{+g}$ we assign the unique $\Phi \in \mathcal{M}(\nabla^+,\nabla^{+g})$  which is defined by
$$g(\Phi(X),Y) = \frac{1}{2}[g(\phi(X),Y) + g(X,\phi(Y))]\quad \forall X, Y \in \mathcal{G}.$$
We set
$$\mathcal{M}+(\nabla^+,\nabla^{+g}) = \left\{\phi\in \mathcal{M}(\nabla^+,\nabla^{+g})\quad | g(\Phi(X),X) > 0 \quad \forall X\neq 0. \right\}$$
\begin{defn} To every $g\in lr(\mathcal{G}$ we assign the following integers
$(1):\quad s^b(G,g) = \min_{\left\{\phi \in \mathcal{M}(\nabla^+,\nabla^{+g}\right\}}\left\{dim(G) - rank(\Phi)\right\},$\\
$(2):\quad s^b_+(G,g) = \min_{\left\{\phi \in \mathcal{M}+(\nabla^+,\nabla^{+g})\right\}}\left\{dim(G) - rank(\Phi)\right\},$\\
$(3):\quad s^b(G) = \min_{\left\{g \in lr(\mathcal{G})\right\}}\left\{s^b(G,g)\right\},$\\
$(4):\quad s^b_+(G) =\min_{\left\{g \in lr(\mathcal{G})\right\}}\left\{s^b_+(G,g)\right\} \clubsuit $
\end{defn}
Those integers are global geometric invariants of the Lie group $G$.
\subsubsection{The problem $EX(\mathcal{S}: br(\mathcal{G}))$} 
We go to relate the problem $EX(\mathcal{S}: br(\mathcal{G}))$ and the numerical invariants we just defined in the preceding subsubsection.
\begin{thm} In a finite dimensional Lie group $G$ the following assertions are equivalent\\
(1): $ s^b(G) = 0,$\\
(2): the Lie group $G$ admits bi-invariant Riemannian metrics $\clubsuit$
\end{thm}
\begin{cor} In a finite dimensional Lie group $G$ the following assertions are equivalent\\
$(1): s^b_+(G) = 0$,\\
$(2):$ the Lie group $G$ admits bi-invariant positive Riemanian metrics $\clubsuit$
\end{cor}
The demonstrations the theorems above are similar to the demonstrations of Theorems $EXTH$\\
\textbf{Topology continued.} From our investigations emerge two approachs to study the question whether a Koszul connection $\nabla$ is a metric connection.\\
(1) The first approach is based on the computation of the total KV cohomology as in \cite{Nguiffo Boyom-Byande-Ngakeu-Wolak}. This method works in the category $\mathcal{LF}$ of locally flat connections. The key formula is the decomposition of the total KV cohomology space
$$H^2_\tau(M,\nabla) = H^2_{dR}(M) \oplus \mathcal{S}^\nabla_2(M).$$
(2) The second approach is based on a arsenal formed of the functor of Amari and the solutions to the fundamental equation $FE^*(\nabla\nabla^*)$. We just perfomed this method for studying the problem $EX(\mathcal{S}: br(G))$.
\textbf{Actually there exists another approach which is based on the dynamic of discrete groups. This approach is the object of forthgoing work}
\subsection{The problem $DL(\mathcal{S}: br(\mathcal{G}))$}
$s^b(G)$ measures how far from admitting bi-invariant Riemannian metrics is a Lie group $G$.\\
$s^b_+(G)$ measures how far from admitting bi-invariant positive Riemannian metrics is a Lie group $G$.
\section{THE FUNDAMENTAL EQUATION $FE(\nabla\nabla^*)$ AND THE SYMPLECTIC GEOMETRY: THE SYMPLECTIC GAPS}
\textbf{For convenience we go to recall the concerns}.\\
(1): Our purpose is to develop tools which involve both $J_{\nabla\nabla^*}$ and the functor of Amari.\\
(2): Our aim is to implement those tools in the symplectic geometry.\\
(3) We are in search of invariants measuring how far from admitting symplectic structures is a finite dimensional manifold.\\
(4) We are in search of invariants measuring how far from admitting left invariant symplectic structures is a finite dimensional Lie group.\\
Those invariants are called symplectic gaps.\\
We recall the framework of our our purposes.\\
 $\mathcal{LC}(M)$ is the affine space of Koszul connections in $M$.\\
 $\mathcal{SLC}(M)$ is the convex set of torsion free Koszul connections in $M$.\\
 $\mathcal{R}ie(M)$ is the cone of Riemannian metrics in $M$
 $slc(\mathcal{G})$ is the convex set of left invariant symmetric connections in a Lie group $G$.\\
 $lr(\mathcal{G})$ is the cone of left invariant Riemannian metrics in $G$.\\
 We recall the functor of Amari in a differentiable manifold is defined as it follows
 $$\mathcal{SLC}(M)\times \mathcal{R}ie(M) \ni (\nabla,g)\rightarrow \nabla^g \in \mathcal{LC}(M)$$
 Here the Koszul connection $\nabla^g$ is defined as it follows
 $$g(\nabla^g_XY,Z) = Xg(Y,Z) - g(Y,\nabla_XZ).$$
 Mutatis mutandis one gets the left invariant functor of Amari
 $$slc(\mathcal{G}) \times lr(\mathcal{G}) \in (\nabla,g)\rightarrow \nabla^g \in lc(\mathcal{G}).$$
 Here $\nabla^g$ is defined by 
 $$g(\nabla^b_XY,Z) + g(Y,\nabla_XZ) = 0 \quad \forall X, Y, Z \in \mathcal{G}.$$
It is easy to see that $\nabla^g$ is a left invariant Koszul connection in $G$;\\
In a Lie group $G$ we restrict the attention to the following vector spaces
$$Hom_\mathbb{R}(\mathcal{G},\mathcal{G})\subset \mathit{g}(|G|),$$ 
$$\mathcal{M}(\nabla,\nabla^g) = \left\{\phi \in Hom_\mathbb{R}(\mathcal{G},\mathcal{G}) \quad | 
\nabla^g_X\circ\phi - \phi\circ \nabla = 0\quad \forall X \in \mathcal{G}\right\}.$$
In a finite dimensional Lie group $G$ we go to deal with triples\\
 $$(g,\nabla,\phi) \in lr(\mathcal{G})\times lc(\mathcal{G})\times \mathcal{M}(\nabla,\nabla^g).$$
Up to notation change we recall a former construction. To a triple $(g,\nabla,\phi)$ we assign the bilinear forms $q(g,\phi)$ and $\omega(g,\phi)$ which are defined by
$$q(g,\phi) (X,Y) = \frac{1}{2}(g(\phi(X),Y) + g(X,\phi(Y)),$$
$$\omega (g,\phi)(X,Y) = \frac{1}{2}(g(\phi(X),Y) - g(X,\phi(Y)).$$
We implement a former proposition. For every $\phi \in \mathcal{M}(\nabla,\nabla^g)$ there exists a unique pair $(\Phi, \Phi^*) \in \mathcal{M}(\nabla,\nabla^g)$ which is defined by
$$g(\Phi(X),Y) = q(g,\phi)(X,Y),$$
$$g(\Phi^*(X),Y) = \omega(g,\phi)(X,Y).$$
In preceding subsections we have already proved that
$$(q(g,\phi),\omega(g,\phi) \in \mathcal{S}^\nabla_2(G)\times \Omega^\nabla_2(G)$$
\subsection{The symplectic gap of differentiable manifolds}
\textbf{Reminder}
We go to use the solutions to $FE(\nabla\nabla^*)$, namely $J_{\nabla\nabla^*}$ and the functor of Amari. To start we consider the following two maps
$$\mathcal{SLC}(M)\times \mathcal{R}ie(M)\ni (\nabla,g)\rightarrow \nabla^g\in \mathcal{LC}(M),$$
$$J_{\nabla\nabla^g}\ni phi \rightarrow \Phi,Phi^* \in J_{\nabla\nabla^g}.$$
For convenience we recall that for every $\phi\in J_{\nabla\nabla^g}$ the unique pair $(\Phi,\Phi^*)$ is defined by
$$g(\Phi(X),Y) = \frac{1}{2}[g(\phi(X),Y) + g(X,\phi(Y))] $$
$$g(\Phi^*(X),Y) = \frac{1}{2}[g(\phi(X),Y) - g(X,\phi(Y))].$$
We also recall that 
$$(q(g,\phi),\omega(g,\phi) \in \mathcal{S}^\nabla_2(M)\times \Omega^\nabla_2(M)\quad \forall \phi \in J_{\nabla\nabla^g}.$$
Consequently both $\Phi$ and $\Phi^*$ have constant rank.
We go to use $\Phi^*$ for introducing new numerical invriants. 
\begin{defn} The fundamental equation $FE(\nabla\nabla^*)$ for introducing the following functions\\
A. In a finite dimensional differentiable manifold $M$, for every pair $(\nabla,g)\in \mathcal{SLC}(M)\times\mathcal{R}ie(M)$ we put\\
$(1):\quad s^{*b}(\nabla,g) = \min_{[\phi\in J_{\nabla\nabla^g]}}[dim(M) - rank(\Phi*)],$\\
$(2):\quad s^{*b}(M) = \min_{\left\{(\nabla,g) \in \mathcal{SLC}(M)\times \mathcal{R}ie(M)\right\}} \left\{s^{*b}(\nabla,g)\right\}$.\\
B. In a finite dimensional Lie group $G$, for every pair $(\nabla,g)\in slc(\mathcal{G})\times lr(\mathcal{G})$ we put\\
$(3):\quad s^{*b}(\nabla,g) = \min_{[\phi \in \mathcal{M}(\nabla,\nabla^g)]}[dim(G) - rank(\Phi^*)],$\\
$(4):\quad s^{*b}(G) = \min_{\left\{(\nabla,g)\in \sl(\mathcal{G})\times lr(\mathcal{G})\right\}}\left\{s^{*b}(\nabla,g)\right\}.$\\
C. In the union $M\cup G$; for $\nabla \in \mathcal{SLC}(M)$ we put\\
$(5):\quad s^{*b}(M,\nabla) = \min_{[g\in \mathcal{R}ie(M)]} \left\{dim(M) - s^{*b}(\nabla,g)\right\},$\\
For $\nabla \in slc(\mathcal{G})$ we put
$(6):\quad s^{*b}(G,\nabla) = \min_{[g\ni lr(\mathcal{G})]} \left\{dim(G) - s^{*b}(\nabla,g)\right\} \clubsuit$
\end{defn}

\begin{defn} In the defintion above\\
(6) is called the symplectic gap of $\nabla \in slc(\mathcal{G})$,\\
(5) is called the symplectic gap of $\nabla \in \mathcal{SLC}(M)$,\\
(4) is called the symplectic gap of the Lie group $G$,\\
(2) is called the symplectic gap of the manifold $M$ $\clubsuit$
\end{defn}
\subsection{PROBLEM $EX(\mathcal{S}: Symp)$}
We go to solve the problem $EX(\mathcal{S})$ in the symplectic Geometry.
\begin{thm} In a finite dimensional manifold $M$ the following assertion are equivalent\\
$(1): s^{*b(M)} = 0,$\\
$(2)$ the manifold $M$ admits symplectic structures $\clubsuit$
\end{thm}
\textbf{Demonstration}\\
We go to prove that (1) implies (2).\\
By (1) there exists a pair 
$$(\nabla,g)\in \mathcal{SLC}(M)\times \mathcal{R}ie(M)$$
such that
$$s^{*b}(\nabla,g) = 0.$$
We use the functor of Amari for defining $\nabla^g \in \mathcal{LC}(M)$ by
$$ g(\nabla^g_XY,Z) = Xg(Y,Z) - g(Y,\nabla_XZ).$$
There is $\psi J_{\nabla\nabla^g}$ such that 
$$ rank(\Psi^*) = dim(M)$$
We recall that $\Psi^*$ is defined by
$$g(\Psi^*(X),Y) = \frac{1}{2}[g(\psi(X),Y) - g(X,\psi(Y))].$$
Futher the diffrential 2-form $g(\Psi^*(X),Y)$ is $\nabla$-parallel, viz
$$ Xg(\Psi^*(Y),Z) = g(\Psi^*(\nabla_XY),Z) + g(\Psi^*(X),\nabla_XZ).$$
This last property yields the closedness of the differential form
$$\omega (X,Y) = g(\Psi^*(X),Y).$$
In final $(M,\omega)$ is a symplectic manifold.\\
We go to prove that (2) implies (1).\\
By (2) we consider a symplectic form in $M$, namely  $\omega$. By \cite{BCGRS} we choose a special symplectic connection $\nabla \in \mathcal{SLC}(M)$. We consider a Riemannian structure $(M,g)$. Therefore there exists $\psi \in \mathcal{M}(\nabla,\nabla^g)$ such that
$$\omega (X,Y) = g(\Psi^*(X),Y)\quad \forall X, Y \in \mathcal{X}(M).$$
This yields
$$s^{*b}(\nabla,g) = 0.$$
In final 
$$ s^{*b}(M) = 0 \clubsuit$$
\textbf{Left invariant symplectic geometry}.\\
Now we go to adderss the left invariant symplectic geometry in finite dimensional Lie groups.
\begin{thm} In a finite dimensional Lie group $G$ the following assertions are equivalent\\
(1): $s^{*b}(G) = 0,$\\ 
(2): the Lie group $G$ admits left invariant symplectic structures $\clubsuit$
\end{thm}
\textbf{Demonstration}.\\
 We go to prove that (1) implies (2).\\
 By (1) there exists a pair $(\nabla,g) \in slc(\mathcal{G})\times lr(\mathcal{G})$ such that
 $$ s^{*b}(\nabla,g) = O.$$
 Therefore there exits $\phi \in \mathcal{M}(\nabla,\nabla^g)$ such that
 $$ rank(\Phi^*) = dim(G);$$
we recall that $\Phi^* \in \mathcal{M}(\nabla,\nabla^g)$ is defined by
$$ g(\Phi^*(X),Y) = \frac{1}{2}[g(\phi(X),Y) - g(X,\phi(Y))].$$
At one side we have the regular skew symmetric bilinear form 
$$ \mathcal{G}\times \mathcal{G} \ni (X,Y)\rightarrow g(\Phi^*(X),Y)\in \mathbb{R}.$$
At another this form is $\nabla$-parallel, viz
$$ g(\Phi^*(\nabla_XY),Z) + g(\Phi^*(Y),\nabla_XZ) = 0\quad\forall X, Y, Z \in \mathcal{G}.$$
Since $\nabla \in slc(\mathcal{G})$ the form $g(\Phi^*(X),Y)$ is a left invariant symplectic form in $G$.\\
We go to prove that (2) implies (1).\\
By (2) we assume that $G$ admits a left invariant symplectic form $\omega$. \\
We recall that torsion free symplectic connections are called special symplectic connections. For construction of special symplectic connections readers are refered to \cite{BCGRS}. That construction yields the existence of left invariant special symplectic connections in $(G,\omega)$. To proceed we choose a left invariant special symplectic connection $\nabla$ and a Riemannian metric tensor 
$g \in lr(\mathcal{G})$.
The pair $(\nabla,g)$ yields the left invariant connection $\nabla^g$ defined by
$$g(\nabla^g_XY,Z) + g(Y,\nabla_XZ) = 0\quad \forall X, Y, Z \in \mathcal{G}.$$
There exists $\phi \in \mathcal{M}(\nabla,\nabla^g)$ such that
$$\omega (X,Y) = g(\Phi^*(X),Y)\quad \forall X, Y \in \mathcal{G}.$$
we have 
$$s^{*b}(\nabla,g) = dim(G)- rank(\Phi^*) = 0.$$
Consequenty we have 
$$s^{*b}(G) = 0.$$
This ends the demonstration $\clubsuit$ 
\subsection{The problem $DL(\mathcal{S})$ in the symplectic geometry}
\textbf{In a differentiable manifold $M$}\\
$(1): s^{*b}(M)$ measures how far from admitting symplectic structures in a finite dimensional manifold $M$,\\
$(2): s^{*b}(M,\nabla)$ measures how far from being a special symplectic connection is $\nabla \in \mathcal{SLC}(M)$.\\
\textbf{In a Lie group $G$}\\
$(1): s^{*b}(G)$ measures how far from admitting left invariant symplectic structures is a finite dimensional Lie group $G$.\\
$(2): s^{*b}(G,\nabla)$ measures how far from being a special connection of left invariant structures in $\nabla \in slc(\mathcal{G})$.\\
To end this section which is devoted to the symplectic geometry it remind to emphasize other roles played by both $s^b(M)$ and $s^{*b}(M)$.
\begin{defn} We remind that a Riemannian foliation in $M$ is a pair $(\mathcal{F},g)$ formed of a foliation $\mathcal{F}$ and a symmetric bi-linear form $g$ subject to the following requirements.
(1) $ Ker(g) = \mathcal{F},$\\
(2) For every section of $Ker(g)$, namely $X$ one has $L_Xg = 0$.\\
Here $L_X$ is the Lie derivative in the direction $X\clubsuit$
\end{defn}
By replacing a symmetric bi-linear form $g$ by a closed differential 2-form $\omega$ we define the notion of symplectic folliation 
\cite{Nguiffo Boyom(6)}
To be symplectically foliated differs from supporting sympectic foliations. The first is the symplectic geometry along the leaves while the second is the transverse symplectic geometry.\\
\begin{defn} A manifold $M$ is called specially symplectically foliated if it supports a triple $(\mathcal{F},\omega,\nabla)$ where\\
(1): $\nabla$ is a torsion free Koszul connection.\\
(2): $M$ is symplectically foliated by $(\mathcal{F},\omega)$.\\
(3): The pair $(\mathcal{F},\omega)$ is $\nabla$-auto-parallel $\clubsuit$
\end{defn}
\textbf{Warning: The auto-parallelism means the following}\\
$(1)$ If $X$ and $Y$ are tangent to $\mathcal{F}$ then $\nabla_XY$ is tangent to $\mathcal{F}$\\
$(2) \nabla\omega = 0.$\\
\section{THE FUNDAMENTAL EQUATION $FE(\nabla\nabla^*)$ AND THE PROBLEM $EXF(\mathcal{S})$ CONTINUED}
We go to focus on the category of statiscal manifolds. In \cite{Nguiffo Boyom(6)} we have sketched the study of the differential topology of statistical manifolds. For instance we have related some structures of statistical manifold and Riemannian submersions on symplectic manifolds. We go to sketch how to implement $FE(\nabla\nabla^*)$ for constructing webs with rich geometry structures in statistical manifolds.
\subsubsection{Webs in statistical manifolds} 
\textbf{Reminder: A statistical manifold $(M,g,\nabla,\nabla^*)$ is called positive if $g$ is positive definite}. In a special statistical manifold $(M,g,\nabla,\nabla^*)$ we consider the fundamental differential equations $FE(\nabla\nabla^*)$ and $FE^*(\nabla^*\nabla)$. We recall that solutions to those equations are infinitesimal gauge transformations of $TM$ \cite{Petrie-Handal}. \\
We recall the map
$$ J_{\nabla\nabla^*} \ni \phi\rightarrow (\Phi,\Phi^*)\in [J_{\nabla\nabla^*}]^2$$
which defined by
$(a):\quad g(\Phi(X),Y) = \frac{1}{2}[g(\phi(X),Y) + g(X,\phi(Y))] = q(g,\phi)(X,Y),$\\
$(a^*):\quad g(\Phi^*(X),Y) = \frac{1}{2}[g(\phi(X),Y) - g(X,\phi(Y))] = \omega(g,\phi)(X,Y).$\\
The maps $\Phi$ and $\Phi^*$ have constant rank. Furthermore\\
$(i)$: the distribution $Ker(\Phi)$ is a Riemannian foliation in the classical sense \cite{Molino}, \cite{Reinhardt}\\
$(ii)$: the distribution $Ker(\Phi^*)$ is a symplectic foliation in the sense of  \cite{Nguiffo Boyom(6)}\\
We assume that the staistical manifold $(M,g,\nabla,\nabla^*)$ is positive.\\
$(iii)$: The pairs $\left\{Ker(\Phi), im(\Phi) \right\}$ and $\left\{Ker(\Phi^*,im(\Phi^*) \right\}$ are $g$-orthogonal 2-webs in the Riemannian manifold $(M,g)$.\\
The family of $\alpha$-connections of $(M,g,\nabla,\nabla^*)$ is the family of Koszul connections
$$\left\{ \nabla^\alpha = \frac{1+\alpha}{2}\nabla + \frac{1-\alpha}{2}\nabla^*\right\}$$
Here $\alpha \in \mathbb{R}$. Further every dual pair $(M,g,\nabla^\alpha,\nabla^{-\alpha})$ is a staistical manifold.
\begin{defn}
Two infinitesimal gauge transformation $\phi, \psi \in J_{\nabla\nabla^*}$ are called topologically equivalent if $$\left\{ Ker(\Phi), im(\Phi)\right\} = \left\{Ker(\Psi), im(\Psi) \right\}$$
We denote this relation by $\mathcal{TE}(\phi,\psi)$
\end{defn}
This relation means that $\phi$ and $\psi$ yield the same orthogonal 2-web. We emphasize this by a proposition
\begin{prop} 
In a statistical manifold $(M,g,\nabla,\nabla^*)$ the quotient set
$\frac{J_\nabla\nabla^*}{\mathcal{TE}}$
is canonically isomorphic to a set of orthogonal 2-webs in the Riemannian manifold $(M,g)$
\end{prop}
Let $\phi \in J_{\nabla\nabla^*}$
The following statement is a straight consequence of our investigations in the preceding sections.
\begin{thm} Let $(M,g,\nabla,\nabla*)$ be a positive statistical manifold. We assume that the following inequalities hold
$$0 < s^{*b}(\nabla,g) < dim(M),$$
$$0 < s^{*b}(\nabla^*,g) < dim(M).$$
Then the manifold $M$ is specially symplectically foliated by the two triples $(\mathcal{F},\omega),\nabla)$ and $(\mathcal{F}^*,\omega^*,\nabla^*)$. Further $s^{*b}(\nabla,g)$ and $s^{*b}(\nabla^*,g)$ are optimal codimensions for special symplectic foliations in $(M,\nabla)$ and in $(M,\nabla^*)$ respectively $\clubsuit$
\end{thm}
\textbf{The idea of demonstration}. Really we go to repeat arguments already used. The key tool is formed of the spaces of solutions to the fundamental equations $FE^*(\nabla\nabla^*)$ and $FE^*(\nabla^*\nabla)$. We consider the split exact sequences of vector spaces\\
$$ 0\rightarrow\Omega^\nabla_2(M)\rightarrow J_{\nabla\nabla^*}\rightarrow \mathcal{S}^\nabla_2(M)\rightarrow 0,$$
$$ 0\leftarrow \Omega^{\nabla^*}_2(M)\leftarrow J_{\nabla^*\nabla}\leftarrow\mathcal{S}^{\nabla^*}_2(M)\leftarrow 0$$
There exist a pair $(\psi,\phi)\in J_{\nabla\nabla^*}\times\phi \in J_{\nabla*\nabla}$ such that
$$s^{*b}(\nabla,g) = dim(M) - rk(\Psi^*),$$
$$s^{*b}(\nabla^*) = dim(M) - rk(\Phi^*).$$
According to our technical machinery we get the following closed differential 2-forms
$$ \omega (X,Y) = g(\Psi^*(X),Y) = \frac{1}{2}[g(\psi(X),Y) - g(X,\psi(Y))],$$
$$ \omega^*(X,Y) = g(\Phi^*(X),Y) = \frac{1}{2}[g(\phi(X),Y) - g(X,\phi(Y))].$$
Those differential 2-fomrs satisfy the condition
$(a):\quad\nabla\omega = 0,$\\
$(b):\quad\nabla^*\omega^* = 0.$\\
Further
$$Ker(\omega) = Ker(\Psi^*),$$
$$Ker(\omega^*) = Ker(\Phi^*).$$
Furthermore the positivity of $g$ yields the orthogonal decompositions
$$ TM = Ker(\Psi^*) \oplus im(\Psi^*),$$
$$TM = Ker(\Phi^*) \oplus im(\Phi^*).$$
Now $\left\{(a),(b)\right\}$ is equivalent to
$$(a*):\quad \nabla*_X\Psi^*(Y) = \Psi^*(\nabla_XY),$$
$$(b*):\quad \nabla_X\Phi^*(Y) = \Phi^*(\nabla^*_XY).$$
Consequently  \\
(1) the distribution $im(\Psi^*)$ is $\nabla^*$-geodesic;\\
(2) the distribution $im(\Phi^*)$ is $\nabla$-geodesic;\\
(3) Since both $\nabla$ and $\nabla^*$ are torsion free $im(\Psi^*)$ and $im(\Phi^*)$ are completely integrable, (by Theorem of Frobenius)\\
In final the manifold $M$ is specially symplectically foliated by
$\left\{im(\Psi),\omega,\nabla^*\right\}$ and by $\left\{im(\Phi),\omega^*,\nabla\right\}$ This ends the demonstration $\clubsuit$\\
The theorem just demonstrated suggests another constructions of webs in statistical manifolds.

\newpage
\section{SOME HIGHLIGHTING CONCLUSIONS}
\textbf{Fascinating differential equations $FE(\nabla\nabla^*)$ and $FE^*(\nabla)$} \\
\textbf{First of all we have been fascined by the differential operators $D^{\nabla\nabla^*}$ and $D_\nabla$. Also the assessment of the range of the efficiency of the global analysis of the fundamental equations $FE(\nabla\nabla^*)$ and $FE^*(\nabla)$. That has been an unexpected feature. In this work we have addressed the problem $EX(\mathcal{S}$ in the global differential Geometry (viz Geometric structures in the sense of M. Gromov \cite{Gromov 2000}), the problem $EXF(\mathcal{S})$ in the differential topology, (essentially restricted foliations and webs) and the problem $DL(\mathcal{S})$ in the combinatoriel analysis (viz distancelike in the graph $\mathcal{GR}(LC)$). In view of the diversity of the structures $\mathcal{S}$ we have studied one would never have thought that the solutions to those problems would be deeply linked with only two same differential equations}\\
\subsection{The pair $(D^{\nabla\nabla^*},D_\nabla)$ and the pair of functions $(r^b,s^b)$}
Consider the gauge dynamic and the moduli space of gauge structures
$$g(M)\times TM\rightarrow TM,$$
$$ \mathbb{LC}(M) = \frac{\mathcal{LC}(M)}{J^1(Diff(M)}$$
The analysis of $(FE(\nabla\nabla^*),FE^*(\nabla))$ yields the map
$$(r^b\times s^b) : \mathbb{LC}(M)\times g(M)\rightarrow \mathbb{Z}^2.$$   
\subsection{The fundamental equation $FE^*(\nabla)$ in finite dimensional Lie groups}
The source of the gauge function $r^b$ is the global analysis of the fundamental equation $FE^{*}(\nabla)$. The purpose is to restrict the function $r^b$ to the moduli space of left invariant gauge structures. This restriction has been efficient in many circonstances: (i) the existence of the structure of Affine Lie Group; (ii) further it has been used for introducing the  semi-inductive system (respectively the semi projective system) of finite dimensional simply connected bi-invariant affine Lie groups and the semi-projective system of bi-invariant affine Cartan-Lie groups; (iii) an optimization of left invariant locally flat foliations in finite dimensional abstract Lie groups. We use those systems of affine Lie group for introducing Betti numbers of high order. However we do know how to interpret tohse Betti number $b_{p,q}\clubsuit$
\subsection{Hessian defects in abstract Lie groups}
We can also restrict the attention to the following left invariant data.\\
$$lf(\mathcal{G})\subset slc(\mathcal{G})\subset lc(\mathcal{G}),$$
Then what are named clan in finite dimensional KV algebras $\mathcal{A}$ are but left invariant Hessian structures in finite dimensional Lie groups \cite{Matsushima}, \cite{Shima(2)}. So the invariant $r^b(G)$ is a characteristic obstruction to $G$ admitting left invariant locally flat structures. Besides are two relative defects\\
(*) $r^b(G,g)$ is the Riemannian Hessian defect measuring how far from admitting left invariant Hessian structures is a left invariant Riemannian structure $(G,g)$ $\clubsuit$\\
(**) $r^b(G,\nabla)$ is the affine Hessian defect measuring how far from admitting left invariant Hessian structure is an affine Lie group $(G,\nabla)$ $\clubsuit$
\subsection{The fundamental equation $FE^*(\nabla\nabla^*)$ and bi-invariant Riemannian geometry in finite dimensional Lie groups}
The source of the function $s^b$ is the fundamental equation $FE^*(\nabla\nabla^*)$.\\
Regarding the problems $ EX(\mathcal{S})$, $EXF(\mathcal{S})$ and $DL(\mathcal{S})$ in the finite dimensional symplectic geometry the function $s^b$ is an alternative to the function $r^b$. We have emphasized its impacts on some outstanding open problems.\\
(i): The existence of bi-invariant Riemannian metrics in finite dimensional Lie groups $\clubsuit$\\
(ii): An optimization of bi-invariant Riemann foliations in finite dimensional Lie groups $\clubsuit$\\
\textbf{SOME INSIHGTS.}\\
$s^b(G) = 0$ for all finite dimensional semi simple Lie groups $\clubsuit$\\
$s^b_+(G) = 0$ for all finite dimensional compact Lie groups $\clubsuit$\\
$s^b_+(G) > 0$ for all finite dimensional non-unimodular Le groups $\clubsuit$\\
\textbf{THE SIZE.}\\
It must be noticed that our machinery that involves the function $s^b$ yields all bi-invariant Riemannian metrics and all bi-invariant Riemannian foliations in finite dimensional Lie groups.
\subsection{The fundamental equation $FE^*(\nabla)$, the homogeneous Kaehlerian geometry and some canonical affine representations}
We have introduced an alternative to the Gindikin-Piatecii-Shapiro-Vinberg fibrations of homogeneous Kaehlerian manifolds $(M,g,J)$ (the fundamental conjecture). The fibers of our fibration are complex cylinders which are  biholomorphic to the same complex cylinder over a complex twisted torus the first Betti number of which is $r^b(\nabla)$. The context is the algebraic topology of homogeneous Kaehlerian manifolds. Here $\nabla$ is the Levi-civita connection of $(M,g,J)$. By the way the foliation $\mathcal{F}_\nabla$ is a solution to $EXF(\mathcal{S})$ for GFHSDS-foliation in $(M,g,J)$. A notable fact is the existence of canonical unitary affine representation of isotropy subgroups of homogeneous Kaehlerian manifolds.
\subsection{The fundamental equation  $FE(\nabla\nabla^*)$ and the symplectic geometry}
Up to today there  was no effecient invariants measuring how far from admitting symplectic structures is a manifold. There was none helping to decide whether a manifold is symplectically foliated. It is a highlingting fact that the anwers to those questions are furnished by a combination of the functor of Amari and the global analysis of the differential equation $FE^*(\nabla,\nabla^*)$ which gives rise to the function $s^b$. That function is used for investigating the existence of symplectic structures in manifolds and the existence of left invariant symplectic structures in finite dimensional Lie groups.\\
The case of Lie group is a byproduct of the machinery used in the case of differentiable manifolds. This approach is an impact on the symplectic geometry of the information geometry and of its methods.
\subsection{Canonical representations of fundamental groups}
This leads to conditions for a discrete group being fundamental groups of manifolds carrying special gauge structures.
\subsection{In final}
\textbf{The global analysis of the fundamental equation $FE^*(\nabla\nabla^*)$ impacts the symplectic geometry and the bi-invariant Riemannian geometry in Lie. It also deeply impacts the differential topology through the theory of webs in statistical manifolds}.\\
\textbf{The global analysis of the fundamental equation $FE^*(\nabla)$ impacts the locally flat geometry, the Hessian geometry, the geometry of Koszul and the information geometry}.
\section{NEW PERSPECTIVES OF APPLICATIONS}
This short section is devoted to focus on a few new perspectives. The purpose is to address the problems $EX(\mathcal{S})$ and $EXF(\mathcal{S})$ in complex systems, (e.g. in Big Data.) From the fundamental mathematics point of view the problems $EX(\mathcal{S})$, $EXF(\mathcal{S})$ and $DL(\mathcal{S})$ have been the prior motivations of this research work. However we have mentioned other source of motivations which are linked with the applied information geometry in complex systems. This section is devoted to overview those aspects of this work.
\subsection{New invariants of statistical models}
Among notable classical invariants of statistical models are the family of connections of Chentsov as well as their their links with the Fisher information. To each $\alpha$-connecion $\nabla^\alpha$ we have assigned many objects having homological nature such as semi-inductive systems of bi-invariant affinely flat Lie groups. Those objects and their avatars form new statistical invariants. Their roles in the applied information geometry are to be investigated. By the global analysis of the differential equation $FE^*(\nabla)$ local statistical models as in \cite{Amari-Nagaoka} are (affine) homogeneous. Therefore they are nothing but open orbits of affine transformations group as in \cite{Koszul(2)}. By the positivity of the Fisher information those models are homogeneous bounded domains. The machineries we have developed significantly enrich the (differential) topolopgy of regular statistical models as well as well as the topology of statistical manifolds. Conider a regular statistical model $\mathbb{M}$ whose base manifold is denoted by $M$. We have the family of fundamental differential equations
$$\mathbb{R} \ni \alpha\rightarrow FE(\alpha,-\alpha) = FE(\nabla^\alpha\nabla^{-\alpha})$$
The solution to every fundammenatal equation $FE(\alpha,-\alpha)$ is linked with the differential topology of $M$ according to the split short exact sequences
$$ 0\rightarrow\Omega^\alpha_2(M)\rightarrow J_{\alpha,-\alpha}\rightarrow \mathcal{S}^\alpha_2(M)\rightarrow 0, $$
$$ 0\rightarrow\Omega^{-\alpha}_2(M)\rightarrow J_{-\alpha,\alpha}\rightarrow\mathcal{S}^{-\alpha}_2(M)\rightarrow0.$$
Here elements of $\Omega^\alpha_2(M)$ are symplectic foliations in $M$ while those of $\mathcal{S}^\alpha_2(M)$ are Riemannian foliations. We know more about solutions to $FE(\alpha,-\alpha)$. Indeed another feature is that every solution $\phi \in J_{\alpha,-\alpha}$ defines another unique $\Phi \in J_{\alpha,-\alpha}$ such that $Ker(\Phi)$ is a regular $\nabla^\alpha$-geodesic foliation, $Im(\Phi)$ is a regular $\nabla^{-\alpha}$-geodesic foliation and $$TM = Ker(\Phi) \oplus Im(\Phi)$$
Those reminders show how the integrability problem for the differential equation $FE(\alpha,-\alpha)$ is linked with the topology of $M$. They also show that the analysis of those  statistical informations which are encoded by the Ficher information is linked with the topology of $M$. This fact strenghtens the role of the differential operators $D^{\nabla^\alpha\nabla^{-\alpha}}$ in applied statistics and in the applied information geometry.
\subsection{Complex systems and digital science}
We regard the information geometry as a functor from the category of measurable sets (or measurable categories) to the category of statistical models  \cite{Nguiffo Boyom(6)}. We regard the topology of information as a functor from the category of measurable sets to homological algebra \cite{Baudot-Bennequin}, \cite{Gromov}. In \cite{Nguiffo Boyom(6)} we have answered the fundamental question raised by P. McCullagh \cite{McCullagh}: \textbf{What is a statistical model?} See also \textbf{In a search of a structures} by Gromov \cite{Gromov}. The information geometry involves the group of symmetries of measurable sets. We observe that in pure or in applied statistics, informations would become all relevanter as they are linked with an appropiate geometric structure of the measurable set. This observation is of fundamental importance to explore complex systems. Doubtless this was the feeling of both P. McCullagh \cite{McCullagh} and M. Gromov \cite{Gromov}. The digital science (or computing science) is the ingeneering step after conceptual theorems have assured the existence of solutions to the fundamental problems such as $EX(\mathcal{S})$, $EXF(\mathcal{S})$ and $DL(\mathcal{S})$. So the prior is the question whether a complex system (,Big Data,) admits an appropiate geometric structure $\mathcal{S}$. In a singular statistical model $\mathbb{M}$ those informations which are encoded by the Fisher information $g_\mathbb{M}$ are supported by the transverse structure to foliation $Ker(g_\mathbb{M})$. Generically challenges in applied statistics are closely linked with the problem $EXF(\mathcal{S})$. In final challenges in applied statistics are related with both $EX(\mathcal{S})$ and $EXF(\mathcal{S})$. Actually what are callenges in complex systems? Let us magnify a McCullagh pattern. A Big Datum is regarded as a measurable set (or a measurable category) $(\Xi,\Omega)$ whose group of measurable isomorphisms is denoted by $\Gamma$ (see the appendix to this paper). A statistical model for $(\Xi,\Omega)$ is a locally trivial $\Gamma$-fibration $\mathbb{M}$ whose fibers are isomorphic to $\Xi$. Lossely speaking the information geometry is the $\Gamma$-geometry in $\mathbb{M}$.\\  
\subsection{Needs of geometric structures in complex systems: new challenges}
The role played by the digital science in complex systems has an ingineering nature. The performance of the digital science is based on solutions to the problems $EX(\mathcal{S})$, $EXF(\mathcal{S})$ and $DL(\mathcal{S})$. We go to simulate three patterns. \\  
(1) Depending on needs and purposes, the extraction of relevant $\Gamma$-Invariants (in $\mathbb{M}$) may depend on restricted $\Gamma$-invariant geometric structures in $\Xi$. This need involves solutions to the problem $EX(\mathcal{S})$ in $\Xi$. \\   
(2) A dreamed structure $\mathcal{S}$ may not exit in $\Xi$. Nevertheless it may exist a partition of $\Xi$ which agrees with $\mathcal{S}$. This occurence is linked with the problem $EXF(\mathcal{S})$. So is the domain of applied statistics (polls). \\  
(3) Arises the question whether a structure $\mathcal{S}$ agrees with other prior parameters. So one may be facing the evaluation of specific Gaps. Such a situtaion is linked with the problem $DL(\mathcal{S})$. \\
\subsection{BigData-Mathematic-StatisticalModels-HighPerformanceComputing}
Challenges at interface of < Big Data- Mathematic- Statistical Models> are objects of the Information Geometry. Those challenges are linked with the $EX(\mathcal{S})$.\\
Challenges at interface of < DigitalScience-Mathematic- BigData > are linked with Topology-Geometric structures of BigData. Those challenges are linked with problems $EX(\mathcal{S})$ and $EXF(\mathcal{S})$\\
Challenges at interface of < BigData-Mathematic-HighPerformanceComputing > are linked with the problems $EXF(\mathcal{S})$ and $DL(\mathcal{S})$.
\subsection{The fundamental equations $<<FE(\nabla\nabla^*)$, $FE^*(\nabla) >>$ and Complex systems}
\textbf{In final Topology-Geometric structures which arise in the Information Geometry are linked with gauge structures via their global analysis and their geometry. Instances are Hessian structures, symplectic structures, bi-invariant Riemannian structures in Lie groups. Thus through the pair of functions $(r^b,s^b)$ the fundamental equations $FE^*(\nabla)$ and $FE(\nabla\nabla^*)$ are useful for analysing topology-geometric structures of complex systems. That analysis might enrich the high performance computing $\clubsuit$}

\newpage
\section{APPENDIX.}

The function $r^b$ has been helpfull to solve the open problem of complexity of models, viz how far from being an exponential familly is a statistical model.\\
The theory of statistical models has been re-established in \cite{Nguiffo Boyom(6)} where reders will find details. This appendix is devoted to overview the core of that re-establised theory. The axioms of definition of statistical models are illustrated with the diagrams given below.\\
\subsection{THE CATEGORY $\mathcal{FB}(\Gamma,\Xi)$}
We go to recall the framework of the theory of statistical models, \cite{Nguiffo Boyom(6)}
\subsubsection{THE OBJECTS OF $\mathcal{FB}(\Gamma,\Xi).$}
$(\Xi,\Omega)$ is a measurable set and $\Gamma$ is the group of measurable ismorphisms the the set 
$\Xi$.
\begin{defn} An object of the category $\mathcal{FB}(\Gamma,\Xi)$ is a quadruple $[\mathcal{E},\pi,M,D]$ which is composed of a set $\mathcal{E}$, a differentiable manifold, a Koszul connection in $M$ $D$ in$M$ and a projection $\pi$ of $\mathcal{E}$ onto $M$. Below are the links between those items. \\
$(1):$ The map $\pi$ is a locally trivial$\Gamma$-bundle
$$\pi: \mathcal{E}\rightarrow M.$$
The fibers $\mathcal{E}_x$ are isomorphic to the set $\Xi$.\\
$(2):$ The couple $(M,D)$ is an $m$-dimensional locally flat manifold.\\
$(3):$ There is a group action
$$ \Gamma\times [\mathcal{E}\times M]\times \mathbb{R}^m \ni (\gamma,[e,x,\theta])\rightarrow [\gamma.e),\gamma.x,\tilde{\gamma}.\theta]\in [\mathcal{E}\times M\times\mathbb{R}^m ]. $$
This action is subject to the following compatibility requirement\\
$$ \pi(\gamma.e) = \gamma.\pi(e)\quad \forall e \in \mathcal{E}.$$
$(4):$ Every point $x\in M$ has an open neighborhood $U$ which is the domain of a local fiber chart $$\Phi_U\times\phi_U : [\mathcal{E}_U \times U ]\ni (e_x,x)\rightarrow [\Phi_U(e_x),\phi_U(x)]\in [\mathbb{R}^m \times \Xi]\times \mathbb{R}^m.$$
Further $\Phi$ is linked with $\phi$ as it follows \\
$(4.1): (U,\phi_U)$ is an affine local chart of the locally flat manifold $(M,D),$\\
$(4.2): \phi_U(\pi(e)) = p_1(\Phi_U(e)).$\\
$(5):$  Put $$\Phi_U(e) = (\theta_U(e),\xi_U(e))\in \mathbb{R}^m \times \Xi,$$
consider two local charts $(U,\Phi\times\phi)$ and $(U^*,\Phi^*\times\phi^*)$ and assume that
$$U\cap U^* \neq \emptyset;$$
then there exists a unique $\gamma_{UU^*}\in \Gamma $
such that
 $$[\gamma_{UU^*}.\Phi](e) = \Phi^*(e)\quad \forall e\in \mathcal{E}_{U\cap U^*} \clubsuit $$
\end{defn}

\begin{figure}
 \resizebox{1.0\linewidth}{!}{
\centering
  \begin{tikzpicture}[font=\tiny]
    \coordinate (a) at (0,0);
    \coordinate (b) at (2,0);
    \coordinate (c) at (2,1);
    \coordinate (d) at (0,1);

    \foreach \x in {a,b,c,d}{
    \node[point] at (\x) {};
    }
    \foreach \x/\y in {a/E,b/E}{
    \node[] at ($(\x)+(0,+1.2)$) {$\y$};
    }

    \foreach \x/\y in {c/M,d/M}{
    \node[] at ($(\x)+(0,-1.2)$) {$\y$};
    }
    \draw (a) -- (b) -- (c) -- (d) --cycle;
    \draw (a) -- ($(a)!0.5!(b)$) node[above] {\(\gamma\)};
     \draw (c) -- ($(c)!0.5!(d)$) node[above] {\(\gamma\)};
     \draw (d) -- ($(a)!0.5!(d)$) node[left] {\(\pi\)};
     \draw (b) -- ($(b)!0.5!(c)$) node[right] {\(\pi\)};
  \end{tikzpicture}
  }
\caption{Fibration} \label{fig:fibration}
\end{figure}
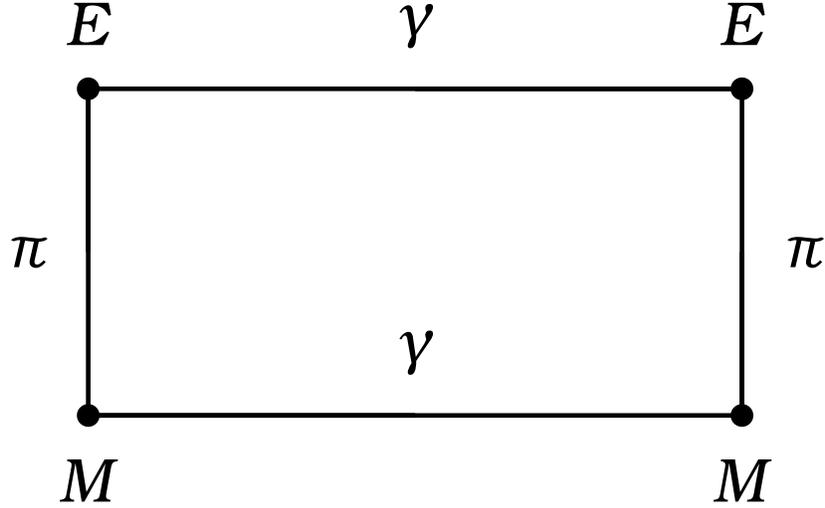
\textbf{COMMENTS.} \\
Figure 1 expresses that the fibration $\pi$ is $\Gamma$-equivariant
Requirement (3) and (4) mean that
$$[\Phi_U(e),\phi_U(\pi(e)] = [[\theta_U(e),\xi_U(e)],\theta_U(e)]$$
Both requirements (4) and (5) yield the following remarks,\\
the following action is differentiable
$$ \Gamma \times M \ni (\gamma,x)\rightarrow \gamma.x \in M,$$
the following action is an affine action
$$\Gamma\times\mathbb{R}^m \ni (\gamma,\theta)\rightarrow \tilde{\gamma}.\theta,$$
the left side member and the right side member of (5) have the following meanings.
$$\gamma_{UU^*}.[\theta_U(e),\xi_U(e)] = [\theta_{U^*}(e),\xi_{U^*}(e)]$$
Consequently (5) implies that for all $x \in U\cap U^*$ one has
$$ \tilde{\gamma}_{UU^*}.\phi(x) = \phi^*(x)$$
Therefore we get
$$\gamma_{UU^*} = \phi^*\circ\phi^{-1};$$
Consider three domains of local chart wich intersect
$U$, $u^*$ and $U^{**}$ 
$$U\cap U^*\cap U^{**} \neq \emptyset$$
then
$$ \gamma_{UU^*}\circ \gamma_{UU^*} = \gamma_{UU^{**}}$$
The requirement (3) means that the fibration $\pi$ is $\Gamma$ equivariant.\\
That ends the comments.

\begin{figure}
 \resizebox{1.0\linewidth}{!}{
\centering
    \begin{tikzpicture}[font=\tiny]
    \coordinate (a) at (0,0);
    \coordinate (b) at (2,0);
    \coordinate (c) at (2,2);
    \coordinate (d) at (0,2);
    %
    %\coordinate (a1) at (0.65,0.45);
    \coordinate (b1) at (2.65,0.45);
    \coordinate (c1) at (2.65,2.45);
    %\coordinate (d1) at (0.65,2.45);

    \foreach \x in {a,b,c,d,b1,c1}{
    \node[point] at (\x) {};
    }
    \foreach \x/\y in {c/B_1,d/A_1,c1/C_1}{
    \node[] at ($(\x)+(0,0.2)$) {$\y$};
    }
    \foreach \x/\y in {a/A,b/B,b1/C}{
    \node[] at ($(\x)+(0,-0.2)$) {$\y$};
    }
    \draw (a) -- (b) -- (c) -- (d) --cycle;
    \draw (a) -- (b) [arrowMe=>] ;
    \draw (b) -- (b1) [arrowMe=>] ;
	\draw (c) -- (b) [arrowMe=>] ;
    \draw (b)--($(c)!0.5!(b)$) node[right] {\(p_1\)};
    \draw (c1) -- (b1) [arrowMe=>] ;
    \draw (d)--($(d)!0.5!(a)$) node[left] {\(\pi\)};
    \draw (d) -- (a)[arrowMe=>] ;
    \draw (c1)--($(c1)!0.5!(b1)$) node[right] {\(p_1\)};
    \draw (b) -- (b1) ;
    \draw (c) -- (c1)[arrowMe=>] ;
    \draw (a) -- (b1) [arrowMe=>] ;
    \draw (d) -- (c1) [arrowMe=>]  ;
	\draw (a) -- ($(a)!0.5!(b)$) node[below] {\(\phi_u\)};
    \draw (d) -- ($(d)!0.5!(c1)$) node[above] {\(\Phi_{u^\star}\)};
	\draw (d) -- ($(d)!0.5!(c)$) node[below] {\(\Phi_u\)};
	\draw (c) -- ($(c)!0.5!(c1)$) node[below] {\(\gamma_{uu^\star}\)};	  \draw (b) -- ($(b)!0.5!(b1)$) node[below] {\(\gamma_{uu^\star}\)};
    %\draw ($(b1)!0.5!(c1)$) -- (b1) node[right] {\(p_i\)};
	\draw (a) -- ($(a)!0.5!(b1)$) node[above] {\(\phi_{u^\star}\)};
    %\draw ($(b1)!0.5!(c1)$) -- (b1) node[right] {\(p_i\)};
\end{tikzpicture}
}
\caption{Equivariance} \label{fig:equivariance}
\end{figure}
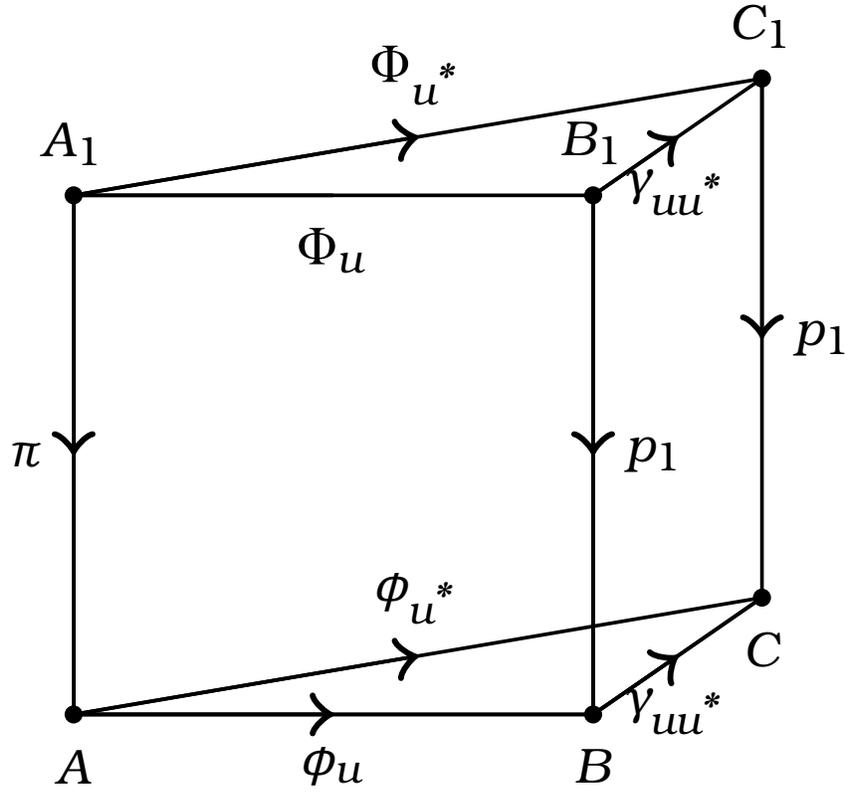

\subsubsection{THE MORPHISMS OF $\mathcal{FB}(\Gamma,\Xi).$}
Without the express statement of the contrary we deal with $\Gamma$-equivariant fibrations
$$\pi : \mathcal{E} \rightarrow M.$$
Let $[\mathcal{E},\pi,M,D]$ and $[\mathcal{E}^*,\pi^*,M^*,D^*]$ be two objects of
$\mathcal{FB}(\Gamma,\Xi)$. \\
Let $\Psi\times\psi$ be a map
$$ [\mathcal{E}\times M] \ni (e,x)\rightarrow (\Psi(e),\psi(x))\in [\mathcal{E}^*\times M^*] $$
\begin{defn} A couple $(\Psi\times\psi)$ is a morphism of the category $\mathcal{FB}(\Gamma,\Xi)$ if the following properties are satisfied \\
$(m.1): \pi^*\circ \Psi = \psi\circ \pi $\\
$(m.2):$ both $\Psi$ and $\psi$ are $\Gamma$-equivariant isomorphisms,
that is to say
$$\Psi(\gamma.e) = \gamma.\Psi(e),$$
$$\psi(\gamma.x) = \gamma.\psi(x) $$
$(m.3): \psi$ is an affine map of $(M,D)$ in $(M^*,D^*)$ $\clubsuit$
\end{defn}

\begin{figure}
 \resizebox{1.0\linewidth}{!}{
\centering
  \begin{tikzpicture}[font=\tiny]
    \coordinate (a) at (0,0);
    \coordinate (b) at (2,0);
    \coordinate (c) at (2,2);
    \coordinate (d) at (0,2);
    \coordinate (a1) at (0.65,0.45);
    \coordinate (b1) at (2.65,0.45);
    \coordinate (c1) at (2.65,2.45);
    \coordinate (d1) at (0.65,2.45);

    \foreach \x in {a,b,c,d,a1,b1,c1,d1}{
    \node[point] at (\x) {};
    }
    \foreach \x/\y in {c/E^*,d/E,c1/E^*,d1/E}{
    \node[] at ($(\x)+(0,0.2)$) {$\y$};
    }
    \foreach \x/\y in {a/M,b/M^*,b1/M^*,a1/M}{
    \node[] at ($(\x)+(0,-0.2)$) {$\y$};
    }

    \draw (a) -- (b) -- (c) -- (d) --cycle;
    \draw (d1) -- (c1) -- (b1) edge[densely dashed] (a1) ;
    \draw[densely dashed] (a) -- (a1) -- (d1) edge[solid] (d);
    \draw (b) -- (b1)  (c) -- (c1) ;

    \draw (b) -- ($(b)!0.5!(c)$) node[right] {\(\pi^\star\)};
    \draw (a) -- ($(a)!0.5!(a1)$) node[left] {\(\gamma\)};
    \draw (a) -- ($(a)!0.5!(d)$) node[left] {\(\pi\)};
    \draw (d) -- ($(d)!0.5!(d1)$) node[right] {\(\gamma\)};
    \draw (d1) -- ($(d1)!0.5!(c1)$) node[above] {\(\Phi\)};
    \draw (b) -- ($(b)!0.5!(b1)$) node[left] {\(\gamma\)};
    \draw (d) -- ($(d)!0.5!(c)$) node[below] {\(\Phi\)};
    \draw (a) -- ($(a)!0.5!(b)$) node[below] {\(\phi\)};
 	\draw (a1) -- ($(a1)!0.5!(b1)$)[densely dashed] node[below] {\(\phi\)};

	\draw (a) -- (b)[arrowMe=>];
    \draw (d) -- (c)[arrowMe=>] ;
    %\draw[densely dashed] (b) -- (a1) -- (c1) edge[solid]node[point,outer sep=0pt] (c1b){}  (b);
    %\draw[densely dashed] (a1) -- coordinate[pos=0.95] (arc1) (c1b);
    %\draw (c1b) -- coordinate[pos=0.1] (arc2)(b1);
    %\draw (arc1) edge[in=210,out=290] (arc2);
    %\node[] at ($(c1b)+(0.15,0)$) {$E$};
    %\begin{scope}[on background layer]
    %\path[fill=gray!20] (a) -- (a1) -- (d1) -- (d) -- cycle;
    %\path[fill=gray!20] (a1) -- (c1) -- (b) -- cycle;
    %\end{scope}
  \end{tikzpicture}
  }
\caption{Moduli} \label{fig:moduli}
\end{figure}

In the next subsection we go to define the category of statistical models for a measurable set $(\Xi,\Omega)$. The framework is the category $\mathcal{FB}(\Gamma,\Xi)$.

\subsection{THE CATEGORY $\mathcal{GM}(\Xi,\Omega).$}
We keep the notation used in the preceding subsections. Our concern is the category $\mathcal{GM}(\Xi,\Omega)$ whose objects are statistical models for a measurable set $\Xi,\Omega)$.
\subsubsection{THE OBJECTS OF $\mathcal{GM}(\Xi,\Omega)$}
\begin{defn} An object of the category $\mathcal{FB}(\Gamma,\Xi)$, namely $[\mathcal{E},\pi,M,D]$ is an $m$-dimensional statistical model for $(\Xi,\Omega)$ if $M$ is an $m$-dimensional manifold and if the following requirements $[\rho_\star]:$ are satisfied\\
$[\rho_1]$: for every local chart $(U,\Phi_U\times\phi_U)$ the subset
$$(\Theta_U\times\Xi)= \Phi_U(\mathcal{E}_U )$$
supports a non negative real valued function $P_U$ having the following properties\\
$[\rho_{1.1}]:$ for every fixed $\xi\in \Xi$ the function
$$\Theta_U \ni \theta \rightarrow P_U(\theta,\xi)$$
is differentiable; \\
$[\rho_{1.2}]:$ for every fixed $\theta \in \Theta_U$ the triple
$$(\Xi,\Omega,P_U(\theta,-)$$
is a probability space; furthermore the operation of integration $\int_\Xi $ commutes with the operation of differentiation $d_\theta = \frac{d}{d\theta}$;\\
$[\rho_{1.3}]:$ let $(U,\Phi_U\times\phi_U,P_U)$ and $(U^*,\Phi_{U^*}\times\phi_{U^*},P_{U^*})$ be as in $[\rho_{1.1}]$ and in $[\rho_{1.2}]$,\\
if $U\cap U^* \neq \emptyset$ then $P_U$, $P_{U^*}$ and
$$\gamma_{UU^*} = \phi_{U^*}\circ \phi^{-1}_U$$
are linked as it follows
$$P_{U^*}\circ\gamma_{UU^*} = P_U.$$
$[\rho_{1.4}]:$ Let $U\subset M$ be an open subset and let $\gamma \in \Gamma$; \\
we assume that both $U$ and $\gamma.U$ are domains of local charts
$$(U,\Phi_{U}\times\phi_{U},P_{U})$$
and
$$(\gamma.U,\Phi_{\gamma.U}\times\phi_{\gamma.U},P_{\gamma.U})$$
satisfying the requirements $\rho_{1.1}$, $\rho_{1.2}$ and $\rho_{1.3}$; then the relations
$$\Phi_{\gamma.U}\circ\gamma = \gamma\circ \Phi_U,$$
$$\phi_{\gamma.U}\circ\gamma = \gamma\circ \phi_U,$$
imply the relation
$$ P_{\gamma.U}\circ\gamma = P_U \clubsuit $$
\end{defn}

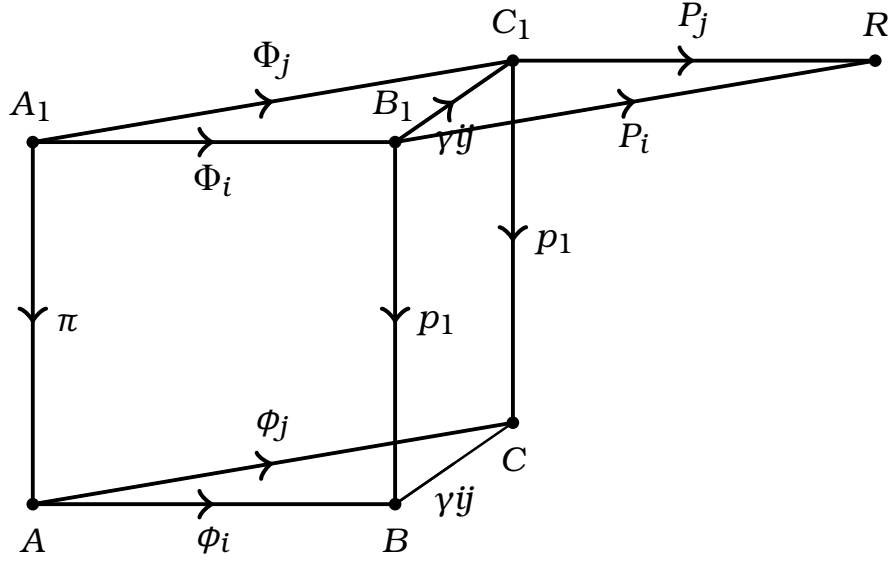
\begin{figure}
\resizebox{1.0\linewidth}{!}{
\centering
   \begin{tikzpicture}[font=\tiny]
    \coordinate (a) at (0,0);
    \coordinate (b) at (2,0);
    \coordinate (c) at (2,2);
    \coordinate (d) at (0,2);
    %
    %\coordinate (a1) at (0.65,0.45);
    \coordinate (b1) at (2.65,0.45);
    \coordinate (c1) at (2.65,2.45);
    \coordinate (r) at (4.65,2.45);
    %\coordinate (d1) at (0.65,2.45);

    \foreach \x in {a,b,c,d,b1,c1,r}{
    \node[point] at (\x) {};
    }
    \foreach \x/\y in {c/B_1,d/A_1,c1/C_1,r/R}{
    \node[] at ($(\x)+(0,0.2)$) {$\y$};
    }
    \foreach \x/\y in {a/A,b/B,b1/C}{
    \node[] at ($(\x)+(0,-0.2)$) {$\y$};
    }
    %\node[] at ($4$) {$R$};
    \draw (a) -- (b) -- (c) -- (d) --cycle;
    \draw (c1) -- (b1) ;
    \draw (b) -- (b1)  (c) -- (c1) ;
    \draw (a) -- (b1);
    \draw (d) -- (c1);
    \draw (c1) -- (r) (c) --  (r);
    %\draw (a) -- coordinate[midway](m) (m) node[below] {\(\lambda_1\)};

    %\tkzDrawSegments[arrowMe=>](a,c)
    \tkzDrawSegments[arrowMe=>](d,a)
    \tkzDrawSegments[arrowMe=>](d,c)
    \tkzDrawSegments[arrowMe=>](c,c1)
    \tkzDrawSegments[arrowMe=>](c1,b1)
    \tkzDrawSegments[arrowMe=>](c,b)
     \tkzDrawSegments[arrowMe=>](a,b)
	\tkzDrawSegments[arrowMe=>](a,b1)
	\tkzDrawSegments[arrowMe=>](d,c1)
	\tkzDrawSegments[arrowMe=>](c,r)
	\tkzDrawSegments[arrowMe=>](c1,r)
    %Les names des figures

    \draw (a) -- ($(a)!0.5!(b)$) node[below] {\(\phi_i\)};
	\draw (d) -- ($(d)!0.5!(a)$) node[right] {\(\pi\)};
    \draw (d) -- ($(d)!0.5!(c)$) node[below] {\(\Phi_i\)};
    \draw (c) -- ($(c)!0.5!(c1)$) node[below] {\(\gamma{ij}\)};
    \draw (b) -- ($(b)!0.5!(b1)$) node[below] {\(\gamma{ij}\)};
    \draw (a) -- ($(a)!0.5!(b1)$) node[above] {\(\phi_j\)};
    \draw (d) -- ($(d)!0.5!(c1)$) node[above] {\(\Phi_j\)};
    \draw (c) -- ($(c)!0.5!(r)$) node[below] {\(P_i\)};
	\draw (c1) -- ($(c1)!0.5!(r)$) node[above] {\(P_j\)};
    \draw (c1) -- ($(c1)!0.5!(b1)$) node[right] {\(p_1\)};
    \draw (c) -- ($(c)!0.5!(b)$) node[right] {\(p_1\)};
  \end{tikzpicture}
  }
\caption{Localisation} \label{fig:localisation}
\end{figure}

\textbf{A COMMENT.}\\
Actually ($[\rho_{1.3}]$) has the following meaning:
$$P_{U^*}[\tilde{\gamma}_{UU^*}.\theta_U(e),\gamma_{UU^*}.\xi_U(e)] =
P_U(\theta_U(e),\xi_U(e))$$
$\forall e \in \mathcal{E}_{U\cap U^*}$. This ends the comment.\\
\begin{defn}
A quadruple $[U,\Phi_U\times\phi_U,P_U,\gamma_{UU^*}]$ (as in the last definition) is called a statistical local chart of the statistical model $[\mathcal{E},\pi,M,D]$  $\clubsuit$
\end{defn}
\textbf{WARING.}\\
To introduce the morphisms of the category $\mathcal{GM}(\Xi,\Omega)$ we firstly plan defining a relevant geometry invariant which does not depend on a fiber atlas $(U_j,\Phi_j\times \phi_j, P_j,\gamma_{ij})$. That is the aim of the next subsubsection.
\subsubsection{THE GLOBAL PROBABILITY DENSITY OF A STATISTICAL MODEL.}
We consider a complete (or maximal statistical) atlas $A_\Phi$ of an object $[\mathcal{E},\pi,M,D]$ of the category $\mathcal{GM}(\Xi,\Omega)$. In detail 
$$A_\Phi = [U_j,\Phi_j,\phi_j, P_j,\gamma_{ij}]$$
The family $U_j$ is an open covering of $M$ and $\mathcal{E}_j\times U_j$ is the domain of the fibered chart $(\Phi_j\times \phi_j)$. Of course $\mathcal{E}_j$ stands for
$\mathcal{E}_{U_j}$.\\
If $U_i\cap U_j \neq \emptyset$ then one has
$$\phi_j(x) = \tilde{\gamma}_{ji}.\phi_ij(x)\quad \forall x\in U_i\cap U_j $$
In particular $A = (U_j,\phi_j)$ is an affine atlas of the locally flat manifold $(M,D)$.\\
One has
$$\Phi_j(\mathcal{E}_{y^*}) = \phi_j(y^*) \times \Xi \quad \forall y^* \in U_j.$$
We set
$$[\mathcal{E}_{y^*},\Omega_{y^*}] = \Phi^{-1}_j [\phi_j(y^*)\times(\Xi,\Omega]$$
The atlas $A_\Phi$ satisfies requirements $[\rho_{1.1}]$, $[\rho_{1.2}]$, $[\rho_{1.3}]$.\\
In $\mathcal{E}_{U_j}$ we define the local function $p_j$ by 
 $$p_j = P_j\circ \Phi_j.$$
We suppose that the following condition holds,
$$U_i\cap U_j \neq \emptyset $$
then by the virtue of $[\rho_{1.3}]$ one has
$$ p_i(e) = p_j(e)\quad \forall e\in \mathcal{E}_{U_i\cap U_j}$$.
Thereby there exists a unique globally defined non negative real function
  $$ \mathcal{E} \ni e\rightarrow p(e)\in \mathbb{R}$$
such that
$p_j $ coincides with the restriction to $\mathcal{E}_j$ of $p$. Furthermore, let $p_x$ be the restriction to a fiber $\mathcal{E}_x$ of the function $p$, then the triple
$$(\mathcal{E}_x,\Omega_x,p_x) $$
is a probability space.\\
\begin{defn} The function $$\mathcal{E}\ni e \rightarrow p(e) \in \mathbb{R}$$
is called the probability density of the model $[\mathcal{E},\pi,M,D]$
\end{defn}

\begin{figure}
 \resizebox{.9\linewidth}{!}{
\centering
   \begin{tikzpicture}[font=\tiny]
    \coordinate (a) at (0,0);
    \coordinate (b) at (2,0);
    \coordinate (c) at (2,1);

    \foreach \x in {a,b,c}{
    \node[point] at (\x) {};
    }
    \foreach \x/\y in {a/$$E_{i}$$,b/E}{
    \node[] at ($(\x)+(0,-0.2)$) {$\y$};
    }
    \foreach \x/\y in {c/$$\mathbb{R}$$}{
    \node[] at ($(\x)+(0,+0.2)$) {$\y$};
    }

    \draw (a) -- (b) -- (c) --cycle;
    \tkzDrawSegments[arrowMe=>](a,b)
     \draw (a) -- ($(a)!0.5!(c)$) node[above] {\(P_i\)};
     \draw (b) -- ($(b)!0.5!(c)$) node[right] {\(p\)};
    %\tkzDrawSegments[>=stealth,arrowMe=>>](a,c)
    \tkzDrawSegments[arrowMe=>](a,c)
    \tkzDrawSegments[arrowMe=>](b,c)
    %\mark(segment(A,D),cross) ;
    %\draw[->>] (a) -- ($((b)+(a))/2$);
  \end{tikzpicture}
  }
\caption{Density} \label{fig:density}
\end{figure}

\textbf{COMMENTS.}\\
$(i):$ We take into account the global probability density $p$. Then an object of $\mathcal{GM}(\Xi,\Omega)$ is denoted by
$$[\mathcal{E},\pi,M,D,p].$$
$(ii):$ The function $p$ is $\Gamma$-equivariant. THIS IS THE GEOMETRY in the sense of Felix Klein.\\
$(iii):$ Really to construct the present theory of models we have not involved any finite dimensional argument. Up to intrinsic difficulties of the infinite dimensional locally flat geometry, our construction walks in any dimension.\\
$(iv)):$ Our formalism cannot be the object of critiques and questions similar to those arisen from the attempt of McCullagh as in \cite{McCullagh}. Indeed depending on concerns and on needs our formalism yields any suitable type of model. The crucial problem is of topological nature. \\
The challenge consists of constructing an open covering of the base manifold, namely
$$\mathcal{U} = \left\{U_j,j\in J\right\}$$
whose simplicial complex supports the following two data. \\
(a) A real valued random 0-chain
$$\left\{p_j:\mathcal{E}_j\rightarrow [0,1],\quad j \in J\right\}.$$
(b) A $\Gamma$-valued 1-cocycle
$$\left\{\gamma_{ij}, i,j \in J\right\},$$
both being subject to the requirements $\rho_\star $. The choice of the 0-cocycle depends on the nature ths measurable set $\Xi,\Omega)$, e.g. Agriculture, Heat, Population? Finance. The choice also depends on the stochastic variables ad hoc.  This ends comments $\clubsuit$\\
We go to describe the morphisms of the category $\mathcal{GM}(\Xi,\Omega)$
\subsubsection{THE MORPHISMS OF $\mathcal{GM}(\Xi,\Omega)$}
\begin{defn}
Let $[\mathcal{E},\pi,M,D,p]$ and $[\mathcal{E}^*,\pi^*,M^*,D^*,p^*]$ be two objects
of the category $\mathcal{GM}(\Xi,\Omega)$. A morphism of the category $\mathcal{FB}(\Gamma,\Xi)$
$$(\Psi\times\psi): [\mathcal{E},\pi,M,D] \rightarrow [\mathcal{E}^*,\pi^*,M^*,D^*]$$ 
is a morphism of $[\mathcal{E},\pi,M,D,p]$ in $[\mathcal{E}^*,\pi^*,M^*,D^*,p^*]$ if
$$p^*\circ\Psi = p \clubsuit $$
\end{defn}
\textbf{A COMMENT.}\\
Now let
$$ \mathsf{G} = \mathcal{A}ut([\mathcal{E},\pi,M,D])$$
be the group of automorphisms $\Psi\times \psi$ of a datum $[\mathcal{E},\pi,M,D]$. If $M$ is finite dimensional then $\mathsf{G}$ is a finite dimensional Lie group.\\
Of course the group $\mathsf{G}$ acts in the categories of structures of statistical models
$$\mathsf{M} = \mathcal{GM}(\Xi,\Omega)$$ 
and
$$ \mathbb{M} = [\mathcal{E},\pi,M,D,p]$$
\begin{defn} The moduli space of structures of statistical models in $[\mathcal{E},\pi,M,D]$ is the quotient datum
$$\mathsf{m} = \frac{\mathsf{M}}{\mathsf{G}} \clubsuit $$
\end{defn}
The challenge is the search of an invariant which encodes the orbits of $\mathsf{G}$ in $\mathsf{M}$. That is what is named the problem of moduli space.\\
The theory of statistical models has a homological nature (in the theory of KV homology) \cite{Nguiffo Boyom(6)}.
\newpage
\bibliographystyle{amsplain}

\begin{thebibliography}{ABCD}
\bibitem[Amari]{Amari} Differential Geometry Methods in Statistics. Lecture Notes in Statistics 28 Springer Verlag, NY 1990
\bibitem[Amari-Armstrong]{Amari-Armstrong} Amari S-I and Armstrong J. Curvature of Hessian manifolds Diff geom Appl 33 (2014) 1-12
\bibitem[Amari-Nagaoka]{Amari-Nagaoka} Methods of information geometry. translations of Mathematical Monographs, AMS-OXFORD vol 191.
\bibitem[Armstrong-Amari]{Armstrong-Amari} J. Armstrong and S. Amari The Pontryagin Forms of Hessian Manifolds, Geometric Sicience of information,242-247, Lecture Notes in Comput.sci.,9389 Springer, Cham, 2015
\bibitem[Arnaudon-Barbaresco]{Arnaudon-Barbaresco} Arnaudon M. and Barbaresco F. Medeans and means in Riemannian geometry: Existence, uniqueness and computation in Matrix Information Geometry; Springer: Heidelberg,Germany 2013 169-197

\bibitem[Arnaudon-Nielsen]{Arnaudon-Nielsen} Arnaudon M. and Nielsen F. Medians and means in Fisher geometry. LMS J. Comput. Math 2012 15, 23-37.
\bibitem[Auslander-Kuranishi]{Auslander-Kuranishi} Auslander L. and Kuranishi M. On the holonomy group of locally euclidean space. Ann of Math 65 (1957) 411-425.
\bibitem[Auslander-Markus]{Auslander-Markus} Auslander L. and Markus L. Holonomy of flat affinely connected manifolds. Ann of Math  62 (1955) 139-152.
\bibitem{Ay-Tuschmann} Ay,N. and Tuschmann, W. Dually flat manifolds and globalinformation geometry; Open. Syst. Inf. Dyn. 9(2),195-200 (2002)
\bibitem[Barndorff-Nielsen]{Barndorff-Nielsen} Barndorff-Nielsen O.E. Information and exponential families in statistical theory.  Wiley, New York]
\bibitem[Barbaresco]{Barbaresco}in GSI2013: Entropy Special Issue (2014). Barbaresco-Nielsen Editors
\bibitem[Baudot-Bennequin]{Baudot-Bennequin} Baudot P. and Bennequin D. The homological nature of Entropy, maxEnt 2014: (1) Proc Amer Institute phys (2014),(2) Entropy/Special Issue 2015. Barbaresco-Nielsen Editors.
\bibitem[Barbaresco]{Barbaresco} Barbaresco F. Koszol information geometry and Souriau geometric temperature/capacity of Lie group themodynamics. Entrpoy 16, (2014) 4521-4565
\bibitem[Benoit]{Benoit} Benoit Y. une vari\'et\'e non affine. Jour Differential Geometry (1992)
\bibitem[Benson-Gordon]{Benson-Gordon} Benson S. and Gordon C.S. Kaehler and symplectic structures on nilmanifolds, Topology 27 (1988) 513-518
\bibitem[Carriere]{Carriere} Autour de la conjecture de L. Markus sur les vari'et'es affines. Inv. Math. 95 (1989) 615-628
\bibitem[BCGRS]{BCGRS} P. Biliavski, M. Cahen, S. Gutt, J. Rawnsky and L. Schwachhofer. Special symplectic connections arXiv[math. SG] May 2006
\bibitem[Byande]{Byande} Byande P.M Des structures affines \`a la g\'eom\'etrie de l'information. La notion de T-plongement. Ed. Omniscriptum (2011)
\bibitem[Cahen-Schwachhofer]{Cahen-Schwachhofer}, M. Cahen and L.J Schwachhofer. Special symplectic connections. arXiv:[math DG] Sept 11 (2009).
\bibitem[Cartan(0)]{Cartan(0)} Cartan E. La Theory des groups finis et contnuset lAnalysis Situs. M\'emorial des Sciences Mathematiques, vol. 42 (Oeuvres compl\`etes vol I 1165-1225.)
\bibitem[Cartan(1)]{Cartan(1)}. Cartan E. (1) Sur les domaines born\'es de l'espace de n variabes complexes, Abh. Math. Sem. Hambourg 11, 116-162
\bibitem[Cartan(2)]{Cartan(2)} Cartan E. Sur les vari\'et\'es \`a connexions affine et la t\'eorie de la relativit'e g\'en\'erale. Ann Sci Ec Norm Sup,(1): 40 (1923) 325-412. (2): 41 (1924) 1-25. 42 (1925) 17-88.
\bibitem[Cartan(3)]{Cartan(3)}. Cartan E. Les espaces \`a connexions conformes Ann Soc Pol Math 2 (1923) 171-221.
\bibitem[Cartan(4)]{Cartan(4)}.Cartan E. Sur les vari\'et\'e \`a connexions projectives. Bull Soc Math France 52 (1924) 205-244.
\bibitem[Chevalley]{Chevalley}. Chevelley C. Theory of Lie groups, I. Princeton 1946, reprinted 1999.
\bibitem[Dorfmeister]{Dorfmeister} Dorfmeister J. Homogeneousn Kaehlerian manifolds admitting a  solvable transitive group of automorphisms. Ann Sci Ecole Norm Sup (4) 18 (1985) 143-180.
\bibitem[Dorfmeister-Nakajima]{Dorfmeister-Nakajima} Dorfmeister J. and  Nakajima K. The fundamental conjecture for homogeneous Kaehlerian manifolds. Acta Math. 161 (1988) 23-70
\bibitem[FGH]{FGH} Fried D., Goldman W. and Hirsch M. Affine manifold with nilpotent holonomy Comment Helv Math 56 (1983) 487-523
\bibitem[FUCK]{FUCK} Fucks D.B. Cohomology of infinite dimensional Lie algebras. Contemp. Soviet Mathematics. Consultant Bureau Publishing, New York 1986.
\bibitem[Gerstenhaber]{Gerstenhaber} Gerstenhaber M. Deformations of Rings and Algebras, Ann of Math. 79 (1964) 59-103
\bibitem[GPSV]{GPSV} Gindikin S.G, Pyateckii I.I and Vinberg E. B. Homogeneous Kaehler Manifolds, Geometry of bounded domains, CIMEIII, ciclo 3, ciclo Urbino, Ed Cremonese Roma (1968) 3-87.
\bibitem[Goldman]{Goldman} Goldman W.M. Complete affine manifolds: a survey. WWW.math.umd.edu/WMG/MFO.
\bibitem[Goldschmidt]{Goldschmidt} On the nonlinear Spencer cohomology of Lie Equation. I-V J. Differential Geometry 13 (1978); 16 (1981).
\bibitem[Gromov]{Gromov} Gromov M. The searc of structures: (1) ECM6, Krakow 2012. (2) maxEnt2014, Proc Amer Inst. Phys (2013).
\bibitem[Gromov 2000]{Gromov 2000} Gromov M. Geometric structures, peopole.mpim.mpg.de/hwbllmnn/arciv/geostr00.pdf
\bibitem[Guillemin]{Guillemin} Guillemin V. The integrability problem for G-structures, Trans Amer Math Soc 116 (1965) 544-560
\bibitem[Guillemin-Sternberg]{Guillemin-Sternberg} An algebraic Model for Transitive Differential Geometry. Bull of Amer Math Soc 70 (1964) 16-47
\bibitem{Guts} Guts A. K. Private communication (1986!)
\bibitem[Hano-Morimoto]{Hano-Morimoto} Hano and Morimoto. Note on the group of an affine transformation of affinely flat manifolds. Nagoya Math Journal 8 (1955) 71-81.
\em[Hochschild]{Hochschild}. G. The Structures of Lie Groups, Holden-day, San Francisco 1965.
\bibitem[Kuranishi]{Kuranishi} Kuranishi M. New prof for the existence of locally complete families of complex structures, Proc. Conference on Complex Analysis, Minneapolis (1964), Springer-Verlag (1965), 142-154.
\bibitem[Kobayashi]{Kobayashi} Kobayashi K. Differential Geometry of complex vector bundles, Publ Math Soc Japan, no 15, Iwanami-Princton Univ Press
\bibitem[Kodaira]{Kodaira} Kodaira K. Complex manifolds and deformation of complex structures,Springer-Verlag.
\bibitem[Katsumi]{Katsumi} Katsumi Y. On Hessian structure on an affine manifold, in Manifolds and Lie groups in Honor of Yozo Matsushima, Progress in mathematics vol 114, 449-459
\bibitem[Kaup]{Kaup}.hyperbolische komplexe Rume. Ann. Institut Fourier 18 (1968) 303-330.
\bibitem[Kim]{Kim} Kim H. Complete left invariant affine structures on nilpotent lie groups, J. Differential Geom 24 (1986) 373-394
\bibitem[Kobayaschi]{Kobayaschi} Kobayaschi S. The Theory of connections. Annali di Matematica 43 (1957) 119-194.
\bibitem{Kontsevich} Kontsevich M. Deformation quantization of Poisson manifolds, Lett. in Phys.66 (2003) 157-216
\bibitem[Koszul(1)]{Koszul(1)} Koszul J-L. D\'eformations de connexions localement plates, Ann Institut Fourier, 18, 1 (1968), 103-114
\bibitem[Koszul(2)]{Koszul(2)} Koszul J-L. Domaines homog\`enes born\'es et orbites des transformations affines, Bull Soc. Math. France 89 (1961) 515-533.
\bibitem[Koszul(3)]{KOSZul(3} Koszul J-L  Sous-groupes discrets des groupes de transformations affines admettant trajectoire ouverte, C.R. Acad. Sc. Paris 259 (1964) 3675-3677
\bibitem[Koszul(4)]{Koszul(4)} J-L Koszul. Vari\'et\'es localement plates et convexit\'e, Osaka J. Math 2 285-290
\bibitem[Koszul(5)]{Koszul(5)} Koszul J-L, Sur les j-algebres (unpublished)
\bibitem[Kumpera-Spencer]{Kumpera-Spencer} Kumpera A. and Spencer D.C. Lie equations, the general theory. Annal of Math. studies, Princeton Univ Press.
\bibitem[KU-RO]{KU-RO} Kuranishi M. and Rodrigues A. Quotients of Lie pseudogroups by invariant fiberings. Nagoya Math J. 24 (1964) 109-128
\bibitem[McCleary]{McCleary} A user's guide to spectral sequences. Cambridge University Press (2001)
\bibitem[Matsushima]{Matsushima} Matsushma Y. Affine structures on complex manifolds, Osaka J. Math 5 (1968) 215-22
\bibitem[McCullagh]{McCullagh} McCullagh P. What is a statistical model? The Annals of Statistics, (2002) Vol(30) NO5 1225-1310.
\bibitem[McDuff]{McDuff} McDuff D. The moment map for circle actionson symplectic manifolds. J. Geome. Phys. 5 (1988) 149-160
\bibitem[Malgrange]{Malgrange} Malgrange B. Equations de Lie I,II Jour of Diff Geometry 6 (1972) 503-522 and 7 (1972) 117-141
\bibitem[Medina-Saldarriaga-Giraldo]{Medina-Saldarriaga-Giraldo}, Medina A;,Saldarriaga O. and Giraldo H. Flat affine or projective geometries on Lie groups. jour of Algebra, 455 (2016) 183-208.
\bibitem[Medina-Revoy]{Medina-Revoy}. Algebres de Lie et produits scalaires invariants, Ann Sci EC Nor Sup 18 (1985) 553-561.
\bibitem[Medina]{Medina}. A. Medina Groupes de Lie munis de metriques bi-invariantes. Tohoku Math. journal 37 1984 405-421.
\bibitem[Milnor]{Milnor} On the fundamental groups of complete affinely flat manifolds, Advances in Math.25 (1977), 178-187
\bibitem[Moerdijk-Mrcun]{Moerdijk-Mrcun} Introduction to Foliations and Lie groupoids, Cambridge studies in advanced mathematics 91
\bibitem[Molino]{Molino} Riemannian foliations; Birkhauser, Boston Massachusetts.
\bibitem[Murray-Rice]{Murray-Rice} Differential geometry and Statistics. Monographs on Statistics and Applied Probability 48, Chapman \& hall/CRC.
\bibitem[Nakajima]{Nakajima}, On j-algebra and homogeneous Kaehler manifolds, Hokkaido Math Journal vol 15 N°1 (1986) 21-47
\bibitem[Nguiffo Boyom(1)]{Nguiffo Boyom(1)} Nguiffo Boyom M.(1) Structures affine homotopes \`a z\`ero des groupes de Lie nilpotents. Jour of Diff Geometry 31 (1990) 859-911
\bibitem[Nguiffo Boyom(2)]{Nguiffo Boyom(2)} Nguiffo Boyom M.(2) Sructures isotopes \`a z\'ero . Ann. della Scu. Norm Sup di Pisa. IV (1993) 91-131.
\bibitem[Nguiffo Boyom(3)]{Nguiffo Boyom(3)} Nguiffo Boyom M. The cohomology of Koszul-Vinberg algebras. Pacific Journal of Mathematics vol 225 N°1 (2006) 119-153
\bibitem[Nguiffo Boyom(4)]{Nguiffo Boyom(4)} Nguiffo Boyom M. The homology of Koszul-Vinberg algebroids and Poisson manifolds I, Banach Center Publ (2001) 99-110
\bibitem[Nguiffo Boyom- Wolak(1)]{Nguiffo Boyom-Wlak(1)} Nguiffo Boyom M. and Wolak R. Affine structure and KV-cohomology, Jour of Geom and Phys 42 (2002)307-317.
\bibitem[Nguiffo Boyom-Wolak(2)]{Nguiffo Boyom-Wolak(2)} Nguiffo Boyom M. and Wolak R. Local structures of Koszul-Vinberg algebroids, Bull Sci. Math. 128 (2004)467-479
\bibitem[Nguiffo Boyom-Wolak(3)]{Nguiffo Boyom-Wolak(3)} Trnsversely Hessian Foliations and Information geometry, Intern Jour of Math (IJM) vol 27 N° 11 (2016) 
\bibitem[Nguiffo Boyom-Byande-Ngakeu-Wolak]{Nguiffo Boyom-Byande-Ngakeu-Wolak}, Nguiffo Boyom M;, Byande P.M, Ngakeu F. and Wolak R. KV cohomology and differential Geometry of locally flat manifolds, Afr Diaspora J. Math (2012) vol 14 197-226
\bibitem[Nguiffo Boyom(6)]{Nguiffo Boyom(6)} Nguiffo Boyom. Foliations-Webs, Hessian Geometry, Information Geometry, Entropy and Cohomology, Entropy 2016, vol 18,N°12,433.
\bibitem[Nguiffo Boyom(7)]{Nguiffo Boyom(7)} R\'eductions Kaehl\'eriennes dans les goupes de lie resolubles et applications. Osaka J. Math vol 47 (2010) 237-283
\bibitem[Nguiffo Boyom(8)]{Nguiffo Boyom(8)} Nguiffo Boyom M. Alg\`ebres sym\'etriques \`a gauche et alg\`ebres de Lie r\'eductives, PhD Fac Sciences University of Grenoble (1968).
\bibitem[Nguiffo Boyom(9)]{Nguiffo Boyom(9)} Affine embeddings of real Lie groups, Proc London Math Soc LNS 26 (1976) 21-39.
\bibitem[Nijenhuis]{Nijenhuis} Nijenhuis A. Sur une classe de proprietes communes à divers types d'algebres, Enseign Mathem 1968,14, 225-277
\bibitem[Nijenhuis-Richardson]{Nijenhuis-Richardson}. Nijenhuis A. and Richardson  W. cohomology and deformations of algebraic structures, Bull Amer Math Soc 70 (1964) 1-29 
\bibitem[Nirenberg]{Nirenberg} Nirenberg L: lectures on Linear partial differential equations  conf. Board of Math Sci no  N°17 Amer Math. Soc 1973.
\bibitem[Palais]{Palais} Palais R. A global formulation of the Lie theory of Transformation Groups, Mem. Amer Math Soc N°22 (1957)
\bibitem[Pennec]{Pennec} Pennec X. GSI2013. Geometric statistics on manifolds and Lie groups. LNCS 8085 Nielsen- Barbaresco (eds) 59-67. ( See also Entropy Special Issue (2014). Barbaresco-Nielsen eds.)
\bibitem[Petrie-Handal]{Petrie-Handal} Petrie T. and Handal J, Connections, Definite forms and Four manifolds, Oxford Mathematical Monograph, Oxford Science Publications 1990
\bibitem[Piper]{Piper} Piper S. Algebraic deformation theory. J. Differential Geometry, 1 (1967) 133-168
\bibitem[Pontryagin]{Pontryagin} Pontryagin. Topological Groups.
\bibitem[Raghunathan]{Raghunathan} Raghunathan S. The deformations of linear connections XXX
\bibitem[Pyatetskii-Shapiro]{Pyatetskii-Shapiro} Pyatetskii-Shapiro I.I. On a problem of E. Cartan, Dokl Akad Nauk SSSR, 124 (1959) 272-273
\bibitem[Reinhardt]{Reinhardt} B.L. Reinhardt. Foliated manifolds with bundle-like metrics. Ann. Math. (2) 69 (1959) 119-132
\bibitem[Shima(1)]{Shima(1)} Shima H. Homogeneous hessian manifolds, in honnor of Yozo Matsushima, Progress in Mathematics, 14, Birkhauser, Boston 1981, 385-392
\bibitem[Shima(2)]{Shima(2)} The Geometry of Hessian Manifolds, World Scientific Publishing Co, NJ 2007.
\bibitem[Singer-Sternberg]{Singer-Sternberg} Singer I. and Stenberg S. The infinite groups of Lie and Cartan, Journal An.Math.(15) 1965, 1-114.
\bibitem[Smilie]{Smilie} Smilie J. An obstruction to the Existence of affine structures, Inv. Math 64 (1981) 411-415
\bibitem[Spencer]{Spencer} Spencer D.C. Overdeterminated systems of partial differential equations. Bull of the Amer Math Soc  75 (1964) 179-223
\bibitem[Vey(1)]{Vey(1)} Vey J. Deformations du crochet de Poisson sur une variété symplectique, Comment Math helv 50 (1975) 421-454
\bibitem[Vey(2)]{Vey(2)} Vey J. Sur la division des domaines de Siegel, Ann Scient Ec Norm Sup 43 (1970) 479-506
\bibitem[Vinberg-Katz]{Vinberg-Katz} Vinberg E.B. and Katz V. Kvaziodnorodnye konusy, Mat. Zametki 1 (1967)347-354
\bibitem[Vinberg(1)]{Vinberg(1)} Vinberg E.B. Convex homogeneous cones,Transl Moscow Math. Soc 12 (1963)340-403
\bibitem[Vinberg(2)]{Vinberg(2)} Vinberg . B. (1963) The theory of homogeneous cone. Trans. Moscow. Math Soc 340-403
\bibitem[Vinberg-Gindikin-Pyatetskii-Shapiro]{Vinberg-Gindikin-Pyatetskii-Shapiro} E.B Vinberg, SG Gindikin and II Pyatetskii-Shapiro. On the classification of and canonical realization of complex homogeneous bounded domains Proc Moscow Math Soc 12. 404-437
\bibitem[Vinberg-Gindikin-Pyatetskii-Shapiro]{Vinberg-Gindikin-Pyatetskii-Shapiro}. Vinberg E.B, Gindikin S.G and Pyatetskii-Shapiro I.I: Homogeneous Kaehler manifolds, Collect CIME.
\bibitem[Wolf]{Wolf} Wolf J. Spaces of constant curvature, Publish of Perish, Boston Mass; (1974).
\bibitem[Zeghib]{Zeghib} A. Zeghib On Gromov theory of rigid transformation groups, a dual approach, Ergod. Theory Dynam. Syst. 2000,20, 935-946.
\end{thebibliography}

\end{document}